\documentclass[11pt]{article}
\usepackage{cite}
\usepackage{amsmath,amsfonts,amssymb,booktabs,tabularx}
\pdfoutput=1
\usepackage[small,bf,hang]{caption}
\usepackage{slashed}
\input epsf.sty
\usepackage{epsfig}
\usepackage[titletoc,toc]{appendix}
\usepackage{rotating,diagbox}
\usepackage{xcolor}

\def\hybrid{
        \topmargin -20pt
        \oddsidemargin 0pt
        \headheight 0pt \headsep 0pt
        \textwidth 6.55in 
        \textheight 9.5in 
        \marginparwidth .875in
        \parskip 5pt plus 1pt \jot = 1.5ex}

\hybrid

 \csname
@addtoreset\endcsname{equation}{section}

\newcommand{\p}{\partial}

\def\moth{\mathsurround=0pt}
\newdimen\zo \zo=0pt

\def\tick{\leaders\hrule height 0.5ex depth 0pt \hskip 0.5pt}
\def\upboxfill{$\moth \setbox\zo\hbox{\tick}%
  \hskip 3pt\hbox to 0pt{$\tick$\hss}\hrulefill \hbox to 7.5pt{$\tick$\hss}$}

\def\dtick{\leaders\hrule height .34pt depth 0.5ex \hskip 0.5pt}
\def\downboxfill{$\moth \setbox\zo\hbox{\dtick}%
  \hskip 2pt\hbox to 0pt{$\dtick$\hss}\hrulefill \hbox to 2pt{$\dtick$\hss}$}

\def\bec{\begin{center}}
\def\ec{\end{center}}

 \def\det{{\rm det\,}}
\def\be{\begin{equation}}
\def\ee{\end{equation}}
\def\bea{\begin{eqnarray}}
\def\eea{\end{eqnarray}}
\def\ba{\begin{array}}
\def\ea{\end{array}}

    \newcommand{\sbkt}[1]{\left[#1\right]}
    \newcommand{\bkt}[1]{\left(#1\right)}
	\newcommand{\PreserveBackslash}[1]{\let\temp=\\#1\let\\=\temp}
\newcolumntype{C}[1]{>{\PreserveBackslash\centering}p{#1}}
\newcolumntype{R}[1]{>{\PreserveBackslash\raggedleft}p{#1}}
\newcolumntype{L}[1]{>{\PreserveBackslash\raggedright}p{#1}}
\begin{document}

\begin{titlepage}

\rightline{\tt MIT-CTP-5105} 
\rightline{\tt UUITP-8/19}
\hfill \today
\begin{center}
\vskip 0.5cm

{\Large \bf {Extremal isosystolic metrics with multiple\\  
\vskip 0.5cm  
bands of crossing geodesics}}

\vskip 2.5cm

 {\large {Usman Naseer$^{1,2,3}$ and Barton Zwiebach$^3$}}

{\em  \hskip -.1truecm 
${}^1$Department of Physics, \\
Harvard University,\\
Cambridge, MA, 02138, USA\\
}
\vskip 10pt
{\em  \hskip -.1truecm 
${}^2$
Department of Physics and Astronomy,\\
     Uppsala University,\\
     Box 516,
     SE-751s 20 Uppsala,
     Sweden\\
}
\vskip 10pt
{\em 
${}^3$
Center for Theoretical Physics \\
Massachusetts Institute of Technology\\
Cambridge MA 02139, USA\\
\tt unaseer@g.harvard.edu, zwiebach@mit.edu \vskip 5pt }

\vskip 2.5cm
{\bf Abstract}

\end{center}

\noindent
\begin{narrower}

\baselineskip15pt

We apply recently developed convex programs to 
find the minimal-area Riemannian metric on $2n$-sided
polygons ($n\geq 3$)
with length conditions
on curves joining opposite sides.  
We argue that the Riemannian extremal metric
coincides with the conformal
extremal metric on the regular $2n$-gon. 
The hexagon was considered by Calabi.  The region covered by the maximal number $n$ 
of geodesics bands extends over most of the surface and exhibits positive
curvature.   As $n\to \infty$ the metric, away from the boundary, 
approaches the well-known round extremal metric on $\mathbb{RP}_2$. 
We extend Calabi's isosystolic variational principle
to the case of regions with more than three bands of systolic geodesics.  
The extremal metric on $\mathbb{RP}_2$ is a stationary point of this functional
 applied to a surface with infinite number of systolic bands.

\end{narrower}

\end{titlepage}

\baselineskip11pt

\tableofcontents

\baselineskip15pt

\section{Introduction}

The subject of minimal area {\em conformal} metrics on Riemann surfaces
is an old one (see, for example,~\cite{ahlfors, strebel}).
 Sometimes called extremal length problems, they 
ask for the conformal metric of least area on a fixed Riemann
surface under some length conditions.  Typically the length conditions
require that certain sets of curves have length bounded below
by some fixed constants.  In order to find the string diagrams that would
define a closed string field theory a minimal area problem on (punctured)
Riemann surfaces was proposed in which a single length condition applies 
to {\em all} non-contractible closed curves on the 
surface~\cite{Zwiebach:1990nh,Zwiebach:1992ie,Wolf:1992bk}.  
All such curves must be
longer than the {\em systole} $\ell_s$, a number that can be conventionally
set equal to one, or $2\pi$, or any fixed constant.  

Consider this conformal minimal area problem for the case of a Riemann surface of genus $g$ with no punctures.  In general the extremal metric is not known. 
In genus $g$ we have a moduli space ${\cal M}_g$ of Riemann surfaces.  
For string field theory we need the minimal area metric for each point on
 ${\cal M}_g$, that
is, for each and every Riemann surface of genus $g$.  If we knew these extremal
metrics and their area we could ask for those Riemann surfaces for which the
extremal area is a minimum over the moduli space.  Those surfaces would be
solutions of the problem of finding the {\em Riemannian} metric of least area on a
genus $g$ surface under the same length condition.  Indeed, as we vary over the
space of Riemannian metrics we change the complex structure, and in this way
we end up finding the Riemann surface (or surfaces) that has the least
area.  This problem of extremal Riemannian metrics has been investigated
in Gromov~\cite{gromov}  and is a problem of 
systolic geometry~\cite{m_katz,mberger,guth}.  
The conformal
specialization has also been explored in~\cite{bavard}.  The case of metrics
with non-positive curvature has been studied in~\cite{katz, katz2}. 

It is clear that on two-dimensional surfaces  
the conformal minimal area problem is richer than the Riemannian version.
In fact, extremal conformal metrics show different behavior.  
While conformal minimal-area metrics can have regions 
where there is exactly one band of 
systolic geodesics~\cite{strebel,Zwiebach:1992ie},  
Riemannian extremal metrics cannot, as explained
by Calabi~\cite{calabi}.  In Riemannian extremal metrics regions with exactly two bands
of geodesics are flat and the geodesics are orthogonal~\cite{calabi}. 
 In the conformal version, as was recently discovered~\cite{Headrick:2018ncs,Headrick:2018dlw}, regions with two bands of geodesics
support positive and sometimes negative curvature.  

Indeed,  recent work with 
M.~Headrick~\cite{Headrick:2018ncs,Headrick:2018dlw}
applied the methods of convex optimization~\cite{boyd} to the conformal
minimal area problem.  The result was a couple of programs, a primal
and its dual, that give new analytic insights into the form of the
 minimal area metrics.  Additionally, they
can be used to find accurate minimal area metrics and
to deal with the minimal area problem of string field theory.  
The dual program involves the maximization
of a functional. Its optimum, by the property of strong duality,
 coincides with the minimum of the primal.   
 
In this paper we explore a set of issues 
motivated by the work of Calabi~\cite{calabi}.
He considered a regular hexagon in $\mathbb{R}^2$ 
and asked for the Riemannian metric of least
area under the condition that all curves joining opposite sides of the hexagon
be longer than or equal to one.  
Reference~\cite{calabi} assumes that the metric
is unique and thus admits the dihedral group $D_6$ of symmetries as isometries.
We will consider this problem in detail, as
well as  the obvious generalization to regular $2n$-sided polygons $P_{2n}$, with $n\geq 3$.   We visualize these regular 
polygons as regions of the complex $z$ plane. 
We claim that the minimal area {\em Riemannian} metric on $P_{2n}$,
 if unique,  is in fact the minimal area {\em conformal} metric on $P_{2n}$.   
The argument is
 presented in section~\ref{riem_conf}. 
It follows that we can use the methods of~\cite{Headrick:2018ncs,Headrick:2018dlw} 
to study the problem posed by Calabi 
as well as    
 its generalization to higher polygons.    
 
 Calabi  also proposed
 a variational principle for the metric, valid for regions
 with exactly three bands of systolic geodesics.  This  
 was studied further in~\cite{bryant},
who gave more details on the differential equations arising from
the variational principle and following~\cite{calabi} proposed an extension to 
regions covered by more than three bands of 
systolic geodesics.  We have found difficulties with this proposal
and have constructed a new variational principle to circumvent these
difficulties.

For the case of the hexagon we confirm that there is a large central 
region $U_3$ covered by
three bands of geodesics, and a region $U_2$ comprising 
neighborhoods of the vertices that  are flat
and are covered
by two orthogonal bands 
of geodesics (see Figure~\ref{fig:hexfinal}). A Gauss-Bonnet argument shows that
the integral of the Gaussian curvature over $U_3$ is equal to $2\pi$, the amount
that corresponds to half a sphere.    
We find an accurate
estimate of the extremal area using the primal and dual programs.\footnote{It is not clear
from~\cite{calabi} if the differential 
equations relevant to the problem were studied numerically.} 
 Moreover, the methods of~\cite{Headrick:2018ncs} give an exact relation ${\cal P} = 4 A/\ell_s$  
between the perimeter ${\cal P}$ 
of the polygon and its area $A$  in the minimal-area metric with systole $\ell_s$.

Since  the possible
behavior of metrics in regions covered by multiple bands of systolic geodesics
is of great interest,  higher polygons are a natural ground for exploration.  
A polygon $P_{2n}$ 
is expected to have a region $U_n$ with precisely $n$ systolic geodesics
going through every point.    For the case of an 
octagon we find a rather large central region 
$U_4\in P_8$ (see Figure~\ref{fig:octfinal}).   There are also small regions $U_3$ and
$U_2$.  We find strong evidence that the integral of the Gaussian curvature over $U_4$ is equal to $2\pi$.   
We briefly discuss the 10-gon or decagon. For this polygon we find strong evidence 
that the integral of the Gaussian 
curvature over $U_5$ is equal to $2\pi$, 
in contradiction with the
claim~\cite{bryant} that the extremal metric in regions $U_k$ with $k \geq 5$
must be flat (see Figure~\ref{fig:10Gon}). 

We observe that the region $U_n\subset P_{2n}$ covered by $n$ bands
fills a  larger part of $P_{2n}$ 
as $n$ increases, and 
 the integrated curvature over $U_{n}$ remains constant at $2\pi$.   
 These patterns make 
sense in the limit of very large $n$.   
We find evidence 
 that, away from the boundary of the polygon, 
the minimal area metric on $P_{2n}$ 
approaches 
 the minimal area Riemannian metric on $\mathbb{RP}_2$ as $n\to \infty$.  This
metric was obtained by Pu\cite{pu}.  It takes the form of a constant-curvature hemisphere, whose boundary is subject to
antipodal identifications.  We can view this metric as one with infinite
number of geodesic bands, each band defined by a pair of antipodal
points on the boundary prior to identification.  
Each geodesic band covers the full hemisphere and
the integrated curvature is indeed $2\pi$.   
The conformal version of the minimal area metric on $\mathbb{RP}_2$ has
also been studied.
This problem asks for the minimal area conformal metric on a disk
$|z| \leq {1\over 2}$ with the condition that any curve going from a point
to its antipode is of length larger than or equal to one.  Reference~\cite{ahlfors}
shows that the extremal metric is indeed the constant curvature metric that makes
the disk into a hemisphere.  This is, of course, the same exact metric that Pu found,
showing that in this case the Riemannian and conformal problems are the same, 
just as we believe they are for all polygons $P_{2n}$.   

The extremal metric on $P_{2n}$ 
approaches the hemisphere metric except at the boundary, where the discrete
nature of the problem appears to make a difference.  The perimeter/area relation
${\cal P} = 4 A/\ell_s$, valid for all polygons $P_{2n}$, does not hold for the hemisphere
metric where we have instead ${\cal P} = \pi A/\ell_s$.   In fact, each corner of the polygon
has a flat neighborhood $U_2$ with the two incident edges meeting orthogonally.
This makes the metric diverge in the original coordinates where incident edges
meet with angles approaching $\pi$ as $n\to \infty$.  We have some sort of infinitely
`serrated' boundary.  More work will be needed to find the detailed picture of
the metric near the boundary.  We also do not know much about line curvature
singularities of the metric as well as the curvature in the rapidly shrinking $U_{n-1}, \ldots , U_2$ regions.

We then examine the dual functional of~\cite{Headrick:2018ncs} 
for the extremal metric on $\mathbb{RP}_2$.  Here the challenge
is to show that this metric provides the optimum of the functional.  
Each band, parametrized by the value $\phi_0$ of the azimuth of the 
starting point (and azimuth $\pi + \phi_0$ for the ending point) provides a local coordinate system with coordinates $x_{\phi_0}$ measuring length along the 
geodesics and $\varphi_{\phi_0}$ which is a constant along each geodesic.  
Since the geodesics are known, $\varphi_{\phi_0}$
 is determined up to reparameterization.   
There is one key constraint following from the
dual functional.  The equation
\be
\sum_{\phi_0}  |d\varphi_{\phi_0}| = 1 \,,
\ee
must hold at {\em every} point of the sphere.  This is not easily satisfied since 
every candidate $\varphi_{\phi_0}$ is a nontrivial function on the sphere.  Moreover, the sum must be turned into an integral.
We found the solution $\varphi_{\phi_0}(\theta, \phi)$ satisfying this
condition and showing that the extremal metric is the optimum of the dual program.

The metric on $\mathbb{RP}_2$ also allows us to test 
the variational principle of~\cite{calabi,bryant}.  We use
two bands of geodesics on the extremal metric to introduce coordinates 
$X,Y$ on the  sphere.   
A third,  arbitrary band, is fixed by a pair of antipodal points and defines 
a length-based coordinate $Z (X,Y)$.  
The least-area
variational principle in~\cite{calabi,bryant} gives a partial differential equation
for $Z$.  Unfortunately, the equation is 
\emph{not} satisfied.  We propose a modification of Calabi's variational
principle for regions with three or more bands of geodesics.  The 
new proposal uses Lagrange multipliers to impose the constraints 
on the norm of one-forms that
arise from each geodesic band.  The equations of motion are
nontrivially modified.   We then show that the extremal metric on $\mathbb{RP}_2$
is a solution for the equations following from the new variational principle.

Here is a summary of our main results:
\begin{enumerate}
\item  Showing (up to a uniqueness assumption) that the Riemannian
and conformal minimal area metrics are the same for polygons $P_{2n}$, $n\geq 3$, 
with dihedral symmetry (section~\ref{riem_conf}).  

\item  Finding an extension of Calabi's isosystolic variational principle 
valid for regions
$U_m$ covered by $m \geq 3$ bands of systolic geodesics (section~\ref{ext-var-principle}). 

\item  Showing that the extremal metric on $\mathbb{RP}_2$ is 
a solution of the dual conformal program (section~\ref{con_met_on_rp2}) as well as a solution
 of the new isosystolic
variational principle (section~\ref{ext-ire-met-new-var}).

\item A  study of extremal metrics on $P_{2n}$
with $n \geq 3$ showing the interplay of multiple bands of geodesics
and suggesting convergence for $n\to \infty$ to the 
extremal $\mathbb{RP}_2$ metric away from the
boundary (section~\ref{pus_metric}).  

\end{enumerate}

\medskip
This paper is organized as follows.  In section~\ref{the_for_gen_pol} we define
the minimal area problem on polygons and 
 establish that the Riemannian 
and conformal problems are equivalent.
In section~\ref{con-pro-for-con-met} we discuss the explicit formulation
of the primal and dual programs for arbitrary regular polygons $P_{2n}$.
The dual program has a height parameter $\nu$ which is fixed at the optimum.
This parameter determines the extremal area via the formula $A = n\nu \ell_s^2$. 
We show that this parameter also determines
the length $L_e$ of each edge of the polygon via the formula $L_e = 2\nu \ell_s$. 
This allows one to establish the perimeter/area relations discussed above.

Section~\ref{cal_bi_hex} gives our results for Calabi's hexagon.  We provide
a detailed accounting of the curvature on the surface, which includes
negative curvature on the boundary of $U_3$.  The critical
area $A=0.840$ is only about 3\% lower than the area of the admissible 
flat-metric on the hexagon.  The pattern of systolic geodesics on the surface
is shown in Figure~\ref{fig:hexfinal}.   Section~\ref{the_cas_of_the_oct} deals with
the case of the octagon and then briefly the decagon, discussing the curvature
on their central regions. 

In section~\ref{pus_metric} we review the extremal metric
on $\mathbb{RP}_2$~\cite{pu},  setting up the notation to compare it with
 the extremal metric on 
$P_{2n}$ for large $n$.  We show
numerical evidence of convergence of the metric to the 
$\mathbb{RP}_2$~ metric as long
as we are away from the boundary.  
 In section~\ref{con_met_on_rp2}
we show that the extremal metric on $\mathbb{RP}_2$ is in fact the optimum
for the variational system defining the dual program of~\cite{Headrick:2018ncs}.
This is a rather nontrivial test of dual program and provides further evidence
of the claimed relation of the polygonal problem   
to the $\mathbb{RP}_2$ problem in the limit $n\to \infty$.

In section~\ref{cal_var_pri} we review the variational
 approach of Calabi~\cite{calabi} explaining
the introduction of potentials associated to (two) systolic bands, that serve as coordinates on the surface, and the description of the area form using the 
potential for yet another systolic band.  We review the variational principle,
and write down the explicit form of the equation of motion.  We use
the $\mathbb{RP}_2$ metric to show that the equations of motion must
be changed to deal with more than three bands of systolic geodesics. 
In section~\ref{iso-var-prin-mult-fol} we give the isosytolic variational
principle valid in regions with three or more bands of systolic geodesics.
We confirm that the extremal metric on $\mathbb{RP}_2$ is a stationary
point of the new functional.

\section{Riemannian and conformal metrics with dihedral symmetry}\label{the_for_gen_pol}

In this section we begin by defining the regular $2n$-gon $P_{2n}$ 
as a region of the complex plane. After describing the
sets of curves whose lengths are constrained, we define conformal
and Riemannian versions of the minimal area problem.
   Assuming the Riemannian problem has a unique solution
we show that it coincides with the solution 
 of the conformal problem.

\subsection{The polygons and the curves to be constrained}

\begin{figure}[!ht]
\begin{center}
\epsfysize=8.0cm
\epsfbox{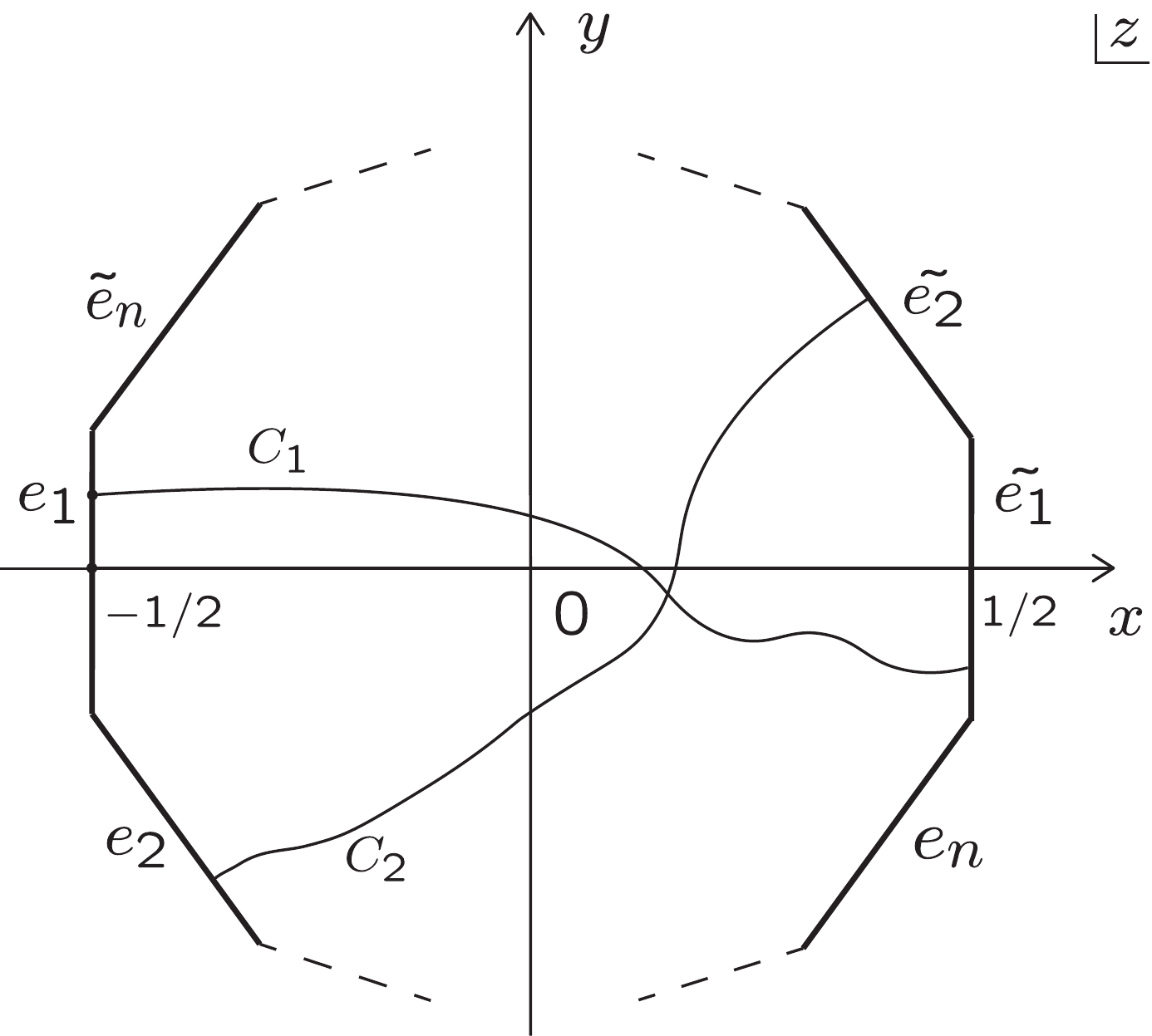}
\end{center}
\caption{\small The canonical presentation of the $2n$-gon $P_{2n}$ on the 
$z$ plane, with center at the origin and with opposite
edges $e_1$ and $\tilde e_1$ bisected orthogonally by the real line at
the points $x = \pm \tfrac{1}{2}$.  The curves starting on $e_1$ and ending 
on $\tilde e_1$ belong to the class $C_1$.  After a rotation of a curve on $C_1$
by an angle $\pi/n$ we obtain a curve in $C_2$, starting on $e_2$ and ending on
$\tilde e_2$. }
\label{fig:polycurves}
\end{figure}

We will consider regular $2n$-gons with $n\geq 3$ and will present them
as a region $P_{2n}$ on the $z = x+ iy$ plane with the center at the origin and 
with two opposite edges orthogonal to the $x$ axis at $x= \pm \tfrac{1}{2}$,
each bisected by the axis (see Figure~\ref{fig:polycurves}).  The  polygons so defined have 
apothem equal to $\tfrac{1}{2}$.
There are $n$ sets of curves, $C_1, \ldots C_n$, 
that are constrained in this polygon. These are the curves that join opposite edges of the polygon.  A typical curve starts
at any point on one edge and ends at an arbitrary point on the {\em opposite} edge. 
With the fiducial  distance between  edges 
 equal to one, it is convenient to set the systole length $\ell_s$ equal to one.
In the chosen presentation we take $C_1$ to be the curves that run 
 between the opposite edges that intersect the $x$ axis; the edge $e_1$ at $x=-\tfrac{1}{2}$
and the edge $\tilde e_1$ at $x = \tfrac{1}{2}$.  The curves on $C_2$
are obtained by a rotation of $\pi/n$ applied to the curves on $C_1$; they
begin on $e_2$ and end on $\tilde e_2$.  The curves on
$C_\alpha$, with $1< \alpha \leq n$ are obtained by a rotation of $(\alpha-1)\pi/n$ applied to the
curves on $C_1$.  They begin on $e_\alpha$ and end on $\tilde e_\alpha$.
The polygon $P_{2n}$ is invariant under transformations of the dihedral group
$D_{2n}$.  These are generated by
 rigid rotations by $\pi/n$ and reflections about any line
joining the midpoints of two opposite edges.

\begin{figure}[!ht]
\begin{center}
\epsfysize=7.0cm
\epsfbox{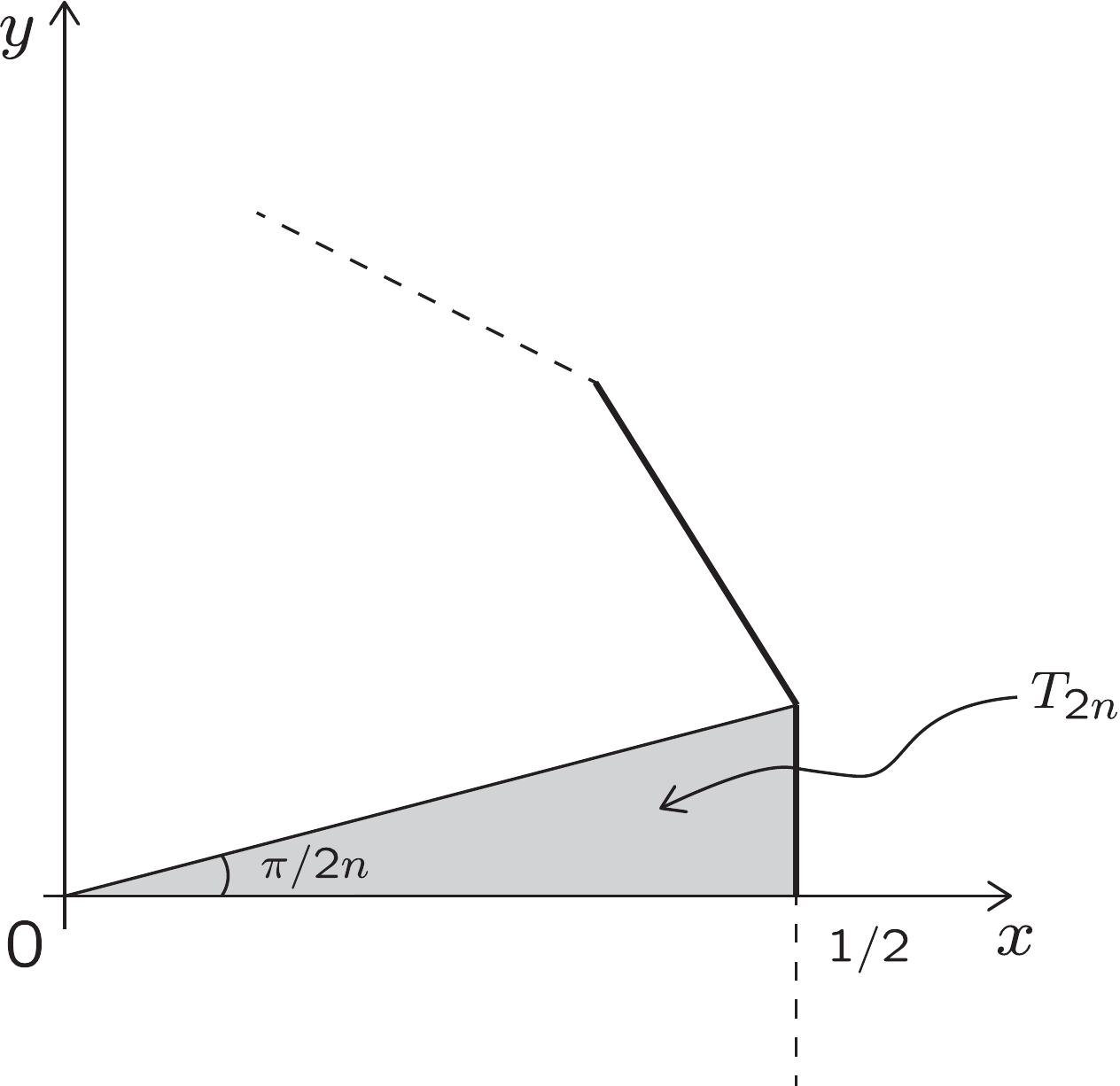}
\end{center}
\caption{\small The right triangle $T_{2n}$ in the polygon $P_{2n}$ is a 
fundamental domain for the extremal metric.  If known
on $T_{2n}$ the metric can be extended to $P_{2n}$ using the dihedral isometries.}
\label{fig:polycurves2}
\end{figure}

For the conformal minimal area problem we search in the space of metrics
conformal to the fiducial metric $|dz|^2$.    Because of convexity of the conformal
problem, the dihedral transformations
must be isometries of the extremal metric and we can therefore
restrict ourselves to search over the space of metrics on $P_{2n}$
with such symmetry.   This also means that we know the metric over
the whole polygon if we know it on the fundamental domain 
of 
the dihedral transformations.  We can choose the fundamental
domain to be the right triangle  $T_{2n}$ with a vertex at the origin, the second
vertex at $(\tfrac{1}{2}, 0)$ and the third vertex at the polygon corner above $x= \tfrac{1}{2}$ 
 (see Figure~\ref{fig:polycurves2}).

The fiducial metric $g^0_{\mu\nu}$ on the polygon
is taken to be the constant unit metric in the $z$ plane: 
\begin{equation}
 g^0_{\mu\nu} dx^\mu dx^\nu = dx^2+ dy^2\,. 
\end{equation}
It follows that the fiducial area form $\omega_0$ is
\be
\omega_0\ = \  d^2 x \sqrt{g^0}= dx dy\,.
 \ee
The general conformal metric 
on the polygon takes the form  
\be
g_{\mu\nu} = \Omega \, g^0_{\mu\nu}\,, 
\ee
with $\Omega$ a function on the polygon.  The area form $\omega$ is then
\be
\omega \ = \ d^2 x \sqrt{g} \ = \Omega \, dx dy \,. 
\ee
The area of $P_{2n}$ is $\int_{P_{2n}} dx dy\,  \Omega \,. $
Parameterizing a closed curve $\gamma$ using $x^\mu (t): [0,1]\to P_{2n}$,  we have
$\text{length} (\gamma) = \int_\gamma  \sqrt{\Omega}  |\dot x |_{{}_0}$, where 
$| \cdot |_{{}_0}$ denotes norm in the fiducial metric.

In the Riemannian version of the problem not much is changed.  The polygon
$P_{2n}$ is presented the exactly the same way, and the set of curves to
be constrained are exactly the same.  The metrics, however,
are no longer required to be conformal to the fiducial metric; they are arbitrary.
Lacking a convexity property 
 it seems challenging to show that
the extremal metric is unique, and therefore, just as in~\cite{calabi}, we will assume it is.
In this case, the extremal metric must be invariant under the dihedral group $D_{2n}$.

\subsection{Riemannian and conformal metrics on $P_{2n}$} \label{riem_conf}

Consider the Riemannian minimal area problem on the regular polygon
$P_{2n}$, presented as a region of the $z$ plane.  We will assume
that the extremal metric $g^0$ is unique and therefore has the dihedral symmetry $D_{2n}$  of the polygon as isometries.  
We want to show that $g^0$
defines the canonical conformal
structure on $P_{2n}$ or, equivalently, $g^0$
is conformal to the flat metric $|dz|^2$ on $P_{2n}$. 

Here are the steps that allow us to reach this conclusion.\footnote{We thank M.~Headrick for his help in constructing this argument.}
\begin{enumerate}

\item The conformal structures for a $2n$-gon 
can be represented as a round disk $|w|\leq 1$ 
with labeled marked points on the
boundary corresponding to  the vertices.  Alternatively,
they  can be represented
as the upper half plane $\hbox{Im}\, w > 0$,  with boundary the real line, and 
marked points on the boundary.  Since the positions of 
three marked 
points can be fixed at will by conformal maps, for a $2n$-gon, the conformal
moduli space ${\cal M}_{2n}$ is of real dimension $2n-3$. 

The moduli space can also be defined as the space of Riemannian metrics
on $P_{2n}$ modulo diffeomorphisms (keeping the boundary and vertices
 unchanged) and Weyl transformations.   Any metric $g$
on $P_{2n}$ defines a conformal structure  $\lambda(g) \in {\cal M}_{2n}$.  

\item 
Now consider the action of the dihedral group on any metric $g$ on $P_{2n}$
and focus on the rotation subgroup, generated by 
the transformation $Q:  z\to  z e^{i\pi/n}$ satisfying $Q^{2n}=1$.  
Note that $Q$ cycles the vertices
and $2n$ is the smallest
possible integer for which $Q$ raised to this power gives the identity.  
 We can think of $Q$ as an active transformation
that rotates the metric by $\pi/n$
 while keeping the labeling of the punctures the same.  
We write $Q: g \to Q(g)$.  It follows that this induces an 
action on the moduli space ${\cal M}_{2n}$ as $Q : \lambda(g) \to \lambda (Q(g))$.  
If the metric $g$ is dihedral invariant, it is also rotational invariant 
and $\lambda(g)$ is a fixed point of the action of $Q$ on the moduli space. 

Any action of $Q$ on the moduli space that leaves a point $\lambda \in {\cal M}_{2n}$
invariant is in fact
leaving a Riemann surface invariant and thus must be realized by a conformal
map, in this case, a projective transformation $\tilde Q$.  
This projective transformation must be such that $\tilde Q^{2n}$ is the identity 
transformation and no lower power of $\tilde Q$ is the identity.

\item  We now claim that a unit disk $R_{2n}$ with $2n$ marked points
admitting a nontrivial conformal self map $\tilde Q$ such that $\tilde Q^{2n}$ is
the lowest power of $\tilde Q$ equal 
the identity is conformal to the unit disk with the $2n$
marked points equally spaced along the boundary of the disk.

For this we present a disk  $R_{2n}$ on the upper half-plane $H = \hbox{Im} \, u \geq 0$ 
with boundary
the real line and marked points on the boundary.  The generator $Q$
of rotational symmetry cycling the $2n$ marked points 
must be implemented as a projective transformation $u\to  \tilde Q(u)$, 
which must satisfy $\tilde Q^{2n} = {\bf 1}$. 
Projective transformations come in three types: hyperbolic, parabolic, and elliptic. 
A hyperbolic transformation has two fixed points on the real line. By conjugation with
another projective transformation the two fixed points can be placed at $0$ and at $\infty$,
and then the transformation acts like $u\to c u$ with $c$ a real positive number.  Such a transformation, if nontrivial ($c\not=1$), cannot have its $2n$-th power equal to the identity. 
So $\tilde Q$ cannot
be hyperbolic.  It cannot be parabolic either, for a parabolic transformation has a single
fixed point, which can be placed at $\infty$ making the transformation take the form
$u\to u+b$ for $b$ real and nonzero.  Such a transformation  
cannot have its $2n$-th power equal to the identity.  The transformation $\tilde Q$ must
therefore be elliptic.  Elliptic transformations have two fixed points $\eta$ and $\eta^*$, one the complex conjugate
of the other.  The map $B: u \to w$ of $u \in H$ to the unit disk $|w| < 1$ is given by
\be
w = {u - \eta\over u - \eta^*} \,. 
\ee
It follows that 
$w=0$ is a fixed point of $\hat Q = B\tilde Q B^{-1}$.  
A map of a round disk to itself preserving the origin must take the form $w \to w e^{i\theta}$,
and the condition $\tilde Q^{2n}= \hat Q^{2n}= 1$ implies $w \to w e^{i\pi/n}$.  This shows that 
the disk $R_{2n}$ can be presented as the unit disk with equally spaced 
$2n$ marked points.  This is what we wanted to show.

\item  Having identified a {\em unique} point in the moduli space
with $Q$ action invariance, it follows that any dihedral invariant
metric $g$ corresponds to this point on moduli space.  
The point is also the image of the canonical flat metric $|dz|^2$
on $P_{2n}$.  This completes our proof.

\end{enumerate} 

\section{Convex programs for conformal metrics on polygons}
\label{con-pro-for-con-met}

 We now  give the explicit formulation
of the primal and dual convex programs that determine the extremal
conformal metric on the regular polygons $P_{2n}$.  By the result
of section~\ref{riem_conf} this extremal conformal metric is also the
Riemannian extremal metric.  
For a detailed
explanation of these primal and dual 
programs see~\cite{Headrick:2018ncs}, for a brief review
see~\cite{Headrick:2018dlw}.   
The dual program height parameter $\nu$, fixed at the optimum,
determines the extremal area via $A = n\nu \ell_s^2$. 
We then prove that $\nu$ also determines
the length $L_e$ of each edge of the polygon via $L_e = 2\nu \ell_s$. 

\subsection{Primal Program} \label{primal_program_sub}

\begin{figure}[!ht]
\begin{center}
\epsfysize=7.0cm
\epsfbox{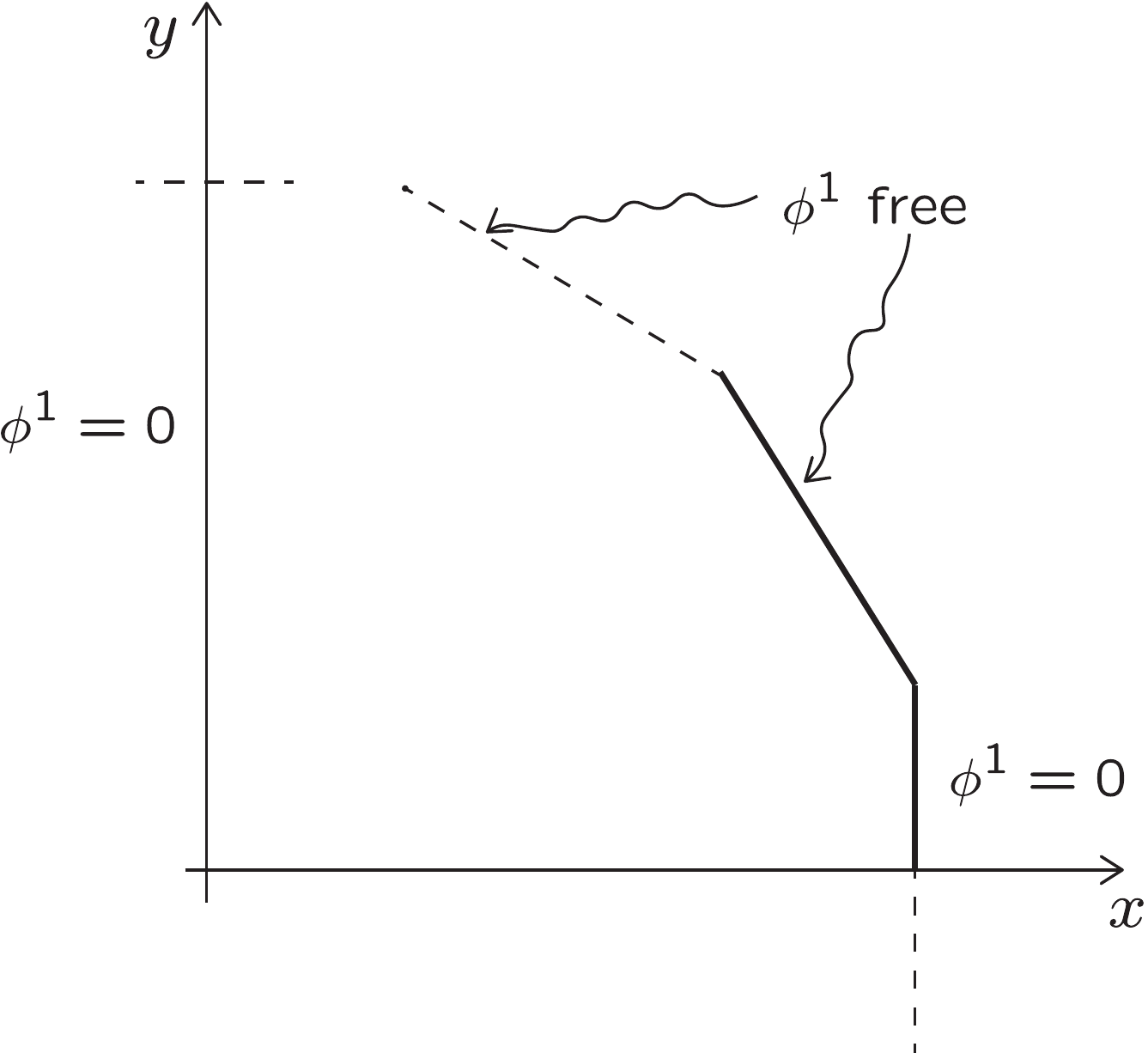}
\end{center}
\caption{\small The part $Q_{2n}$ of the polygon $P_{2n}$ on the first quadrant is
the fundamental domain for the function $\phi^1$ entering 
the calibration
$u^1$.   
The function $\phi^1$ vanishes
on the edge $\tilde e_1$ and on the $y$ axis.  There is no boundary condition
for $\phi^1$ over the other edges on the first quadrant. }
\label{fig:polycurves3}
\end{figure}

The primal program makes use of calibrations, one for each class
$C_k$ of curves on the Riemann surface.  The calibrations are 
real  
closed one-forms
 with some 
specified periods and with norm
less than or equal to one.  The norm is defined with respect to the conformal
metric on the surface.  
The calibrations 
$u^\alpha$
 take the form
\be\label{eq:cal1}
u^\alpha =  \omega^\alpha  + d\phi^\alpha\,,
\ee
where $\omega^\alpha$ is a one-form with the requisite
periods for the curves $C_\alpha$ and the function $\phi^\alpha$ 
defines the trivial part of the calibration.
The first calibration $u^1$ associated with $C_1$  is given by 
\be u^1 =  dx  + d\phi^1\,. 
\ee
The choice $\omega^1 = dx$ is consistent in that, as required,
\be
\int_{\gamma \in C_1}  dx  = 1 \,.
\ee
As discussed for the case of the Swiss cross or the torus~\cite{Headrick:2018dlw},
the symmetry of this calibration under $(x,y) \to (x, -y)$ implies
there is no contribution to $\omega^1$ from $dy$.  Since the integral $\int d\phi^1$
over any curve in $C_1$ must vanish, the value of $\phi^1$ on 
$e_1$ must be a constant equal to the value of $\phi^1$ on $\tilde e_1$.
We can choose this value to be zero:
\be
\phi^1 \Bigl|_{e_1} =  \phi^1\Bigl|_{\tilde e_1} = 0 \, . 
\ee
Since the polygon is symmetric under reflections about the $x$ and $y$ axes,
the function $\phi^1$ must define a representation under these transformations.
One finds, just like in~\cite{Headrick:2018dlw}  that
\be
\phi^1 \sim  \ ( - , + )  \,.    
\ee
meaning that $\phi^1$ is odd under $(x, y) \to (-x, y)$  (the first sign) and 
even   
under $(x, y) \to (x , -y)$ (the second sign).   We can therefore choose to work with a fundamental domain $Q_{2n}$ that is the part of the polygon on the first quadrant $x,y \geq 0$.
The antisymmetry under $x$ flips 
also implies that $\phi^1$ vanishes
on the part of the $y$ axis inside the polygon:
\be
\phi^1 (x=0, y)  = 0 \,.
\ee
Thus on the first quadrant, $\phi^1$ vanishes on the $y$ axis and on 
$\tilde e_1$.   
On the rest of the boundary
the value of $\phi^1$ must be left arbitrary.   This completes the 
specification of  $\phi^1$ (see Figure~\ref{fig:polycurves3}).

Because of the rotation symmetry of the polygon, the calibrations
$u^\alpha$ and their associated exact parts $d\phi^\alpha$, with $\alpha=2, \cdots, n$ 
are determined by $u^1$ and $\phi^1$.  
Let $(x_\alpha, y_\alpha)$ denote the points obtained by 
a counterclockwise rotation  of $(x,y)$ by an angle $\theta_\alpha= \alpha{\pi\over n}$:
\be
(x_{(\alpha)}, y_{(\alpha)} )  = \bigl(x \cos \theta_\alpha- y \sin \theta_\alpha \ , 
\ x \sin\theta_\alpha + y \cos \theta_\alpha \bigr)  \,,  \quad \theta_\alpha =  \alpha {\pi\over n} \,. 
\ee
The various calibrations are mapped into each other by these rotations, thus
we have
\be
\label{uk_from_u1}
u^{\alpha+1} (x_{(\alpha)}, y_{(\alpha)} ) =  u^1 (x, y) \,.
\ee
More explicitly this gives
\be\label{eq:calib}
u^{\alpha+1} (x, y ) =  \cos\theta_\alpha \, dx + \sin\theta_\alpha dy  + d\phi^1 (x_{(-\alpha)}, y_{(-\alpha)})  \,.
\ee
Comparing this result with  the notation $u^\alpha = \omega^\alpha + d\phi^\alpha$ we have
\be\label{eq:omega}
\omega^{\alpha+1} =  \cos\theta_\alpha \, dx + \sin\theta_\alpha dy\,, \quad
\hbox{and} \quad \phi^{\alpha+1}(x,y)  = \phi^1 (x_{(-\alpha)}, y_{(-\alpha)})\,. 
\ee

If we know $\phi^1$ over $Q_{2n}$, we know $\phi^1$ everywhere, 
and therefore we know all the $\phi^\alpha$ everywhere.  The full set of calibrations is then
known.  
Due to the choice of $\omega^\alpha$ 
and the boundary conditions on $\phi^\alpha$, all calibrations $u^\alpha$
satisfy the integral conditions on the curves in $C_\alpha$.

The program uses the  fiducial  norm squared of an arbitrary one-form, given by: 
\be
\bigl| \alpha_x dx + \alpha_y dy \bigr|_0^2   
\equiv  \alpha_x^2 + \alpha_y^2 \,.
\ee
The primal program  then becomes
\begin{equation}\label{thirdprogram-vm99}
\begin{split}
&\text{Minimize }\, \int_{P_{2n}} dx dy \,\Omega\,  \ \quad \hbox{over}   
\ \Omega, \ \phi^\alpha  , \qquad  \alpha = 1 , \cdots, n \\
& \hbox{subject to:} \ \ \ 
 \bigl|u^\alpha|_0^2\,   -\Omega\le0, \quad   \alpha = 1 , \cdots, n.
\end{split}
\end{equation}
Using  
the dihedral symmetry  the objective can be written as
\be
2\cdot (2n) \int_{T_{2n}} dx dy \,\Omega\,  ,
\ee
and the optimization is carried out applying the constraints $\Omega \geq\bigl|u^\alpha|_0^2$
over the region $T_{2n}$.

\subsection{Dual  Program}  \label{dual_program_subsection}

This time we need to determine functions $\varphi^\alpha$ on the 
polygon $P_{2n}$ with prescribed discontinuities $\nu^\alpha$ across 
chosen curves $m^\alpha\in C_\alpha$.  At the optimum the curves of constant $\varphi^\alpha$
(in a region where $\varphi^\alpha$ is not a constant) 
represent the saturating curves in the class $C_\alpha$.

Consider first the function $\varphi^1$ associated with the class $C_1$.
For this function we select the discontinuity to occur on the real $x$ axis,
along the curve $x \in [-\tfrac{1}{2}, \tfrac{1}{2}], y=0$.  We set 
\be
\varphi^1 (x, 0^+) = \tfrac{1}{2} \nu  \,,  \qquad 
\varphi^1 (x, 0^-) = - \tfrac{1}{2} \nu \,,
\ee
This is a function with discontinuity $\nu$:  The value of $\varphi^1$ increases
by $\nu$ as $y$ changes from $0^-$ to $0^+$.   The function $\varphi^1$ 
transforms as a representation under the $\mathbb{Z}_2 \times \mathbb{Z}_2$
group where the first factor refers to $(x,y)\to (-x,y)$ and the second
to $(x,y)\to (x, -y)$:   
\be
\label{varphi_signs}
\varphi^1  \sim  ( + , -)\,.    
\ee
 Consistent with the discontinuity, $\varphi^1$ changes sign under
 $y$ flips and is invariant under $x$ flips.  
Note that when calculating $\varphi^1$ derivatives in the program we do not include
the contribution from the discontinuity.
 It follows from the above equation that $Q_{2n}$ is a fundamental domain for $\varphi^1$.  If known there, it is known
all over the polygon. 

\begin{figure}[!ht]
\begin{center}
\epsfysize=7.0cm
\epsfbox{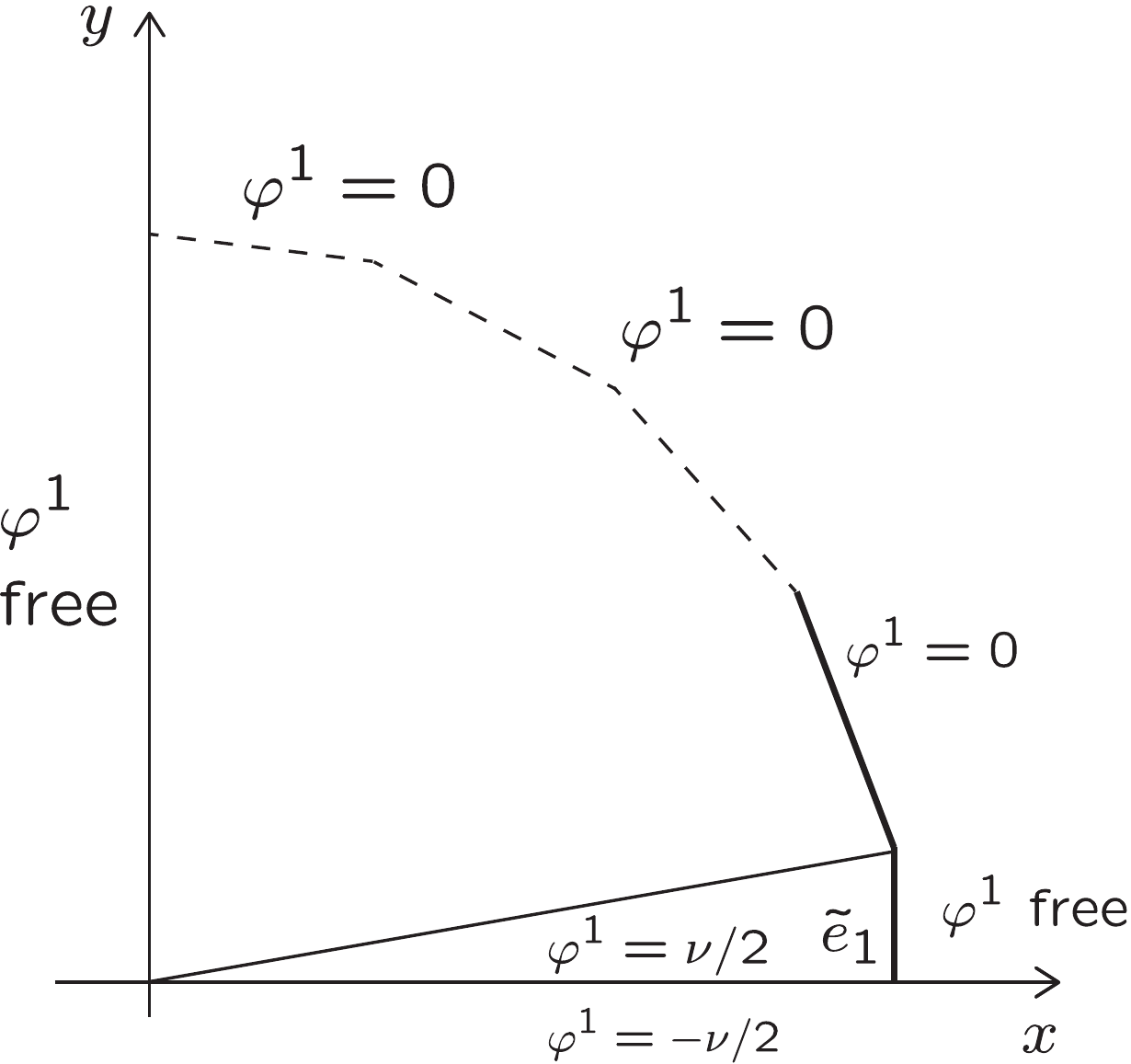}
\end{center}
\caption{\small The part $Q_{2n}$ of the polygon $P_{2n}$ on the first quadrant is
the fundamental domain for the function~$\varphi^1$.  The function is discontinuous
by $\nu$ across the $x$ axis.  It is free to vary over the edge $\tilde e_1$
where the curves in $C_1$ end, and on the $y$ axis.  $\varphi^1$
vanishes on the other polygon edges on~$Q_{2n}$.}
\label{fig:varf2}
\end{figure}

Working over $Q_{2n}$   
we have the boundary condition
$\varphi^1= \nu/2$ on the $x$ axis. Over the edge $\tilde e_1$ The function
$\varphi^1$ is free to vary over the edge $\tilde e_1$ and over
the $y$ axis.   Since no $C_1$ geodesics can end on the
other edges of the polygon and $\varphi^1$ must be continuous
except at the real axis,  $\varphi^1$ vanishes on all polygon edges on the first quadrant other   
than $\tilde e_1$ (see Figure~\ref{fig:varf2}).

As in (\ref{uk_from_u1}), the other functions $\varphi^\alpha$, 
$\alpha = 2, \ldots , n$, are determined
from $\varphi^1$ via rotations:
\be
\label{varphi_rot}
\varphi^{\alpha+1} (x_{(\alpha)}, y_{(\alpha)} ) =  \varphi^1 (x, y) \,.
\ee
Using the transformations (\ref{varphi_signs}) of $\varphi^1$ under 
reflections one can quickly write $\varphi^\alpha$ on $T_{2n}$, for 
all $\alpha = 2, \cdots, n$ in terms of the values of $\varphi^1$ on $Q_{2n}$.   Since all the $\varphi^\alpha$ are related to each other 
by rotations, the discontinuity parameters $\nu^\alpha$  are all the same
and equal to $\nu$.   

\begin{figure}[!ht]
\begin{center}
\epsfysize=7.5cm
\epsfbox{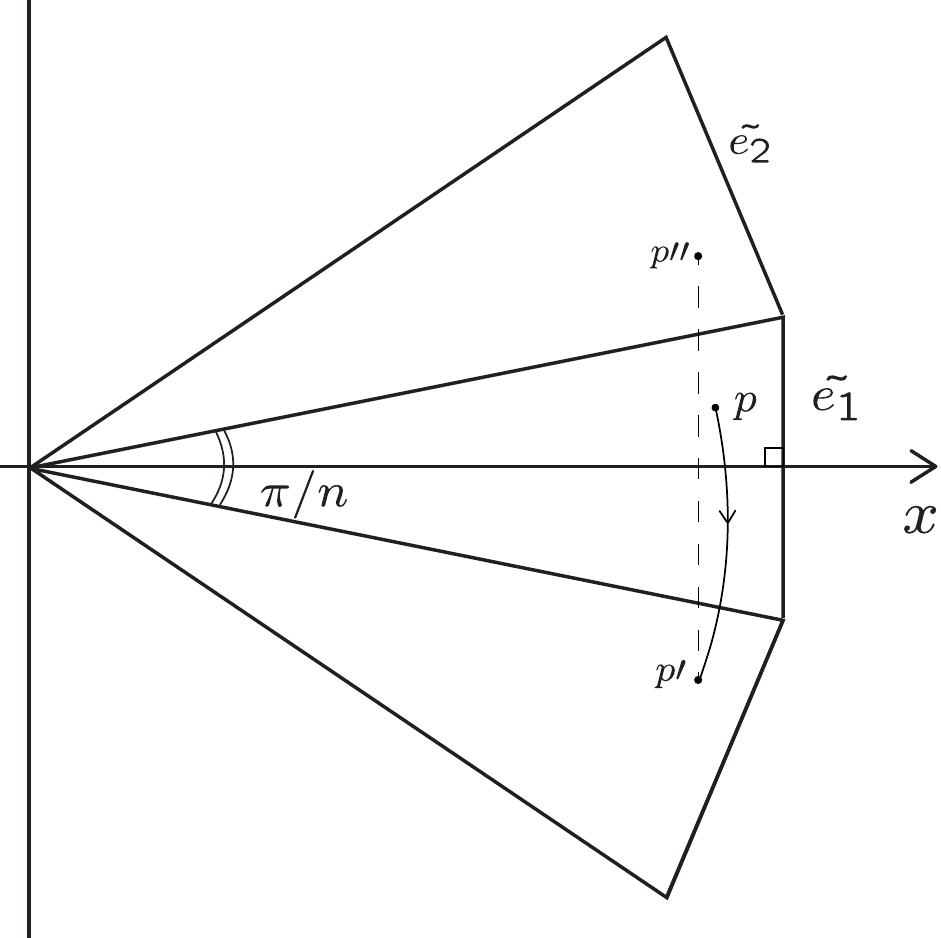}
\end{center}
\caption{\small In order to write the dual program using the fundamental
domain $T_{2n}$ we need to evaluate all $\varphi^\alpha$ on this domain.
As illustrated here $\varphi^2$ at $p$ can ultimately be related to the value
of $\varphi^1$ at $p''$, the point obtained from $p$ by a clockwise rotation 
by $\pi/n$ followed by a reflection about the $x$ axis. 
}
\label{fig:varf1}
\end{figure}

An example, shown in Figure~\ref{fig:varf1},
applies to the evaluation of $\varphi^2$ on $T_{2n}$.  We have,
by application of (\ref{varphi_rot}) for $k=1$, 
\be
\varphi^2 (p) =  \varphi^1 (p') \,,
\ee
where $p'$ is obtained from $p$ by a clockwise rotation by $\theta_1 =\pi/n$.
We also have, on account of (\ref{varphi_signs}),
\be
\varphi^1 (p') = - \varphi^1 (p''), 
\ee
where $p''$ is the reflection of $p'$ about the $x$ axis.  Therefore
\be
\varphi^2(p) = - \varphi^1 (p'')\,,
\ee
expressing, as desired, $\varphi^2$ on $T_{2n}$ in terms of $\varphi^1$ on
$Q_{2n}$.  Note, incidentally, that $p$ and $p''$ are related by reflection
about the line joining the origin to the vertex at the top end of $\tilde e_1$. 

The dual program as stated in~\cite{Headrick:2018dlw}, eqn.(2.7), takes the general form
\begin{equation}\label{thirddual}
\begin{split} 
 &\quad\text{Maximize}\ \  
  \, 2\sum_\alpha \nu^\alpha \ell_\alpha-\int_{M'}
\omega_0\Bigl(\sum_\alpha\left|d\varphi^\alpha\right|_0\Bigr)^2\   \ \     \hbox{over $\nu^\alpha$ (constants), $\varphi^\alpha$ (functions)}
\\
&\quad \hbox{subject to} \ \ \ \ 
 \Delta\varphi^\alpha|_{m_\alpha}=-\nu^\alpha\,, \\
 &  \hskip70pt  \varphi^\alpha|_{\partial M}=\ 0\, ,  \quad \forall \alpha\in J \,. 
\end{split} 
\end{equation}
In our application $\ell_\alpha = 1$ and $\nu^\alpha = \nu$ for all $n$ values of $\alpha$.
Moreover, $\omega_0 = dx\wedge dy$ and $M'$ is the polygon $P'_{2n}$ with the prime 
reminding us not to include the discontinuities in the evaluation of derivatives.
While the above statement indicates that $\varphi^\alpha$ is to vanish at the boundaries
of the manifold, this requires some 
qualification:    
For the 
class $C_\alpha$, the edges $e_\alpha, \tilde e_\alpha$ where the curves begin and end
 do not correspond to boundaries, 
 and the values of $\varphi^\alpha$ are unconstrained there.  
All other edges {\em are} boundaries and $\varphi^\alpha$ must vanish on them.  
With these comments the program becomes:
\begin{equation}\label{dualpolygon}
\text{Maximize}\ \  
  \, 2n\nu-\int_{P_{2n}'}
dx dy \, \Bigl(\sum_\alpha\left|d\varphi^\alpha\right|_0\Bigr)^2\   \ \     \hbox{over $\nu$ (constant), $\varphi^1$ (function)}
\end{equation}
Under the symmetry conditions that we impose, the integral over $P'_{2n}$ is equal to $4n$ times
the integral over $T_{2n}$:
\be
\int_{P_{2n}'}
dx dy \, \Bigl(\sum_\alpha\left|d\varphi^\alpha\right|_0\Bigr)^2 = \, 4n \, \int_{T_{2n}}
dx dy \, \Bigl(\sum_\alpha\left|d\varphi^\alpha\right|_0\Bigr)^2\,.
\ee
  As discussed above, this integral can be evaluated 
if we know $\varphi^1$ over $Q_{2n}$.  Thus $\nu$ and $\varphi^1$
over $Q_{2n}$ (with its boundary conditions), are the maximization variables.  

As shown in \cite{Headrick:2018ncs} in the general dual problem
the extremal area $A$ and the extremal value of the $\nu^\alpha$ 
are related as $A = \sum_\alpha \nu^\alpha \ell_\alpha$.  For our case, with all 
$\nu$'s identical and $\ell_s =1$  this gives for the polygon $P_{2n}$ 
an extremal area $A_{2n}$ of value
\be
\label{annuval}
A_{2n} =   n   \, \nu \,.
\ee

\subsection{Edge length and height parameter}\label{edg_len_and_hei}

In this section we show that the edge length $L_e$ of the polygon on the extremal metric is also related
to the value of the $\nu$ parameter.  In this way $L_e$ is also related
to the extremal area.

\begin{figure}[!ht]
\leavevmode
\begin{center}
\epsfysize=11.0cm
\epsfbox{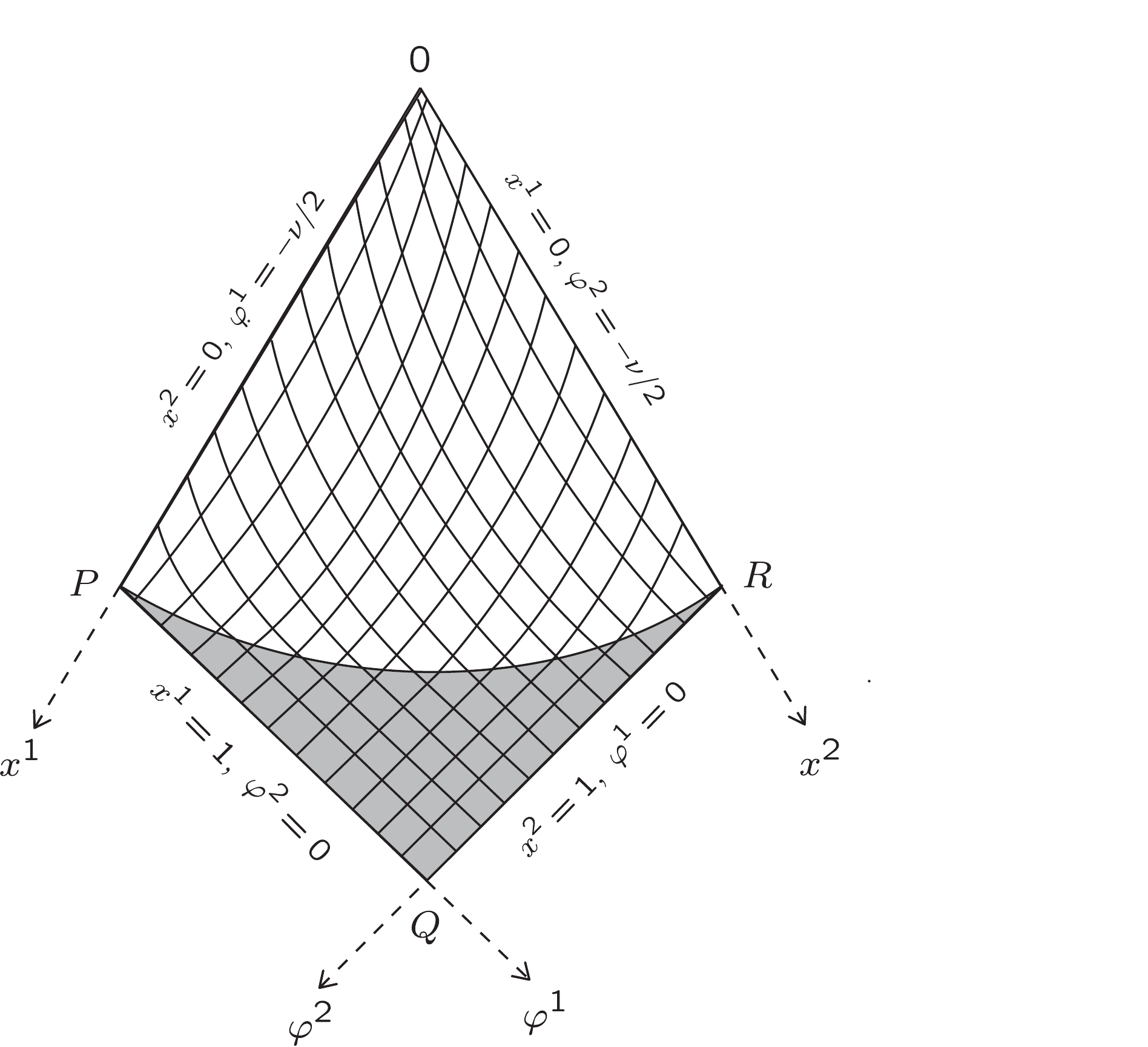}
\end{center}
\caption{\small  In the regions of the polygon where there are exactly
two bands of systolic geodesics the bands are orthogonal.  The curves
ending on the edge $e_1$ (the PQ segment) belong to the first band.  
They are parameterized by $x^1$, and
are curves of constant $\varphi^1$.  
The curves
ending on the edge $e_2$ (the QR segment) belong to the second band.  
They are parameterized by $x^2$, and are curves of constant $\varphi^2$.  
Indicated are also the geodesics with $\varphi^1= -\nu/2$ and with 
$\varphi^2 = -\nu/2$, with endpoints the midpoints of hexagon edges.}
\label{ffcl}
\end{figure}

We consider a vertex of the polygon $P_{2n}$.  As we will discuss further later,
for each vertex of the polygon 
there exists a region containing the vertex and half of the edges emerging 
from that vertex in which there are exactly
 two bands of systolic geodesics.  As established in~\cite{calabi} the two bands are orthogonal and the metric is flat in this region.  This situation is illustrated in Figure~\ref{ffcl}. 
 There we show edges $e_1$ and $e_2$ attached to a vertex $Q$
 and the geodesics that end on
   them. 
 The figure
 uses the conformal frame in which the angle at the vertex $Q$
 is $\tfrac{\pi}{2}$, rather than the original frame where the angle is
 the internal angle of the polygon at a vertex. 
  Point $P$ is the midpoint of $e_1$ and point 
 $R$ is the midpoint of $e_2$.  
 The line $PO$ is the central geodesic of the first band, the band that ends on $e_1$,
 while the line $RO$ is the central geodesic of the second band, 
 the band that ends on $e_2$. The region bounded by the curved line from $P$ to $R$
 and the edges $e_1$ and $e_2$ is the region with two bands of geodesics.
 The curve $PR$ is the boundary geodesic of another band and above this curve
  we have a region with three or more bands of geodesics.
 It is crucial for the argument below that the edge lines $PQ$ and $QR$ are boundaries
 of a two-band region.  Thus the boundary geodesic $PR$ runs from one edge midpoint to the next.  This was assumed but not proven in~\cite{calabi} and is not proven here either.

As discussed in~\cite{Headrick:2018ncs} we can construct Gaussian
coordinate systems adapted to a given set of systolic geodesics.  
We have coordinates $(x^1 , \varphi^1)$ associated with the systolic
geodesics that end on $e_1$ 
and $(x^2, \varphi^2)$ associated with the systolic geodesics
that end on $e_2$.  
The extremal metric
in this common region can be written in two equivalent forms
\begin{equation}
\label{Gaussian}
\begin{split}
ds^2 = & \  (dx^1)^2+{1\over (h_1)^2} (d\varphi^1)^2\,, \\
ds^2 = & \  (dx^2)^2+{1\over (h_2)^2} (d\varphi^1)^2\,. \\
\end{split}
\end{equation}
Clearly, $x^1$ parameterizes by length the geodesics
that end on $e_1$ and $x^2$ parameterizes by length the 
geodesics that end on $e_2$.  Since the geodesics have length one, and 
we define the range of $x^1$ and $x^2$ symmetrically, the edge $e_1$
is a line with $x^1=1/2$ and the edge $e_2$ is a line with $x^2=1/2$.
Here $h_1 ( x^1, \varphi^1)$ and $h_2 (x^2, \varphi^2)$ are two functions
that quantify the density of systolic geodesics.   Moreover,
as was also proven in~\cite{Headrick:2018ncs}, over the region
where 
only two bands exist, these functions satisfy
\begin{equation}
\label{solution3}
| h_1|+ | h_2| \,  = \, 1\, .    
\end{equation}
We will now show that these functions are in fact constants over the two band region.   
For this purpose a coordinate system defined by
 $x^1$ and $x^2$ above is useful.  
Recalling that $|dx^1 | = |dx^2| =1$ the inverse metric takes the form
\begin{equation}
g^{-1} = \begin{pmatrix} 1 & f \\ f & 1 \end{pmatrix} 
\,,  \quad \  f \ = \ \langle dx^1, dx^2\rangle \ = \ 
\langle \hat u^1 , \hat u^2 \rangle   \,.
\end{equation}
Here $\hat u^i$ is the vector associated with the one-form
 $u_i$ and it is
tangent to the geodesics in the $i$-th band.   Since the geodesics are orthogonal $f=0$  and the metric simplifies tremendously:
\be
\label{newcsvm}
ds^2 \ = \ (dx^1)^2  + (dx^2)^2  \,. 
\ee
We now exploit the compatibility of the three forms of the metric.
Consider the function $\varphi^1 (x^1, x^2)$ expressing the coordinate
$\varphi^1$ in terms of the values of the coordinates $x^1$ and $x^2$.
Then back on the first metric in (\ref{Gaussian}) and comparing with 
(\ref{newcsvm}) we have 
\be
\ (dx^1)^2+{1\over (h_1)^2} (\partial_1\varphi^1 \, dx^1 + \partial_2 \varphi^1 \, dx^2)^2  \ = \ (dx^1)^2  + (dx^2)^2\,.
\ee
Two conditions follow:
\be
\partial_1 \varphi^1 = 0 \,,  \quad  \partial_2 \varphi^1  =  \pm h_1  \,. 
\ee
Using the second form for the metric in (\ref{Gaussian}), analogous conditions are obtained 
for $\varphi^2$:
\be
\partial_2 \varphi^2 = 0 \,,  \quad  \partial_1 \varphi^2  =  \pm h_2  \,. 
\ee
The first set of equations tells us that $\varphi^1$ is just $\varphi^1(x^2)$,
with no $x^1$ dependence.  It also tells us that 
\be
\label{ghej}
{d\varphi^1\over dx^2}   (x^2)  =  \pm h_1 (x^1, \varphi^1(x^2)) \,. 
\ee  
Similarly $\varphi^2$ is just $\varphi^2(x^1)$,
and 
\be
\label{vmbb}
{d\varphi^2\over dx^1}   (x^1)  =  \pm h_2 (x^2, \varphi^2(x^1)) \,. 
\ee  
In (\ref{ghej}) the left-hand side does not depend on $x^1$ and
in (\ref{vmbb}) the left-hand side does not depend on
$x^2$.   The same must hold for the right-hand sides and 
therefore 
$h_1 = \tilde h_1 (x^2)$, is just
a function of $x^2$ and 
$h_2 = \tilde h_2 (x^1)$.  But then the condition (\ref{solution3}) gives
\be
|\tilde h_1 (x^2)|  + |\tilde h_2 (x^1) | = 1 \,.
\ee
Since $x^1$ and $x^2$ can be varied independently, the only possible solution is
that both functions are constants.  We have thus learned that both $|h_1|$ and 
$|h_2|$
are constants that add up to 1.

Note now that the geodesics with $\varphi^1= -\nu/2$ and $\varphi^2= -\nu/2$
hit the polygon
boundary at the 
middle of the edges, points $P$ and $R$, respectively; this is simply the way we defined the functions $\varphi$.  
Let $L_e$ denote the length of the edges in the extremal
metric. Because of the symmetry of this metric the 
half-length $\tfrac{1}{2} L_e$ is both the distance from $P$, where $\varphi^1= -\nu/2$,
 to the vertex $Q$, where $\varphi^1=0$,
as well as the distance from $R$, where $\varphi^2 = -\nu/2$,
to the vertex $Q$, where $\varphi^2=0$.  For the first, $x^1$ is constant, and 
for the second $x^2$ is constant. We thus have, integrating
the length element over these intervals:
\be
\begin{split}
\tfrac{1}{2} L_e = & \ \int_{-\nu/2}^0  {d\varphi^1\over |h_1| } =  {\nu\over 2|h_1|} \,,  
 \\
\tfrac{1}{2} L_e = & \  \int_{-\nu/2}^0  {d\varphi^2\over |h_2| } =  {\nu\over 2|h_2|}  \,. 
\end{split}
\ee  
This implies $|h_1| = |h_2| = \tfrac{1}{2}$.   As a result,  
\be
\label{le-equal-two-nu}
L_e = 2 \nu \,.  
\ee
This is what we wanted to show and holds for all $2n$-gons.

 Given the above result for the edge length $L_e$,
 the perimeter  ${\cal P}_{2n}$ of the polygon is given by 
 ${\cal P}_{2n} = 2n L_e =  4 n \nu$. Recalling  from (\ref{annuval}) that the extremal area is $A_{2n} = n\nu$ we find the perimeter/area relation:
 \be
 \label{avp-exact}
A_{2n}  \, = \,  \tfrac{1}{4}\,  {\cal P}_{2n} \,,  \qquad  \ell_s=1\,.
 \ee
 This is for the value $\ell_s =1$ of the systole.  If we introduce an arbitrary
  systole $\ell_s$ back
 into this formula by scaling of the extremal metric we 
\be
A_{2n}  \, = \,  \tfrac{1}{4}\, \ell_s\,  {\cal P}_{2n} \,.
 \ee
 This formula is as if the surface was built (which is not)
 by piling flat little triangles of height $\ell_s/2$ 
all over the length of the perimeter.  
  In fact for the flat metric $\rho=1$ on 
 a regular $2n$-sided polygon, one quickly notes that the area $\bar A_{2n}$
 and the perimeter $\bar{\cal P}_{2n}$ are related by 
 \be
\bar A_{2n}  \ = \ \tfrac{1}{2} \cdot a \cdot \bar {\cal P}_{2n} \,,
 \ee
where $a$ is the apothem.  Since we take $a = {1\over 2}$ in our canonical
presentation we have 
\be
\label{flat_area_perimeter_apo-half}
\bar A_{2n}  = \tfrac {1}{4} \, \bar {\cal P}_{2n}\, ,  \qquad    a = \tfrac{1}{2}\,, 
\ee
 the same relation that holds for the extremal metric with $\ell_s=1$.  
 It follows that with $\ell_s=1$
and apothem ${1\over 2}$:
\be
{{\cal P}_{2n}  \over \bar {\cal P}_{2n} } =  {A_{2n} \over \bar A_{2n}}   <  1\,.
\ee
In passing from the flat metric to the extremal metric, the fractional 
reduction in 
area  equals the fractional reduction of the perimeter, or the fractional
reduction of the edge length.

\section{Calabi's hexagon} \label{cal_bi_hex}
In this section we discuss the extremal problem for Calabi's hexagon $P_6$. 
We give the setup of the primal and the dual programs in detail.
We then present our numerical solution, which gives quantitative support for the qualitative picture provided by Calabi.  We also discuss the distribution of 
curvature on the surface showing there is positive curvature on the three-band
region and line curvature on the boundary of this region.

\subsection{Setup for convex programs}
\begin{figure}[!ht]
\leavevmode
\begin{center}
\epsfysize=8cm
\epsfbox{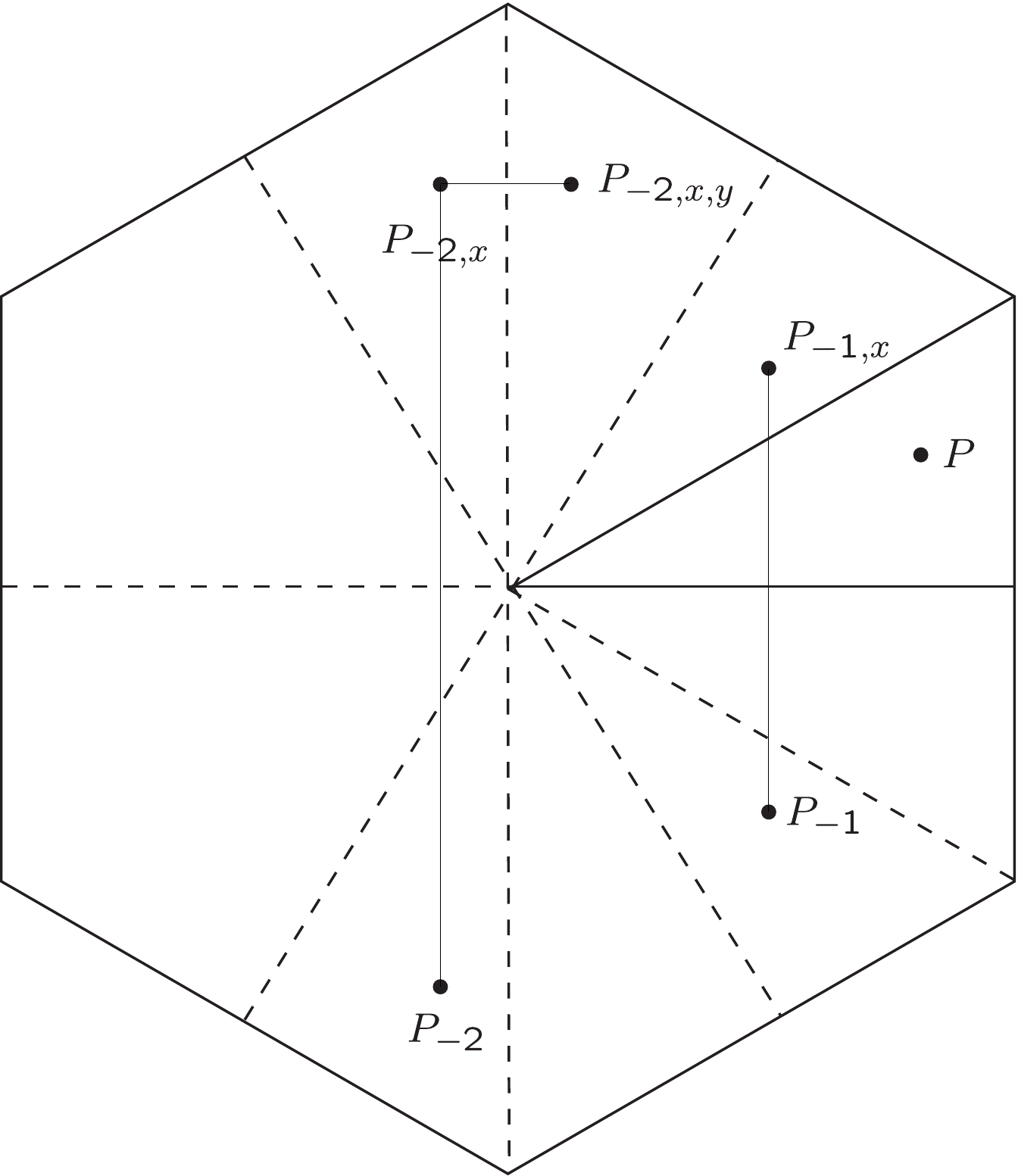}
\end{center}
\caption{\small  The values of $\phi^2$ and $\phi^3$ over $T_6$ are related to values of $\phi^1$ over $Q_6$. The point $P_{-1}$ is a rotation of $P\in T_6$ by $-\tfrac{\pi}{3}$. $P_{-1,x}\in Q_6$ is the reflection of $P_{-1}$ about the $x$-axis. The point $P_{-2}$ is a rotation  of $P$ by $-\tfrac{2\pi}{3} $. $P_{-2,x}$ is the reflection of $P_{-2}$ about the $x$-axis. $P_{-2,x,y}\in Q_6$ is the reflection of $P_{-2,x}$ about the $y$-axis.}
\label{fig:hexafns}
\end{figure}

Here we formulate the primal and the dual programs for the hexagon explicitly. We will follow the general theory developed in sections~\ref{primal_program_sub} and 
\ref{dual_program_subsection} and specialize to $n=3$.
From eqs. (\ref{eq:cal1}) and (\ref{eq:omega}) the three calibrations can be written as
\be\label{eq:hexacalib}
\begin{split}
u^1&\ =\ dx  + d \phi^1, \\
u^2&\ = \ \tfrac{1}{2}dx+\tfrac{\sqrt{3}}{2} d y + d \phi^2, \\
 u^3& \ =\   -\tfrac{1}{2}dx+\tfrac{\sqrt{3}}{2} d y + d \phi^3.
\end{split}
\ee
Due to the hexagonal symmetry we only need to compute the minimal area metric on $T_6$. 
For the primal program we need the three functions $\phi^1,\phi^2$ and $\phi^3$ on $T_6$. Using the symmetries of the hexagon one can define $\phi^2$ and $\phi^3$ in terms of $\phi^1$ if we know the value of $\phi^1$ over $Q_6$, the part of $P_6$  in the first quadrant.    
We can explain this clearly using Figure~\ref{fig:hexafns}. The value of $\phi^2$ at a point $P\in T_6$ is equal to the value of $\phi^1$ at the point $P_{-1}$, which is same as the value of $\phi^1$ at the point $P_{-1,x}\in Q_6$. 
This lets us determine $\phi^{2}$ everywhere on $T_6$ with the formula: 
\be
  \phi^2\bkt{x,y}=\phi^1\bkt{\tfrac12 x+\tfrac{\sqrt{3}}{2} y, \tfrac{\sqrt{3}}{2} x-\tfrac12 y}.
\ee
In a similar fashion, the value of $\phi^3$ at point $P\in T_6$ is equal to the value of $\phi^1$ at the point $P_{-2}$. This is same as the value of $\phi^{1}$ at the point $P_{-2,x}$  which is negative of the value of $\phi^1$ at the point $P_{-2,x,y}\in Q_6$. This lets us determine $\phi^3$ over $T_6$ in terms of $\phi^1$ over $Q_6$:
\be 
\phi^3\bkt{x,y}=-\phi^1\bkt{\tfrac{1}{2} x-\tfrac{\sqrt{3}}{2} y, \tfrac{\sqrt{3}}{2} x+\tfrac12 y}\,.
\ee
The boundary conditions for $\phi^1$ on $Q_6$ are those described
in Figure~\ref{fig:polycurves3}:   $\phi^1$ is zero on $\tilde e_1$ and on the $y$
axis, and it is free on the other edges on $Q_6$.
The program is then explicitly defined by
\begin{equation}\label{eq:pphexa}
\begin{split}
&\text{Minimize }\, 12 \int_{T_{6}} dx dy \,\Omega\,  \quad \hbox{over}
\  \ \Omega\ \  \hbox{and}\ \  \phi^1\,, \\
& \hbox{subject to:} \ \ \ 
 \bigl|u^i|_0^2\,   -\Omega\le0, \quad  i = 1,2,3.  
 \end{split}
\end{equation}
The fiducial norms of the one-forms are given by
\be
\begin{split}
\bigl|u^1|_0^2& = \bkt{1+\p_x \phi^1}^2+\bkt{\p_y\phi^1}^2\, ,
\\
\bigl|u^2|_0^2&= \bkt{\tfrac12+\p_x \phi^2}^2+\bkt{\tfrac{\sqrt{3}}{2}+\p_y\phi^2}^2\, ,
\\
\bigl|u^3|_0^2&= \bkt{-\tfrac12+\p_x \phi^3}^2+\bkt{\tfrac{\sqrt{3}}{2}+\p_y\phi^3}^2\, .
\end{split}
\ee

\bigskip
We now give the explicit setup for the dual program.
Having three homology classes of curves we need to define three functions $\varphi^1, \varphi^2$, and $\varphi^3$ over $P_6$. Each function must have a discontinuity across some arbitrarily selected curve in the corresponding class. We choose $\varphi^1$ to be discontinuous across the $x$-axis:
\be
\varphi^1\bkt{x,0^-}=-\tfrac{\nu}{2}\quad\hbox{and}\quad\varphi^1\bkt{x,0^+}=\ \tfrac{\nu}{2}.
\ee
The function $\varphi^1$ is free to vary on the edge $\tilde e_1$ and on the $y$ axis and vanishes on
the other edges in $Q_6$ (as was summarized in Figure~\ref{fig:varf2}). 
The hexagonal symmetry can then be used to determine $\varphi^2,\varphi^3$ in terms of $\varphi^1$ as we did in the case of the primal program. 
This time, however,  $\varphi^1$ is odd under reflections about the 
$x$ axis and even under reflections about the $y$ axis.   This gives, 
\be
\begin{split}
  \varphi^2\bkt{x,y}& \ =\ -\varphi^1\bkt{\tfrac12 x+\tfrac{\sqrt{3}}{2} y, \tfrac{\sqrt{3}}{2} x-\tfrac12 y} \,, \\
\varphi^3\bkt{x,y}& \ =\ -\varphi^1\bkt{\tfrac{1}{2} x-\tfrac{\sqrt{3}}{2} y, \tfrac{\sqrt{3}}{2} x+\tfrac12 y}.
\end{split}
\ee
The dual program is then explicitly written as
\begin{equation}\label{eq:hexadual}
\text{Maximize}\ \  
  \, 6\nu-12 \int_{T_6}
dx dy \, \Bigl(\left|d\varphi^1\right|_0+\left|d\varphi^2\right|_0+\left|d\varphi^3\right|_0\Bigr)^2\   \ \     \hbox{over $\nu$ and $\varphi^1$}.
\end{equation}
A numerical solution of the program requires discretization (see appendix~\ref{discretization_appendix}).

\subsection{Numerical results, metric description, and consistency checks}
\begin{figure}[!ht]
\leavevmode
\begin{center}
\epsfysize=14.0cm
\epsfbox{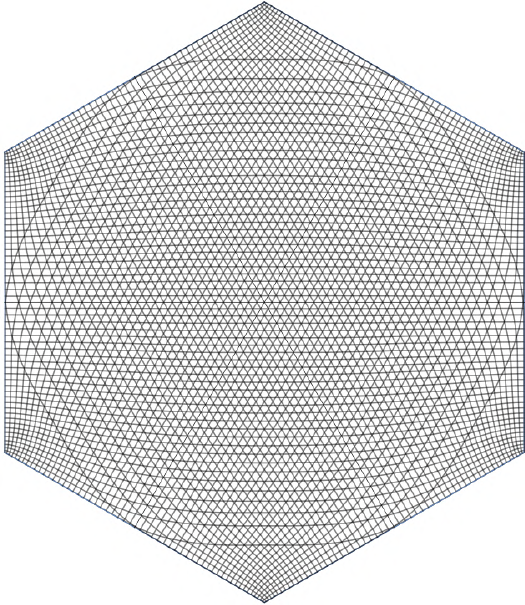}
\end{center}
\caption{\small  Saturating geodesic curves in the minimal area metric
for Calabi's hexagon. The surface
is covered by a connected region $U_3$ with three geodesic bands and six regions
comprising the two-band part $U_2$,  each region around a vertex.   
The metric is flat on $U_2$ but has Gaussian curvature on $U_3$ and
line curvature supported on the boundary of $U_3$.}
\label{fig:hexfinal}
\end{figure}

We now present our results from the numerical evaluation of the primal and the dual programs for the hexagon.  From the solution of the dual program one can construct the saturating geodesics in classes $C_1, C_2$, and $C_3$ as curves of constant 
$\varphi^1, \varphi^2$, and $\varphi^3$, respectively. Doing so results in
the picture of $P_6$ shown in Figure~\ref{fig:hexfinal}. 
There are regions $U_3$ and $U_2$ with three and two bands of geodesics, 
respectively. The central region of the hexagon defines the whole of $U_3$. 
Finite neighborhoods of each vertex are in $U_2$.  There are no regions $U_1$ or $U_0$.  
The boundary of the region $U_3$ touches the edges of the polygon 
at  the midpoints.  It is also clear from the figure and the data that, as expected, 
the geodesics in $U_2$ intersect at right angles -- a property  proven 
in~\cite{calabi}. 
  This means that at the vertices of the polygon
we have effectively an internal angle of $\pi/2$, rather than the apparent
$2\pi/3$.   The extremal 
conformal metric in fact diverges as we approach the vertices,
and the conformal map turning the angle from $2\pi/3$ to $\pi/2$ renders
the metric finite.  

Each band of geodesics is defined by a pair of opposite edges $e_i$ and $\tilde e_i$, and each
edge determines two vertices at its endpoints.  
 Each band
has two {\em boundary} geodesics, each departing from a vertex in $e_i$
and ending on a vertex in $\tilde e_i$. As conjectured in~\cite{calabi}, 
each  boundary geodesic follows the (neighboring) edge until the midpoint 
before veering off into the
inside of the polygon.  It is clear it cannot veer off before reaching the midpoint,
as then we would have a region near the midpoint of each edge with just one
band of geodesics, which is impossible in a Riemannian  isosystolic problem.   If it were
to veer off after the midpoint there is no obvious contradiction, but the situation
would be less symmetric, with $U_3$ overlapping with the edges over finite
regions.  In section~\ref{edg_len_and_hei} 
we showed that the edge length $L_e$
is equal to $2\nu$ using the assumption that boundary geodesics
veer off at midpoints.  We will see that the data below gives evidence for this relation.

\medskip
A few more observations follow from the use of the Gauss-Bonnet
formula giving the Euler number $\chi$ of a surface with a boundary
\be\chi  = {1\over 2\pi} \int K\, dA + {1\over 2\pi} \int k ds + {1\over 2\pi}
\sum_i (\pi -\theta_i) \,, \ee
Here $K$ is the Gaussian curvature, $k$ is the curvature on the boundary,
which vanishes when the boundary is a geodesic, and 
$\theta_i\in(0, 2\pi)$    is the internal angle at corners that appear on the boundary. 
Since $U_3$ is topologically a disk it has $\chi=1$. 
 The boundary of $U_3$ is a connected curve 
 comprised of boundary geodesics that 
 join at edge midpoints without forming corners. It follows that  the second
 and third contributions in the Gauss-Bonnet formula vanish and
 \be
{1\over 2\pi}  \int_{U_3} K \, dA =  1 \,. 
 \ee
Now consider applying Gauss-Bonnet to the whole surface, which is still a 
topological disk. 
The internal angle at each of the six corners is $\theta_i = {\pi\over 2}$.  
Moreover $\partial U_3$ has
six arcs, each going from the midpoint of an edge to the midpoint of the
next edge.   These arcs are geodesics on $U_3$ but are not geodesics
on $U_2$.  As a result they carry line curvature singularity on the surface.  
The amount of curvature in each arc is computed from the turning angle 
of the tangent vector to the arc as seen from the flat metric in $U_2$.  
The turning angle is in
fact ${\pi\over 2}$ (a fourth of $2\pi$) 
but as a $U_2$ boundary the sign of $k$ is negative.  We thus have, for each
arc and for each corner:   
\be
 {1\over 2\pi} \int_{\rm arc} k ds  =  -\tfrac{1}{4} \,,   \qquad 
 {1\over 2\pi} (\pi -\theta_i) = +\tfrac{1}{4} \,.
\ee
Since we have six arcs and six corners the contributions to the Euler number
of the whole surface cancel: 
\be
\chi =   1  +  6 (-\tfrac{1}{4})  + 6 (\tfrac{1}{4})   = 1 \,.
\ee
It is also possible to check that  over each
band ${\cal B}$ of geodesics the total integrated curvature  vanishes:
\be
 \int_{\cal B}  K\, dA + \int k ds  = 0 \,.
\ee
This must happen because we can glue the beginning and ending edges of the band to form an annulus,  a surface with $\chi=0$.   This is simple
to confirm. 
Since any band in the hexagon contains the full $U_3$
the bulk contribution above  is $2\pi$.    
The integral $\int k ds$ 
represents the contribution of line curvature in the annulus.  There is no contribution from the boundary
of the band since those are geodesics in the band.  
But each band contains in its interior four arcs
with line curvature, the arcs on the four vertices that define the band.  We saw
that each
one contributes $-{\pi\over 2}$ to $\int k ds$, and therefore this gives a total
of $-2\pi$ canceling the bulk contribution, as expected.

Before looking at the numerical results it is useful to consider 
the reference  
flat metric $\rho=1$, which
 is admissible.  The area $\bar A_6$ of the hexagon in this metric is
\be
\bar A_6 = \tfrac{\sqrt{3}}{2} \simeq 0.866\hskip1pt 025\,.
\ee
The perimeter $\bar {\cal P}_6$ of the polygon in this metric is 
\be
\bar {\cal P}_6 =  2 \sqrt{3} \simeq  3.464\hskip1pt 102\,.
\ee
Note the exact relation 
$\bar {\cal P}_6 =  4 \bar A_6$ anticipated in (\ref{flat_area_perimeter_apo-half}).
Recall also the dual program prediction (\ref{avp-exact}):
$ {\cal P}_6 =  4 A_6$, valid for the extremal metric.

We have below data for various quantities computed with the primal and
with the dual programs at various values of the lattice resolution $N_c$
(essentially the number of subdivisions of the apothem, see Appendix A  
for more details).
To get some additional information
 from the data, we analyze convergence and 
 try extrapolation  to $N_c\to \infty$ in the 
same manner as done in~\cite{Headrick:2018dlw}.  From the data for $N_c=2,4,\cdots, 128$ we define the error $e_{q}(N_c)$ in a quantity $q$ at resolution $N_c$ by
\be
e_{q}\bkt{N_c} \equiv q\bkt{N_c}-q\bkt{128},\quad \quad N_c=2,4,\cdots 64.
\ee
We then fit $e_q(N_c)$ to a decaying exponential function,  
\be
e_{q}(N_c)\, = \, \frac{a}{N_{c}^{b}}\, ,
\ee
to  find the values of $a$ and $b$, with $b$ giving information on the 
convergence rate. The extrapolated value $q_*$ of $q$ at infinite resolution
is estimated by 
\be
q_{*}= q\bkt{128}-e_{q}\bkt{128}.
\ee

Table~\ref{tab:hexprimdata} contains the data for the area $A_6$, perimeter
 ${\cal P}_6$, and 
metric $\rho^2 = \Omega$ at the origin, as computed from the primal 
program at various values of  $N_c$, including the $N_c\to \infty$ extrapolation.   
The area goes monotonically
down, as it should upon minimization with lattices that fit evenly into each other.
The perimeter is evaluated by integrating the metric along the edges.
This is difficult to do accurately, since the discretized metric 
tries to diverge as we approach the vertices of the hexagon.  
The convergence of the perimeter data is
certainly not monotonic.  Nor is the convergence monotonic for the metric
at the center of the hexagon.

The data from the dual program is presented in table~\ref{tab:hexdualdata}. The dual program directly solves for the minimal area and the value of $\nu$. 
It also gives us the functions $\varphi^\alpha$, $\alpha = 1, 2, 3$, from which we compute the metric $\rho=\sum_{\alpha} \left|d \varphi^\alpha \right|_0$.
The table shows the value of the area $A_6$, the perimeter ${\cal P}_6$, the metric
$\rho^2= \Omega$ at the origin, and $2\nu$.  

Let us look first at the area. 
Based on the bounds obtained by the $N_c=128$ data
of the two programs  we find that the extremal area of the hexagon is in the range
\be
0.8400< A_6<0.8414\, .
\ee
The error analysis shows that convergence is slow for the primal but
faster in the dual ($b\sim 0.9$ and $b=1.6$, respectively). 
The dual extrapolation suggests that the value $A_6= 0.8401$
could be correct to four significant digits.  This area, in fact, is only
3\% lower than the flat-metric area $\bar A_6 = \tfrac{\sqrt{3}}{2}$.

The perimeter ${\cal P}_6$ can be used to test the expected relation
${\cal P}_6 = 4A_6$ which, as explained before, is derived rigorously once
we assume   that
 boundary geodesics veer off at edge midpoints.  The equation holds
 to great accuracy for $N_c=128$ in the dual program 
($3.3647/0.84002=4.0055$).  
 It holds less accurately in the primal at the same resolution
 ($3.3275/0.84135=3.95495$), a 1\% error.    The value of the metric at the origin
 in both the primal and dual programs seem consistent with each other. 
 Finally the relation $A_6 = 3 \nu$, that follows from the dual program, 
 is seen to hold very accurately.

\begin{table}[!ht]
\centering
\begin{tabular}{|l| c| c| c|}
\hline
$N_c$ & $A_6$ & ${\cal P}_6$ & $\rho^2\bkt{0,0}$\\
\hline
$2$ & $0.86602$ & $3.4641$ & $1.00000$ \\ \hline
$4$ & $0.86297$ & $3.3888$ & $1.04921$ \\ \hline
$8$ & $0.85389$ & $3.2774$ & $1.14853$ \\ \hline
$16$ & $0.84819$ & $3.2746$ & $1.17448$ \\ \hline
$32$ & $0.84458$ & $3.2894$ & $1.17072$ \\ \hline
$64$ & $0.84249$ & $3.3092$ & $1.16571$ \\ \hline
$128$ & $0.84135$ & $3.3275$ & $1.16060$ \\ \hline
$\to \infty$& $0.84028$ & $3.3275$ & $1.16167$
 \\ \hline
\end{tabular}
\caption{Primal program numerical data for the hexagon area $A_6$, perimeter
${\cal P}_6$, and metric $\rho^2(0,0)$ at the origin at various choices of lattice resolution $N_c$.}
\label{tab:hexprimdata}

\end{table}

\begin{table}[!ht]
\centering
\begin{tabular}{|l| c| c| c| r|}
\hline
$N_c$& $A_6$ & ${\cal P}_6$ & $\rho^2\bkt{0,0}$&$2 \nu$\\
\hline
$2$ & $0.80535$ & $2.9554$ & $1.15614$ &$0.5369$\\ \hline
$4$ & $0.82766$ & $2.9627$ & $1.17678$&$0.5518$\\ \hline
$8$ & $0.83581$ & $3.0499 $ & $1.15032$&$0.5571$\\ \hline
$16$ & $0.83866$ & $3.1755$ & $1.15557$&$0.5592$\\ \hline
$32$ & $0.83964$ & $ 3.2944$ & $1.15295$&$0.5598$\\ \hline
$64$ & $0.83997$ & $  3.3873 $ & $1.15659$&$0.5599$\\ \hline
$128$ & $0.84002$ & $ 3.3647  $ & $1.15680$&$0.5597$\\ \hline
$\to \infty$& $0.84007$ & $3.4251$ & $1.15656$&$0.5598$ \\ \hline
\end{tabular}
\caption{Dual program numerical data for the hexagon 
area $A_6$, perimeter
${\cal P}_6$, metric $\rho^2(0,0)$ at the origin,
  and $2\nu$ at various choices of lattice resolution $N_c$.}
\label{tab:hexdualdata}
\end{table}

\section{The cases of the octagon and decagon} \label{the_cas_of_the_oct}

In this section we discuss in some detail the octagon and then, more briefly,
the decagon.  
There are a few reasons to do this explicitly.  There are no known
extremal metrics with regions covered by a finite number $k$ of geodesic
bands with $k>3$.   
Thus the octagon
 and decagon provide, respectively,  tractable examples with regions having four 
 and five bands of geodesics.  
 The octagon has regions $U_4$, $U_3$, and  $U_2$. 
The possibility has been raised~\cite{bryant}  
that in the isosystolic problem any
$U_4$ could be forced to have zero Gaussian curvature, but our results
indicate that the integrated curvature on $U_4$ should be that of half a sphere.
(We are not able to assess the 
curvature on the region $U_3$.)  
We will also see that the region $U_4$ is rather
large, while $U_3$ and $U_2$ are small.   
It was stated in~\cite{bryant}  
that in the isosystolic problem any
$U_5$ must have zero Gaussian curvature. Our results for the decagon, however,
indicate that the integrated curvature on $U_5$ should be that of half a sphere.

\subsection{Setup for the convex programs}
 For the primal program the four calibrations are given by:
\be\begin{split}
u^1&=d x+ d\phi^1\, ,\\
u^2&=\tfrac{1}{\sqrt{2}} dx+\tfrac{1}{\sqrt{2}} dy+d\phi^2\, ,  \\
u^3&= dy+d\phi^3\,, \\
u^4&= -\tfrac{1}{\sqrt{2}} dx+\tfrac{1}{\sqrt{2}} dy +d\phi^4.
\end{split}
\ee
We need the four functions $\phi^1,\phi^2,\phi^3$, and $\phi^4$ on $T_8$. Due 
to the symmetries of the octagon we can define all 
functions on $T_8$   
in terms of just $\phi^1$ on $Q_8$, the piece of the octagon on the first quadrant. 
The precise relations are: 
\be
\begin{split}
\phi^2\bkt{x,y}&= \phi^1\bkt{\frac{x+y}{\sqrt{2}}, \frac{x-y}{\sqrt{2}}}\,, \\
\phi^3\bkt{x,y}&=\phi^1\bkt{y,x}\, ,\\
\phi^4\bkt{x,y}&=-\phi^1\bkt{\frac{x-y}{\sqrt{2}}, \frac{x+y}{\sqrt{2}}}.
\end{split}
\ee
The boundary conditions for $\phi^1$ on $Q_8$ are those described
in Figure~\ref{fig:polycurves3}:   $\phi^1$ is zero on $\tilde e_1$ and on the $y$
axis, and it is free on the other edges on $Q_8$.
The primal program for the octagon is then:
\begin{equation}\label{eq:ppocta}
\begin{split}
\text{Minimize }\, &16 \int_{T_{8}} dx dy \,\Omega\,  \quad \hbox{over}
\ \Omega\quad \hbox{and}\quad \phi^1 \,,\\
 \hbox{subject to:} \ \ \ 
 &\bigl|u^i|_0^2\,   -\Omega\le0,    \ \ \  i = 1,2,3,4.
 \end{split}
\end{equation}

For the dual program, the functions $\varphi^{2}, \varphi^3$ and $\varphi^4$ are related to $\varphi^1$ as follows:
\be
\begin{split}
\varphi^2\bkt{x,y}&= -\varphi^1\bkt{\frac{x+y}{\sqrt{2}}, \frac{x-y}{\sqrt{2}}}\, , \\
\varphi^3\bkt{x,y}&=-\varphi^1\bkt{y,x}\, , \\
\varphi^4\bkt{x,y}&=-\varphi^1\bkt{\frac{x-y}{\sqrt{2}}, \frac{x+y}{\sqrt{2}}}.
\end{split}
\ee
On $Q_8$ the boundary conditions on $\varphi^1$ are:  it is equal to $\nu/2$ on the
$x$ axis, it free to vary on the edge $\tilde e_1$ and on the $y$ axis, and vanishes on
the other octagon edges of $Q_8$. 
The dual program is then explicitly written as:
\begin{equation}\label{eq:octadual}
\text{Maximize}\ \  
  \, 8\nu-16 \int_{T_8}
dx dy \, \Bigl(\left|d\varphi^1\right|_0+\left|d\varphi^2\right|_0+\left|d\varphi^3\right|_0+\left|d\varphi^4\right|\Bigr)^2\   \ \     \hbox{over $\nu$ and $\varphi^1$}.
\end{equation}
The discretization scheme we use for the octagon is different from 
that of the hexagon and admits a simple generalization to 
$2n$-gons (see Appendix~\ref{discretization_appendix}).

\subsection{Numerical results and patterns of bands}

\begin{figure}[!ht]
\leavevmode
\begin{center}
\epsfysize=13.5cm
\epsfbox{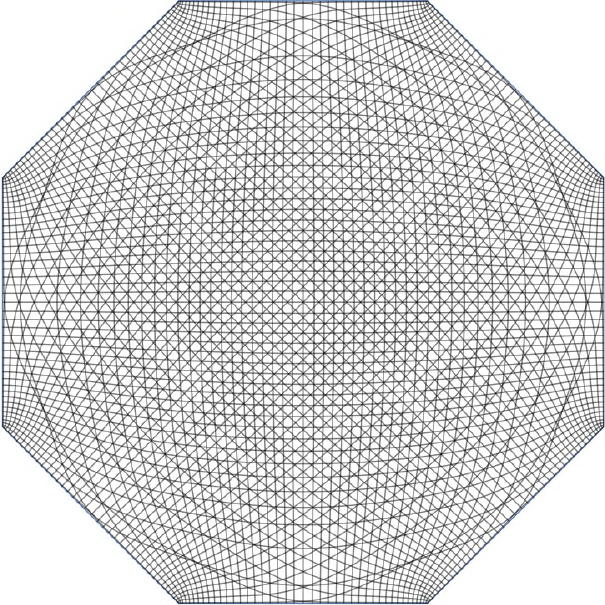}
\end{center}
\caption{\small  The saturating geodesics for the octagon $P_8$. 
The surface
is covered by a connected region with four bands of geodesics, eight regions
 with three bands, and eight regions with two bands.}
\label{fig:octfinal}
\end{figure}

Now we present the results from the primal and the dual program for the octagon. The picture of saturating geodesics on $P_8$ which is constructed from our solution is given in Figure~\ref{fig:octfinal}.

\begin{figure}[!ht]
\leavevmode
\begin{center}
\epsfysize=13.0cm
\epsfbox{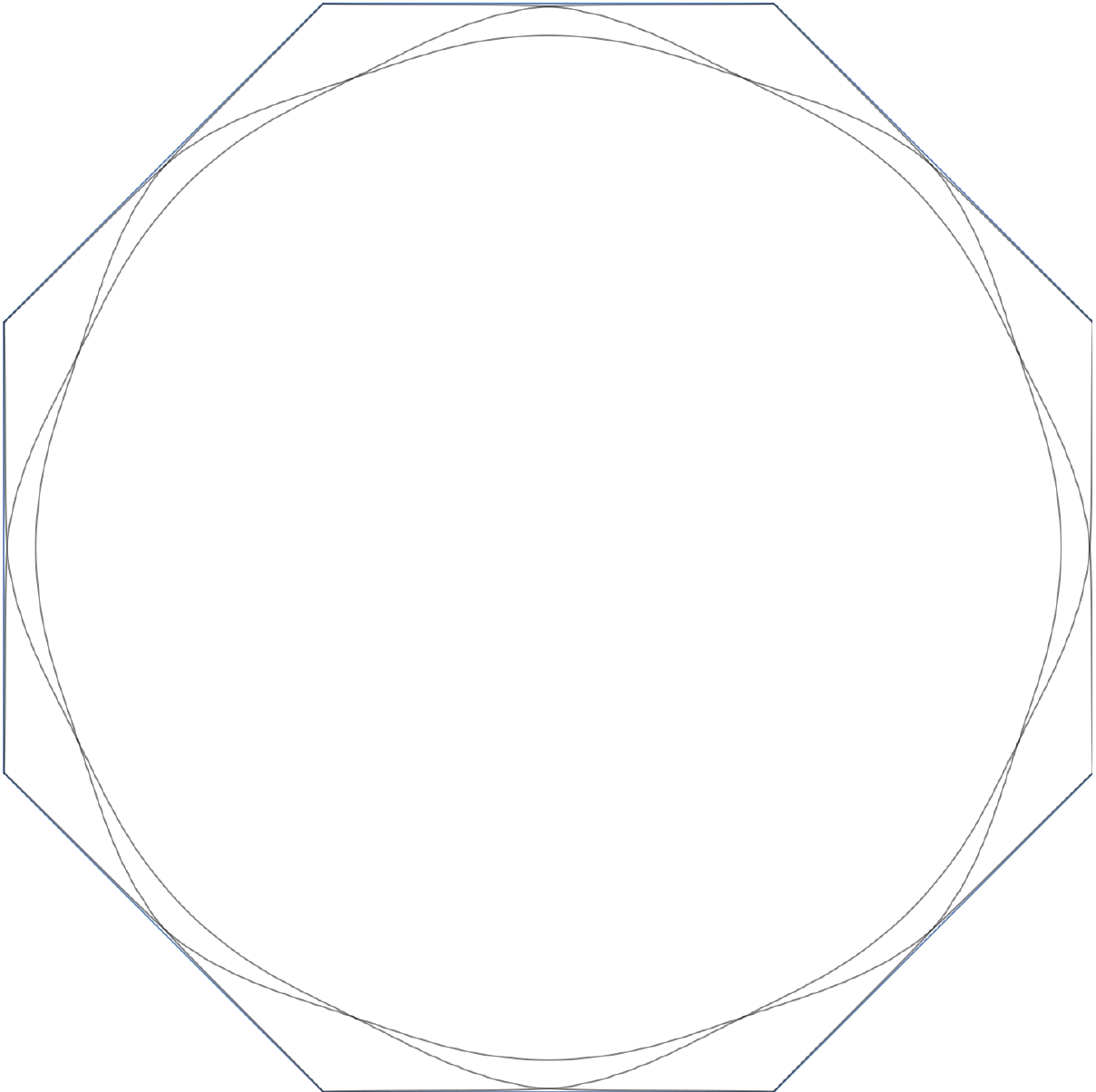}
\end{center}
\caption{\small  The boundary saturating geodesics for $P_8$.  The surface
is covered by a connected region with four bands of geodesics, eight regions
 with three bands, and eight regions with two bands.}
\label{fig:LastGeo8Gon}
\end{figure}

The central region is covered by four geodesic bands  and has positive curvature. There is a small region with three geodesics near each of the edges. 
The regions with two geodesics are flat. 
These regions are easier to see in Figure~\ref{fig:LastGeo8Gon} where we draw the boundary geodesics in classes $C_1,C_2,C_3$, and $C_4$.

The best result for the area in the primal program, giving an upper bound, 
was obtained with $N_c=128$.
For the dual program, giving us a lower bound, we could only work up to 
$N_c = 64$.  Together they give 
\be
0.7776  <  A_8  <   0.7804  \,.
\ee
As in the case of the hexagon, the convergence is faster in the dual program,
and extrapolation suggests that $A_8= 0.778$ could be correct to three significant
digits.  The flat metric area of the octagon is $\bar A_8 = 2 \tan (\pi/8)\simeq 0.82843.$
This means that the extremal area is about 6.1\% lower than the area of the flat-metric.
This is a larger decrease compared to the case of the hexagon.  For the 
perimeter, we find that the
relation ${\cal P}_8= 4 A_8$ holds in the dual program with an accuracy of 3\%.
Finally the relation $A_8 = 4 \nu$ holds very accurately.

 Let us now consider curvature.  For this we can first consider the
 region $U_4$ which is topologically a disk.  To figure out the extent
 of the various $U_k$ on the surface one must draw all boundary geodesics
 for the four bands (a total of 8).  The boundary
 of $U_8$ is given by eight geodesic arcs.  Each arc corresponds to a 
 different boundary geodesic.  
 Figure~\ref{fig:LastGeo8Gon}  
 makes it very
 plausible  that
 the arcs join smoothly without corners.   That being the case, the boundary
 of $U_4$ is a single geodesic and the Gauss-Bonnet formula implies that
 \be
{1\over 2\pi}  \int_{U_4} K \, dA =  1 \,. 
 \ee
The `dome' $U_4$ carries the curvature of a half sphere.   
The union $U_4 \cup U_3$ is also a topological disk. Its boundary
is the outer boundary of $U_3$ and it is clear that it is also bounded
by a smooth geodesic (actually built of sixteen geodesics joining smoothly).
Since it is also a disk and we know the bulk integral of curvature over $U_4$
it follows from Gauss-Bonnet applied to the disk $U_4 \cup U_3$ that  
\be
{1\over 2\pi}  \int_{U_3} K \, dA  + {1\over 2\pi} \int_{\partial U_4} k ds  = 0 \,. 
\ee
Here the second integral is that of the line curvature possibly located on
the boundary of $U_4$.  While we cannot tell the values, it would seem
sensible to have negative line curvature on the boundary (as in the hexagon
and the Swiss-cross) and positive bulk curvature on $U_3$, the two of them
canceling exactly.  This is of course also needed when we think of extending
the disk to the whole surface.  As we saw for the hexagon, the inclusion of
$U_2$ does not contribute to the Euler number, for the line curvature
on the boundary of $U_2$ cancels with the corner contribution at the
corresponding vertex of the polygon.

\subsection{The decagon briefly noted}

\begin{figure}[!ht]
\leavevmode
\begin{center}
\epsfysize=14.0cm
\epsfbox{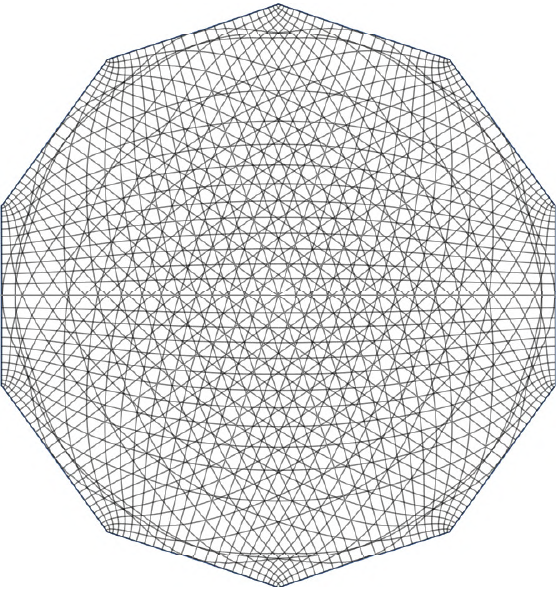}
\end{center}
\caption{\small  The saturating geodesics for $P_{10}$. The surface
is covered by a connected region with five bands of geodesics, ten regions with four bands, ten regions with three bands and ten regions with two bands.}
\label{fig:10Gon}
\end{figure}

\begin{figure}[!ht]
\leavevmode
\begin{center}
\epsfysize=13.5cm
\epsfbox{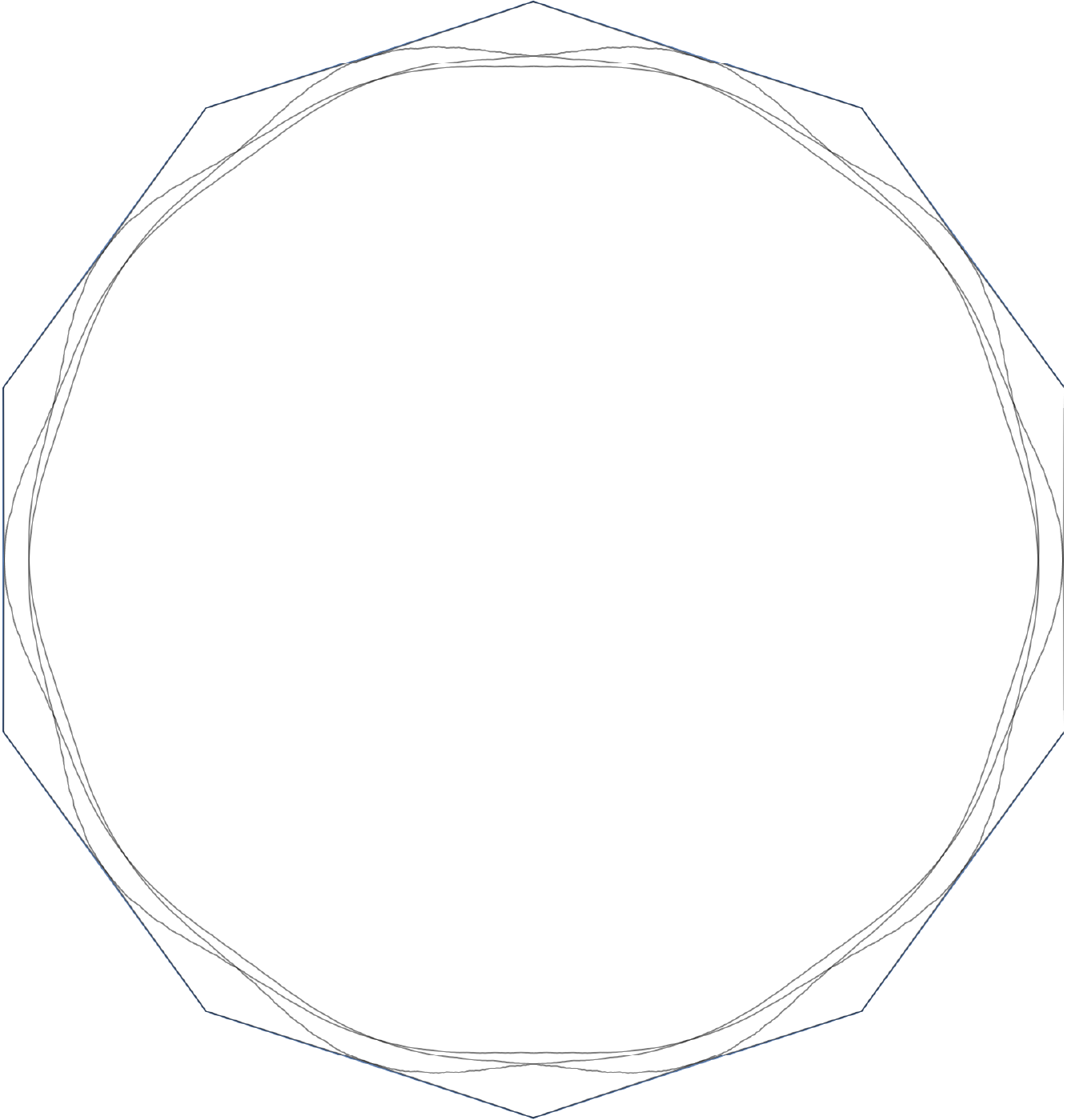}
\end{center}
\caption{\small  The boundary saturating geodesics for $P_{10}$.  The surface
is covered by a connected region $U_5$, ten regions comprising $U_4$, ten regions comprising $U_3$  and ten regions comprising $U_2$.}
\label{fig:LastGeo10Gon}
\end{figure}

We can easily apply our program to find the extremal metric
on the decagon. The picture of saturating geodesics constructed from our solution is given in Figure~\ref{fig:10Gon}. In Figure~\ref{fig:LastGeo10Gon}  we draw the boundary geodesics in classes $C_1,C_2,C_3,C_4$, and $C_5$.
The region $U_5$ comprising most of $P_8$ is
covered by five geodesic bands  
and has positive curvature. 
As in the previous cases, the boundary of $U_5$, clearly visible in
Figure~\ref{fig:LastGeo10Gon},  is composed of geodesics
segments that join without corners. Therefore,
 \be
{1\over 2\pi}  \int_{U_5} K \, dA =  1 \,. 
 \ee
 There are ten small regions of four geodesics. Near each of the edges there is a small region with three geodesics which is bigger than the region with four geodesics. Near each corner we have a flat region with two geodesics. 

 The best results for area both in the primal program
 and in the dual program are obtained with $N_c=64$.  We have
 \be
 0.7444<A_{10}<0.7500.
 \ee
As in the case of the hexagon and octagon, convergence is faster in the dual program and extrapolation suggests that $A_{10}=0.746$. The flat metric area of the decagon is $\bar{A}_{10}
=\tfrac52\tan\bkt{\pi/10}\simeq 0.812299
$.
This means that the extremal area is about $8.2\%$ lower than the area of the flat-metric, which is a larger decrease as compared to both the hexagon and the octagon. We also find that the relations ${\cal P}_{10}=4 A_{10}$ and $A_{10}= 5 \nu$ hold very accurately.

\section{Higher polygons and the extremal metric on $\mathbb{RP}_2$} \label{pus_metric}

One can naturally ask for the limit of the extremal metric on $P_{2n}$ as $n$ approaches infinity.  In this section we
review results of Pu~\cite{pu} and Ahlfors~\cite{ahlfors} that 
consider, respectively, the Riemannian and conformal versions
of a minimal area problem on a round disk with a length constraint on
all curves joining opposite points on the boundary.  We study 
the extremal 
metrics on $P_{2n}$ for large $n$ and conclude that
away from the boundary the metrics seem to converge to that of 
the disk problem. 

\subsection{Riemannian and conformal extremal metrics
on $\mathbb{RP}_2$}

Among the few minimal area metrics that are explicitly known, 
that for the real-projective surface $\mathbb{RP}_2$ in two-dimensions stands out. 
In this manifold there are non-contractible closed curves and 
the Riemannian version of the  minimal 
area  problem is well defined.
The optimal {\em Riemannian} metric was determined by Pu~\cite{pu}.  The surface
with this metric can be described as a round hemisphere with the standard
round metric and antipodal points on the circular boundary identified.
The problem has also been posed and discussed in the conformal
setting~\cite{ahlfors}.  
In this version one can consider the disk $|z | \leq \tfrac{1}{2}$,
including its boundary, 
and ask for the minimal area conformal metric such that any curve joining
an arbitrary point $z$ on the boundary to its opposite $-z$ be longer than
or equal to one.  The minimal area (conformal) metric is the same as above
Riemannian metric;  a half-sphere
with the round metric.  Joining any two antipodal points on the boundary
there are infinite systolic geodesics, the set of all great half-circles joining the points.
In fact these geodesics cover the complete hemisphere, and we can think of them
as the band of geodesics associated with the pair of points.  Furthermore, if we consider
the bands of geodesics for
{\em all pairs} of opposite points, we conclude that the surface is in fact covered by an 
infinite number of bands of geodesics.  

With the systole length equal to 1, the radius $R$ of the
half sphere extremal metric  is such that the systole $\pi R$ is equal to one, so
\be
R = {1\over \pi}\, . 
\ee
The area $A$ of the surface in the minimal area metric is
\be
\label{aminpu}
A =  2\pi R^2 =  {2\over \pi } \simeq  0.63662\,. 
\ee
The systolic area $\sigma(g) $ for an arbitrary metric $g$ is  defined as the area $A(g)$
of this metric divided
by the square of the systole $\ell_s(g)$ in this metric:
\be
\sigma(g)  \equiv  {A(g)\over \ell^2_s(g)} \,.
\ee
The extremal systolic area is the minimal area when the systole is chosen
to be equal to one.  
Pu's result for $\mathbb{RP}_2$ implies that for any metric $g$ on this surface
\be
\sigma(g)   \geq {2\over \pi} \,.
\ee

As $n$ becomes large the minimal area problem on $P_{2n}$ seems
to turn into the $\mathbb{RP}_2$ minimal area problem.  
Indeed,  with apothem fixed at $\tfrac{1}{2}$, 
 as $n\to \infty$ the polygon approaches a circular disk of radius $\tfrac{1}{2}$.
Moreover, as $n \to \infty$ the length condition on curves joining opposite edges seems
to become    
a length condition 
on curves  joining opposite points.

We will show in this section that this intuition is partially correct.  We find
evidence that the metric on the polygons approaches that of the hemisphere
away from the boundary and the area approaches the expected area of a 
round metric.  We find, however, that the perimeter of the polygon
does not approach the perimeter of the hemisphere on the round metric.
The perimeter in the polygonal case is  larger than the perimeter 
in the round metric.  This is presumably due to the local geometry at the
(now infinite number of) vertices, where the edges meet orthogonally.
Recall the polygon $P_{2n}$ is covered by regions 
$U_2,  U_3, \cdots, U_n$ and all except $U_n$ become smaller
as $n\to \infty$.  It seems plausible that the part of the surface 
that is the complement of $U_n$ 
has vanishing area as $n\to \infty$.  While
the `serrated' boundary of the polygon is larger than the perimeter in the
round hemisphere, it is also plausible that the length of the boundary of $U_n$ approaches 
the perimeter of the round hemisphere as $n\to \infty$.

In order to make comparisons let us consider the metric on a round 
half-sphere.   In unit-free coordinates $z$ a sphere of  radius $R$
is described by the metric over the full complex plane 
\be
\label{eord}
ds^2 =  {4R^2|dz|^2  \over (1 + |z|^2  )^2 } \,,   \ \ \ 
\hbox{radius $R$ sphere with equator} \ |z|= 1\,.
\ee
With the north pole at $z=0$ and the south pole at $z=\infty$ one
can quickly check that the equator is the curve $|z|=1$ and half
of the sphere (the northern hemisphere) is the region $|z|\leq 1$.
Since our polygons $P_{2n}$ of apothem $\tfrac{1}{2}$ converge into the
$|z| \leq \tfrac{1}{2}$ region as $n\to \infty$, we need to let $z \to 2 z$ in
the above metric to find that
\be
\label{eordnew}
ds^2 =  {16R^2|dz|^2  \over (1 + 4|z|^2  )^2 } \,, \ \ \   \hbox{radius $R$ sphere with equator} \ |z|= \tfrac{1}{2}\,.
\ee
In the conventional notation $ds^2 = \rho^2 |dz|^2$ we read
\be
\rho^2 (z) =  {16 R^2 \over (1+ 4|z|^2)^2} \,.
\ee
For the extremal metric with systole equal to one, we saw above that $R = 1/\pi$.
Therefore,
the extremal metric on the disk $|z|\leq \tfrac{1}{2}$ is 
\be
\label{rho-disk-metric}
\rho^2 (z) =  {16 \over \pi^2} {1  \over (1+ 4|z|^2)^2} \,.
\ee
In particular, we have  that at the origin
\be
\rho^2 (0)  = {16\over \pi^2} \simeq 1.62114\, . 
\ee

\subsection{Area, perimeter, and metric 
for large $n$ 
 polygons}

We can now use the primal and dual programs to determine 
the minimal area metrics for the following $2n$-gons:
\be
2n =  8 , \ 16, \ 32 , \ 64, \ 128 \,.  
\ee
In each of the polygons
we will calculate the relevant quantities for lattice resolution
$N_c = 4, 8, 16, 32$.
Our Mathematica program for the primal succeeded in giving us answers
for all these cases.  On the other hand the dual program had more trouble
and only gave us partial data.

The theoretical prediction  (\ref{avp-exact}) 
of the dual program relating area and perimeter 
for any $2n$-gon,  
\be
{\cal P}_{2n} =  4 A_{2n}  \,, 
\ee
is supported by the data.  
On the other hand 
for the extremal metric in $\mathbb{RP}_2$ (with systole one) we have
$A = 2/\pi$ and perimeter ${\cal P} = 2\pi R = 2$, giving 
\be
{\cal P}  = \pi A \,, \quad  \hbox{on the extremal metric on } \ \mathbb{RP}_2 \,. 
\ee
The perimeter ${\cal P}_{2n}$  
 as $n\to\infty$ is larger than the one expected in $\mathbb{RP}_2$.
We believe that the explanation for this is the ``serrated" form of the boundary
in which, the flat region about the corners of the polygon (in the 
non-singular 
 metric presentation) has the sides meeting
orthogonally.  We do not have a detailed 
model of the polygonal boundary that 
can explain how to obtain agreement with the area /perimeter relation
of the extremal $\mathbb{RP}_2$.  

\begin{figure}[!ht]
\leavevmode
\begin{center}
\epsfysize=7cm
\epsfbox{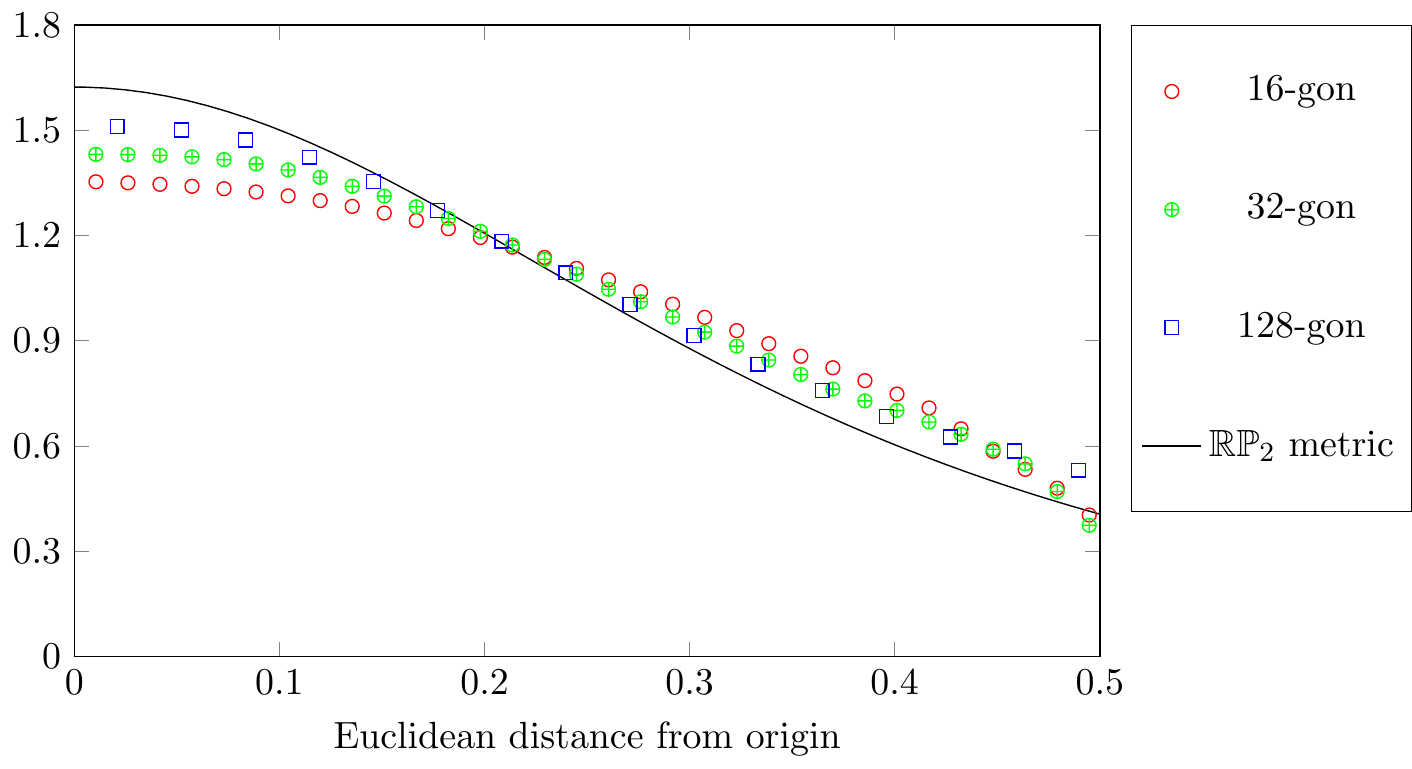}
\end{center}
\caption{\small  
Extremal metric along the horizontal axis for various polygons and comparison with the $\mathbb{RP}_2$ extremal metric. We see see evidence of convergence
to the $\mathbb{RP}_2$ metric.}
\label{fig:Pu_Polygons_hor}
\end{figure}

Modulo the subtlety associated with $U_2$ and the corners we find good
 evidence that the metric converges to the extremal metric on $\mathbb{RP}_2$
 away from the boundary. In Figure~\ref{fig:Pu_Polygons_hor} we plot the metric along the base of $T_{2n}$ for various polygons and compare it with 
 the $\mathbb{RP}_2$ result. It is clear that the metric for higher polygons is converging to the $\mathbb{RP}_2$ extremal metric.

\begin{figure}[!ht]
\leavevmode
\begin{center}
\epsfysize=7cm
\epsfbox{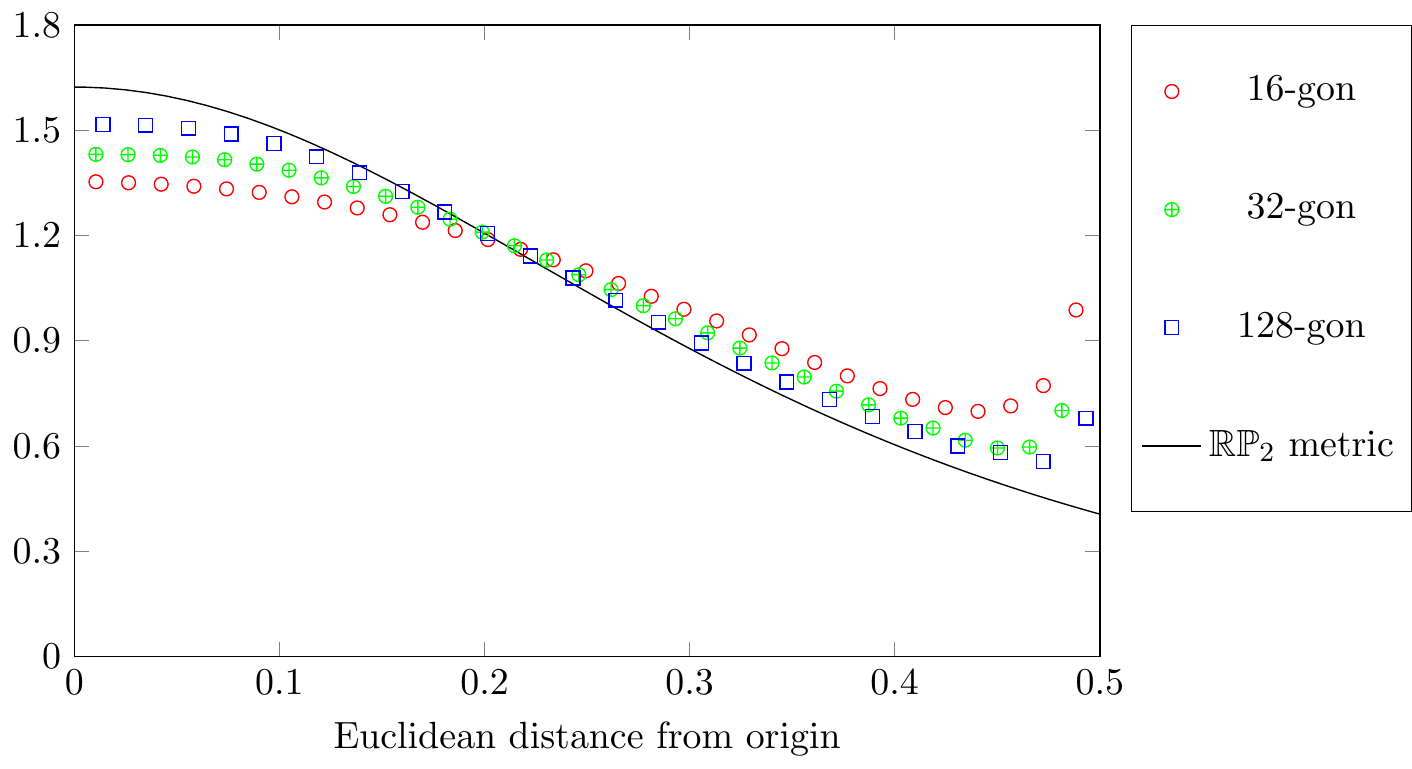}
\end{center}
\caption{\small  
Extremal metric along the hypothenuse of $T_{2n}$ for polygons $P_{2n}$ and comparison with the extremal metric on $\mathbb{RP}_2$. We see evidence 
of convergence away from the vertex of the polygons.
}
\label{fig:Pu_Polygons_dia}
\end{figure}

In Figure~\ref{fig:Pu_Polygons_dia} we plot 
the extremal metric for various polygons along the hypotenuse of 
$T_{2n}$   
and compare it to the extremal
$\mathbb{RP}_2$ metric.  We only plot to a fiducial distance of $0.5$  from the origin. 
As we approach the vertex of the polygon 
the metric starts increasing and then diverges at the vertex. For higher polygons the onset of this divergence  shifts closer 
to the vertex. In the $n\to \infty$ limit we expect the metric for polygons to approach the $\mathbb{RP}_2$ metric everywhere except at the vertex where we encounter a singularity.

For the area, extrapolation from the primal gives quantities $A_{2n}$ that are
\be
A_8 =  0.778\,,  \quad  A_{16} =  0.700\,, \quad  A_{32} =0.677 \,,  \quad A_{64} =   0.66\, ,  \quad A_{128} =   0.654\, .
\ee  
An extrapolation of this suggests $A_{\infty} \simeq  0.65$ not far from the minimum 
in (\ref{aminpu}): $\tfrac{2}{\pi} \simeq 0.63662$.  This is evidence that the area
is converging to the expected value in the correspondence.

\section{Conformal metric on $\mathbb{RP}_2$ from the dual program}
\label{con_met_on_rp2}

In this section we demonstrate that the  
extremal conformal or Riemannian   metric on the real projective plane is a  
critical point of  dual variational principle
of~\cite{Headrick:2018ncs}.  This demonstration will entail
several challenges, as we will have to deal with a situation
in which we have an infinite number of bands of geodesics.
As we recall, any two antipodal points determine a band of
geodesics that in fact covers the full surface.  

\subsection{Setting up the calculation}  

We will use the
azimuthal angle $\phi_0$ of a point on the boundary
to label the band of geodesics that begins   
at this point and ends on the antipode $\phi_0 + \pi$.  Each
geodesic on this band is called a $\phi_0$-geodesic.  
 For each band
we will have  coordinates $x_{\phi_0}$ and $\varphi_{\phi_0}$, where
the latter is the  
function that appears in the dual program.
Note that $x_{\phi_0}$ and $\varphi_{\phi_0}$
 are functions all over the sphere, and the geodesics 
 must be lines of
constant $\varphi_{\phi_0}$.  

Our challenge is to find these functions.  The function $x_{\phi_0}$ 
is the length parameter along $\phi_0$-geodesics and is simple to 
determine.
The 
function $\varphi_{\phi_0}$ is significantly harder to find
because it must satisfy an additional constraint in addition
to being constant along geodesics.   
To explain this, consider
the metric using the coordinates $(x_{\phi_0}, \varphi_{\phi_0})$, that 
as explained in~\cite{Headrick:2018ncs} takes the form:
\be
\label{rp2band_metric}
ds^2 \ = \ (dx_{\phi_0})^2  + {1\over  h^2_{\phi_0} } (d \varphi_{\phi_0} )^2 \,.  
\ee
It is not too hard to find a  particular function $\tilde\varphi_{\phi_0}$ 
that, as required, is a constant along the geodesics, as well as the 
associated $\tilde h_{\phi_0}$ making the above $ds^2$ the round metric on
the hemisphere.  The correct
$\varphi_{\phi_0}$, however, will be a reparameterization of $\tilde \varphi_{\phi_0}$. 
As emphasized in \cite{Headrick:2018ncs} (eqn.\,(7.10)) the 
various bands $C_\alpha$ that are active at {\em any} point on the surface 
must satisfy the sum rule
\be
\sum_\alpha |d\varphi^\alpha | \ = \ 1 \,.
\ee
Here $\alpha$ labels the active bands at any point.  In our case
the role of the index $\alpha$ is played by the continuous variable 
$\phi_0\in [0, \pi]$.  Still using the sum notation, we have
\be
\sum_{\phi_0} |d\varphi_{\phi_0} | \ = \ 1 \,.
\ee
The form (\ref{rp2band_metric}) of the metric implies
 that $|d\varphi_{\phi_0} |= |h_{\phi_0}|$
and therefore the sum rule becomes
\be
\label{srvmvg}
\sum_{\phi_0} |h_{\phi_0} | \ = \ 1 \,.
\ee
The left hand side is a function over the surface, but that function must
be a constant.    
Under a reparameterization $\varphi_{\phi_0} \to \tilde \varphi_{\phi_0} (\varphi_{\phi_0})$ the metric invariance requires
\be
\label{h-repara}
{|d\varphi_{\phi_0}|\over |h_{\phi_0}|} = {|d\tilde \varphi_{\phi_0}|\over |\tilde h_{\phi_0}|}  
\quad \to \quad 
 |h_{\phi_0}| =  |\tilde h_{\phi_0}|\, \Bigl| {d\varphi_{\phi_0}
\over  d\tilde\varphi_{\phi_0}} \Bigr| \,. 
\ee
Since the function $h_{\phi_0}$ transforms nontrivially under a $\varphi_{\phi_0}$
reparameterization, it is clear that if an original choice of $\varphi_{\phi_0}$ does
not not make the sum rule (\ref{srvmvg}) work, a reparameterization of it may.

To find the continuous form of the sum rule we imagine working with a polygon with $2n$ sides, and thus $n$ foliations, in the limit of large $n$.  We will take
$n$ equally spaced values of $\phi_0$  in the range from $0$ to $\pi$. 
With this discretization the sum rule becomes
\be
\sum_{k=1}^n | h_{\phi_k} | =   1 \, , \quad \hbox{with} \ \ \ \ \phi_k = k\,  {\pi\over n} \,,  \  \ k = 1 , \ldots , n\,.  
\ee
Passing to an integral is done by multiplying the above equation by $\Delta \phi = {\pi\over n}$
\be
\sum_{k=1}^n  |h_{\phi_k} | \,  \Delta \phi  
\,  = \,   {\pi\over n}  \,. 
\ee
For very large $n$, using a continuous 
$\phi_0\in (0, \pi)$,  we have the integral constraint: 
\be
\label{integral_constraint}
\int_0^\pi  d \phi_0\,  |h_{\phi_0}|  \ \simeq \  {\pi\over n} \,, \qquad  n \gg 1 \,. 
\ee
This is our main constraint.  The most nontrivial part is ensuring that
the integral over all the bands of the position dependent $|h_{\phi_0}|$ is
a constant over the sphere.  After that is done, there is the issue of the 
constant being equal to the constant on the right-hand side.  This will
require finding the value of the height function $\nu$ as a function of $n$,
for large $n$.  This value of $\nu$ must 
be such that the correct area
of the surface is obtained from the large $n$ limit of $n\nu$ as explained
in (\ref{annuval}).

\subsection{Working out the details}

In this section $z = x+ i y$ is the complex coordinate in which
the disk is $|z| \leq {1\over 2}$.  With radius $R = 1/\pi$, corresponding
to a systole of 1, the metric (\ref{rho-disk-metric})~is:  
 \be
\label{shp-met-eqn}
ds^2 =  {16\over \pi^2}  {|dz|^2 \over (1+ 4 |z|^2)^2}\,,   \quad |z| \leq \tfrac{1}{2} \,, \quad 
z = x + i y \,. 
\ee
To use some of the previously derived results, we must relate these $(x,y)$ coordinates
to the cartesian coordinates $\hat x, \hat y$ 
defined as usual from the spherical angles $\theta, \phi$ and the radius $R$:
\be
\begin{split}
\hat x =\ &   R \sin\theta \cos \phi \,,\\
\hat y  =\ &   R \sin\theta \sin \phi\,. \\
\end{split}
\ee
One can easily show that for $R= 1/\pi$, our case of interest, the two sets of coordinates
are related as follows:
\be
\begin{split}
\hat x =\    {1\over \pi} \,  {4 x \over 1 + 4 |z|^2 }  \,, \ \ \ \ \ 
\hat y =\   {1\over \pi} \,  {4 y \over 1 + 4 |z|^2} \,. \\
\end{split}
\ee
One can also prove the useful identities:
\be
1 + 4 ( x^2 +  y^2)  =  \sec^2 \tfrac{\theta}{2} \,,  \qquad  
\cos \theta = {1- 4 |z|^2\over 1+ 4|z|^2} \,. 
\ee

\begin{figure}[!ht]
\leavevmode
\begin{center}
\epsfysize=5.5cm
\epsfbox{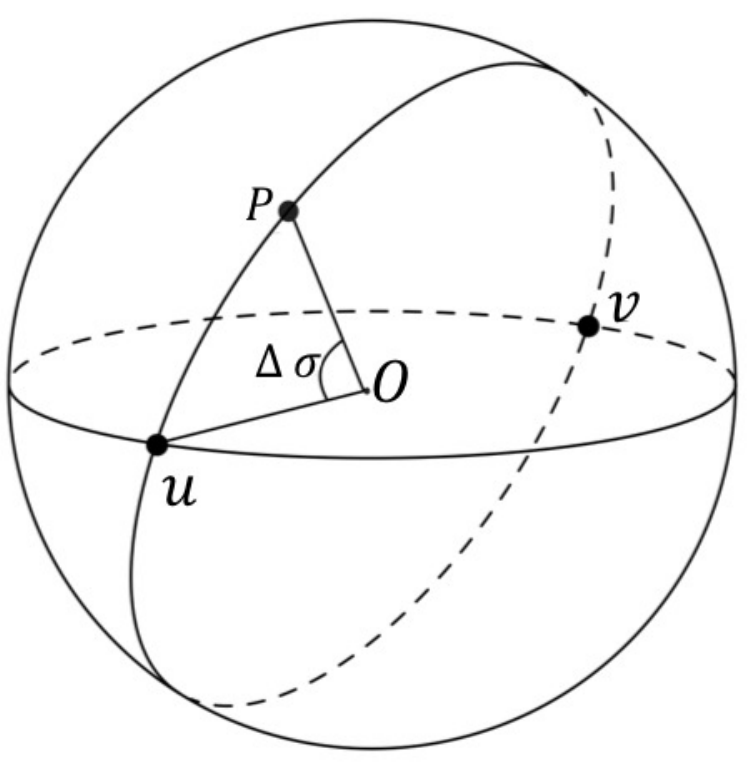}
\end{center}
\caption{\small  
A geodesic joining antipodal points $u$ and $v$ on the equator. 
}
\label{fig:gCircle}
\end{figure}

\medskip
\noindent
\underline{Defining coordinates for the band labelled by $\phi_0$.} 
\nobreak
On the round hemisphere of radius $1/\pi$ representing the extremal
metric with systole one, 
 we can consider the so-called $\phi_0$-geodesics that begin on the equator at 
 $\phi_0$ and
 end on the equator at~$\phi_0 + \pi$.  
They comprise a band of geodesics that covers the full hemisphere. 
 Using coordinates $\theta$ and $\phi$ on the sphere,
$x_{\phi_0}$ is a coordinate that parameterizes these geodesics by length. In Figure~\ref{fig:gCircle}, the geodesic starts at point $u$, corresponding to $\phi_0$ and ends at $v$, corresponding to $\phi_0 + \pi$.  
The point $P$ lies on a geodesic joining $u$ and $v$ and has spherical coordinates $\bkt{\phi,\theta}$. The length of the geodesic arc connecting $u$ and $P$ is then equal to $x_{\phi_0}$ at $P$.
This arc length can be expressed in terms of the radius $1/\pi$ of the sphere and the angle $\Delta\sigma$ shown in the figure:
\be
x_{\phi_0}=\Delta \sigma/\pi\, .
\ee
Let us denote by $P$ and $u$ the vectors from origin to points $P$ and $u$ respectively 
\be
P\ =\tfrac{1}{\pi}\bkt{\sin\theta\cos\phi,\sin\theta\sin\phi,\cos\theta},\qquad
 u\ =\tfrac{1}{\pi}\bkt{\cos\phi_0,\sin\phi_0,0}.
\ee
The angle $\Delta\sigma$ is then given by
\be
\cos\Delta\sigma= \frac{P\cdot u}{|P||u|}=\sin\theta\cos\bkt{\phi-\phi_0}\, .
\ee
The coordinate $x_{\phi_0} (\theta, \phi)$ is just the arc-length: 
\be
\label{Xphidef}
x_{\phi_0} (\theta, \phi) =  \tfrac{1}{\pi}  \cos^{-1} \bigl( \sin \theta \cos (\phi-\phi_0) \bigr) \,.
\ee
 We thus have
\be
\label{Xphidefd}
\cos \pi x_{\phi_0} (\theta, \phi) =   \bigl( \sin \theta \cos (\phi-\phi_0) \bigr) \,.
\ee
In terms of the disk coordinates $(x, y)$, this is
\be
\label{pi-x-phi-def}
\cos \pi x_{\phi_0}   \ = \ {4 ( x \cos\phi_0 +  y \sin\phi_0) \over 1 + 4 |z|^2  } \,. 
\ee
Let us find the second coordinate $\tilde \varphi_{\phi_0}$.
This coordinate must be such that each $\phi_0$ geodesic is a curve
 of constant $\tilde\varphi_{\phi_0}$.
Geodesics on the sphere are great circles
which can be parameterized by a constant $c$ and an angle $\phi_0$ as follows
\be
\cot\theta=c\, \sin\bkt{\phi-\phi_0}\, .
\ee
At the equator we have latitude $\theta = \tfrac{\pi}{2}$, 
the left-hand side vanishes  and
 the right-hand side implies that  the geodesic meets the equator at the azimuthal angles $\phi=\phi_0$ and $\phi_0+\pi$.  
Now consider following the geodesic as $\phi$ increases from 
the value $\phi_0$.  The latitude angle $\theta$ decreases from $\pi/2$ and
reaches a minimum when the right-hand side is a maximum at 
$\phi= \phi_0 + \tfrac{\pi}{2}$.  This minimum latitude will be the
searched for $\tilde\varphi_{\phi_0}$.  This implies 
$\cot \tilde\varphi_{\phi_0} = c$ and
the equation of the geodesic is 
\be
\label{max-circle-sp-vm}
\cot\theta=\cot \tilde\varphi_{\phi_0}\, \sin\bkt{\phi-\phi_0}\, .
\ee
For our problem we only consider geodesics on the upper hemisphere.
For a chosen $\phi_0$ these fall into two classes, separated by the
great circle going through the north pole:  those that are lifts of the half-circle
$\phi\in [\phi_0, \phi_0+ \pi]$ and those that are lifts of the half circle $\phi\in[\phi_0+ \pi, \phi_0+ 2\pi]$.   Since the left-hand side of (\ref{max-circle-sp-vm})
must always be positive to be on the upper hemisphere, 
$\cot\tilde\varphi_{\phi_0}$
must alternate sign in the two classes and will range in the interval
\be
\label{range-tilde-varphi}
\tilde \varphi_{\phi_0} \in \ \bigl(-\tfrac{\pi}{2} , \tfrac{\pi}{2}\bigr )\, .
\ee
For our conventions it will be useful to alter the sign in (\ref{max-circle-sp-vm}) and to have  
\be
\label{max-circle-sp-vmbb}
\cot\theta=-\cot \tilde\varphi_{\phi_0} \, \sin\bkt{\phi-\phi_0}\, .
\ee
This amounts to just changing the role of the two classes of curves.  From this
it follows that 
\be
\label{max-circle-sp-vmbb}
\tan\tilde\varphi_{\phi_0}=-\tan \theta\, \sin\bkt{\phi-\phi_0}\, .
\ee
This formula displays $\tilde \varphi_{\phi_0}$ as a function that is 
constant over each geodesic half circle.  
One can quickly show that in the disk coordinates $(x,y)$ this becomes
\be\label{eq:varphiDisk}
\tilde \varphi_{\phi_0}  = \,-  \tan^{-1} \biggl( {  4 (  y \cos \phi_0 -  x \sin \phi_0) \over
1 - 4 |z|^2 } \biggr) \,. 
\ee

\noindent
\underline{Calculating the metric in the $x_{\phi_0}$ and $\tilde \varphi_{\phi_0}$ coordinates.} 

We have now expressions for 
$x_{\phi_0}$ and $\tilde\varphi_{\phi_0}$.
To find $dx_{\phi_0}$ we take the differential of  equation (\ref{pi-x-phi-def}). 
A little computation gives
\be
- \tfrac{\pi}{ 4} \sin (\pi x_{\phi_0} )  \, d x_{\phi_0}  \ = \ 
{\alpha  \, d x + \beta \, d y\over (1 + 4|z|^2)^2} \,.
\ee
Here we have defined the functions $\alpha$ and $\beta$ as follows:
\be
\begin{split}
\alpha =\ &   ( 1+ 4( y^2 -  x^2)) \cos\phi_0 - 8  x y \sin\phi_0\,, \\
\beta =\ &   ( 1+ 4( y^2 -  x^2)) \sin\phi_0 - 8  x y \cos\phi_0\,.
\end{split}
\ee
A quantity that will appear many times is $\alpha^2 + \beta^2$, which can be written in four useful ways:
\be
\label{aux_eqn}
\begin{split}
\alpha^2 + \beta^2 \ = &  \ 1 + 16 |z|^4 - 8 \bigl[ ( x^2 -  y^2 ) \cos 2\phi_0 
+ 2 x  y \sin 2\phi_0\bigr]\,, \\
\ = &  \ (1 - 4 |z|^2)^2  + 16  (  y \cos \phi_0 -  x \sin\phi_0)^2\,,  \\
\ = &  \ (1 + 4 |z|^2)^2  - 16  (  x \cos \phi_0 +  y \sin\phi_0)^2\,, \\
\ = & \  (e^{i\phi_0} - 4 z^2 e^{-i\phi_0} ) (e^{-i\phi_0} - 4 \bar z^2 e^{i\phi_0} ) \,, \end{split}
\ee
with 
\be
e^{i\phi_0} - 4 z^2 e^{-i\phi_0}  = \alpha + i \beta \,. 
\ee
Calculation of $d\tilde \varphi_{\phi_0}$ from (\ref{eq:varphiDisk}) 
takes a little work, and with the
use of (\ref{aux_eqn}) gives
\be
d\tilde \varphi_{\phi_0}  =  {4(\beta d x - \alpha d y) \over \alpha^2 + \beta^2 } \,. 
\ee
The metric function $\tilde h_{\phi_0}$ can be obtained from the general relation
between calibrations and $\varphi$ functions~\cite{Headrick:2018ncs}:
\be
u^\alpha = -  {{}^* d \varphi^\alpha \over |d\varphi^\alpha|}\,.
\ee
For our case this reads:
\be
dx_{\phi_0}  = -  {{}^* d \tilde\varphi_{\phi_0}  \over |d\tilde\varphi_{\phi_0}|}\,. 
\ee
Since $|d\tilde\varphi_{\phi_0}| =| \tilde h_{\phi_0}|$ when we write the metric as
\be
\label{lgbbb}
ds^2 \ = \ (dx_{\phi_0})^2  + {1\over \tilde h^2_{\phi_0} } (d \tilde\varphi_{\phi_0} )^2 \,,  
\ee
we find
\be
{}^* d \tilde \varphi_{\phi_0}  \ = \ - |\tilde h_{\phi_0} | \, dx_{\phi_0} \,.  
\ee
Using the earlier calculation of $d \tilde \varphi_{\phi_0}$ and $d x_{\phi_0}$ 
and recalling that ${}^* ( a dx + b dy) =  - b dx + a dy$, we get
\be
|h_{\phi_0} |= \,   \pi  \sin ( \pi x_{\phi_0})  { (1 + 4 |z|^2)\over \alpha^2 + \beta^2 }\,.
\ee
This can be simplified.  Indeed, a short calculation shows that 
\be
 { (1 + 4 |z|^2)\over \alpha^2 + \beta^2 } = {1\over (\sin  \pi x_{\phi_0}) ^2} \,,  
\ee
giving
\be
\label{tilde-h-first}
|\tilde h_{\phi_0} |= \,   {\pi \over  \sin  \pi x_{\phi_0}}  \,. 
\ee
The metric (\ref{lgbbb}) 
 is therefore given by
\be
ds^2 =   (dx_{\phi_0})^2  + {1\over \pi^2 } (\sin \pi x_{\phi_0})^2\,   (d \tilde\varphi_{\phi_0} )^2 \,.
\ee
Thinking of the origin of the $\phi_0$ geodesics as a new north pole, we can
see that  $\pi x_{\phi_0}\in [0, \pi] $ is in fact the relevant
polar angle $\tilde\theta_{\phi_0}$.  With this the metric takes the form
\be
ds^2 = R^2 \bigl[ (d\tilde\theta_{\phi_0})^2  + 
\sin^2 \tilde\theta_{\phi_0}\,   (d \tilde\varphi_{\phi_0} )^2 \bigr]  \,,  \ \ \ R = 1/\pi\,.
\ee
This is clearly the expected round metric on the hemisphere. 
While we had
defined $\tilde \varphi_{\phi_0}$  as the minimum latitude of the geodesic,
 in this picture $\tilde \varphi_{\phi_0}$ is seen equivalently 
 to be an azimuthal angle.
The function $\tilde \varphi_{\phi_0}$ that labels the various $\phi_0$
geodesics is the departure angle at the origin $x_{\phi_0}=0$.  This is clearly
a natural variable and its range $(-\tfrac{\pi}{2} , \tfrac{\pi}{2})$ (see~(\ref{range-tilde-varphi}))
is consistent with this interpretation. Here $\tilde \varphi_{\phi_0}=0$
labels the geodesic departing orthogonal to the boundary.

\medskip
\noindent
\underline{Reparameterization and final steps.}

A little experimentation shows that  
the function $|\tilde h_{\phi_0}|$ does not satisfy the sum rule
constraint (\ref{integral_constraint}): the integral on the left-hand side fails
to be position independent on the surface.  As explained before,
the  definitive coordinate $\varphi$ is to be obtained 
by a reparameterization of $\tilde \varphi$.  A bit of trial and
error quickly led to a satisfactory solution: 
\be
\varphi_{\phi_0} =  \tfrac{\nu}{2}  \sin \tilde \varphi_{\phi_0}\,,   \quad  \varphi_{\phi_0} \in 
\bigl( - \tfrac{\nu}{2} , \tfrac{\nu}{2} \bigr) \,.
\ee
By setting the coefficient equal to $\nu/2$ we guarantee that 
 function has the required discontinuity by $\nu$: the top-most and bottom-most
 geodesic in the band are identified, and on those the function $\varphi$ takes
 values $\nu/2$ and $-\nu/2$).  Using (\ref{h-repara}) we have 
 \be
\label{h-reparam}
|h_{\phi_0}| =  |\tilde h_{\phi_0}|\, \Bigl| {d\varphi_{\phi_0}
\over  d\tilde\varphi_{\phi_0}} \Bigr|  =  \frac{\pi \nu}{2} \, { \cos \tilde \varphi_{\phi_0}\over \sin \pi x_{\phi_0}} \,.  
\ee
This result simplifies considerably when referred to $(x,y)$ coordinates.
All square roots disappear and we find 
\be\label{eq:|hphi|}
|h_{\phi_0} |= \frac{\pi \nu}{2} \cdot { 1 - 16 |z|^4 \over 1 + 16 |z|^4} 
 \cdot  {1\over 1- a \cos 2\phi_0 - b \sin 2\phi_0}  \,, 
\ee
with:
\be\label{eq:abconstraints}
a\equiv  {8 (x^2 - y^2) \over 1+ 16 |z|^4}\,,\qquad  \ b \equiv  { 16 x y \over  1+ 16 |z|^4}\,, 
\ee
With 
this value of $|h_{\phi_0}|$   the constraint (\ref{integral_constraint}) becomes
\be
\label{lgbbspr}
\frac{\pi \nu}{2} \cdot{ 1 - 16 |z|^4 \over 1 + 16 |z|^4} 
\int_0^\pi    \,  {d \phi_0\over 1- a \cos 2\phi_0 - b \sin 2\phi_0}  \simeq  {\pi\over n} \,. 
\ee
The left hand side must be $z$ independent for this to work.  The integral
is readily evaluated giving: 
\be\label{eq:confint}
\int_0^\pi {} {d \phi_0\over 1- a \cos 2\phi_0 - b \sin 2\phi_0} = {\pi \over \sqrt{1 - (a^2 + b^2 )} } = \pi \cdot { 1 + 16 |z|^4 \over 1 - 16 |z|^4}  \,. 
\ee
As a result, the $z$ dependence of the left-hand side of (\ref{lgbbspr}) 
disappears and we are
left with 
\be
\frac{\pi \nu}{2} \pi  \simeq  {\pi\over n} \quad \to \quad  \nu = {2\over \pi n} \,. 
\ee
With this choice of $\nu$ we satisfy the constraint for large $n$. 
Recall that the area functional at the critical point is $A = \sum_{k=1}^n  \nu_{k} \ell_s$.
By symmetry, all $\nu_k$ parameters are to $\nu$, and having set 
$\ell_s =1$,  we have
\be
A = n \nu = n \cdot {2\over \pi n}  \ = \  {2\over \pi}\,.
\ee
This is the correct value for the extremal 
area $A = 2 \pi R^2 = 2 \pi/ \pi^2 = 2/\pi$.  This shows that we
have a maximum of the dual functional in the limit $n\to \infty$,
and it reproduces the expected result.

\begin{figure}[!ht]
\leavevmode
\begin{center}
\epsfysize=4.2cm
\epsfbox{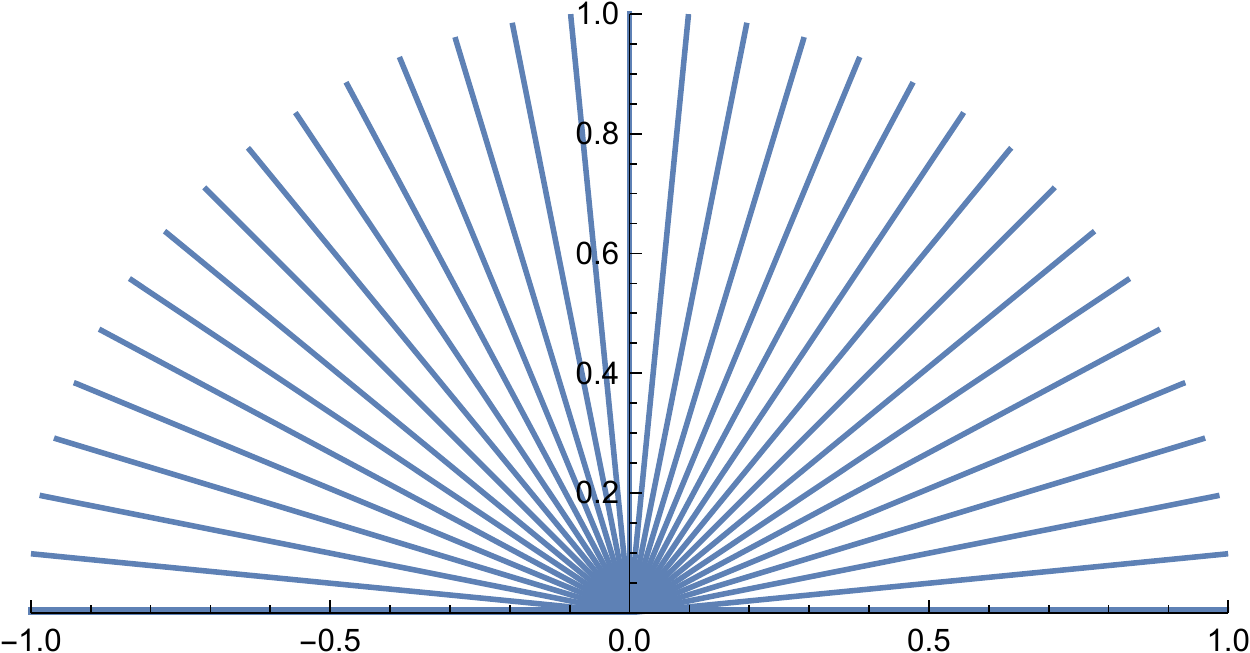} \hskip10pt
\epsfysize=4.2cm
\epsfbox{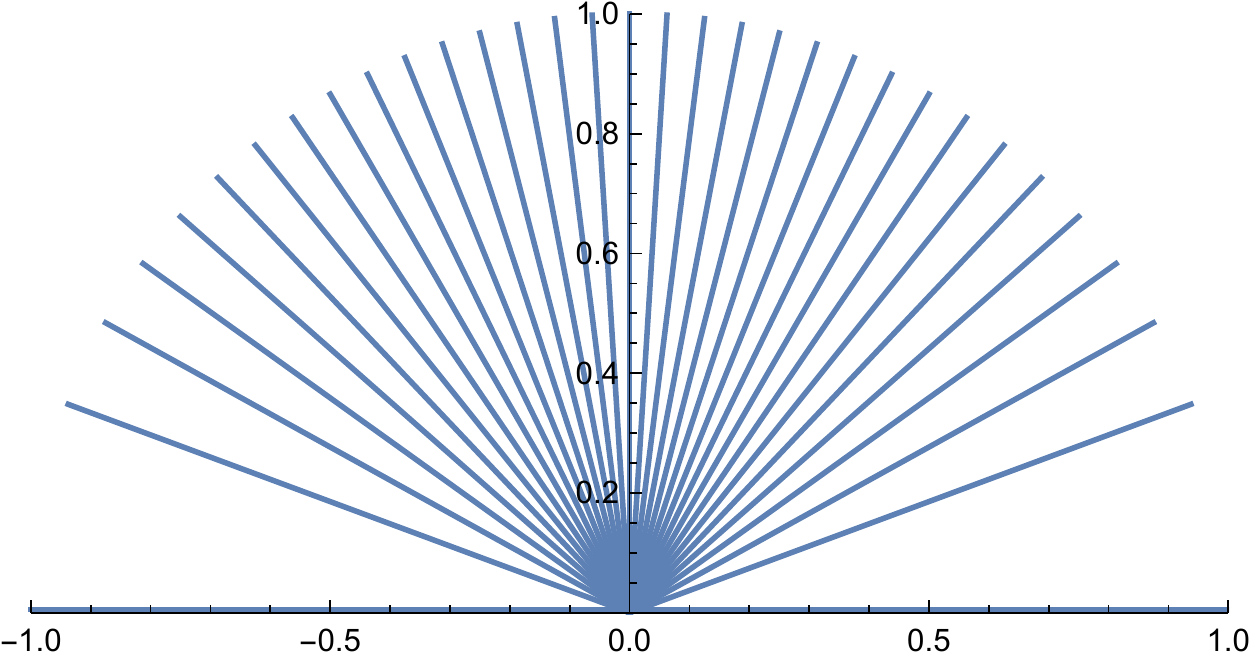}
\end{center}
\caption{\small A closeup of $\phi_0$ geodesics as they depart 
$x_{\phi_0}=0$, placed at the origin of the figure. The 
boundary is represented by the horizontal lines. Left:  With 
$\tilde\varphi_{\phi_0}$ the geodesics depart at identically
spaced azimuthal angles.  Right:  With $\varphi_{\phi_0}$
the geodesics are squeezed around the vertical axis
and are separated out near the boundary. }
\label{ffjnvmcl}
\end{figure}

\begin{figure}[!ht]
\leavevmode
\begin{center}
\epsfysize=7.5cm
\epsfbox{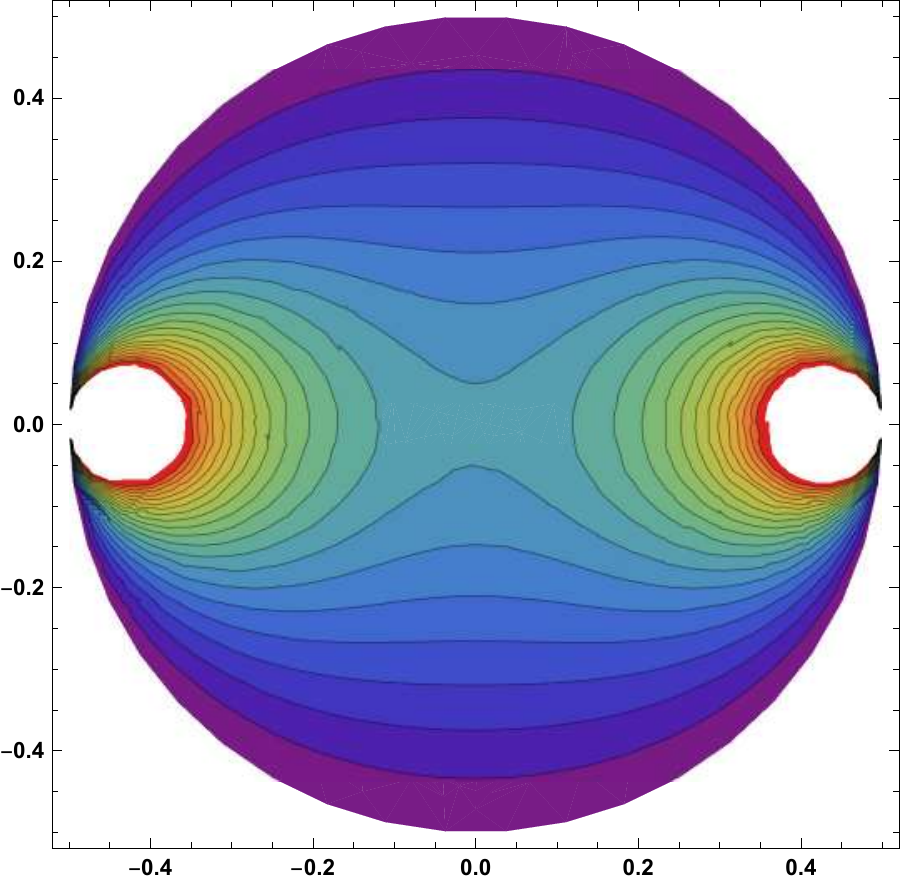} \hskip10pt
\epsfysize=7.5cm
\epsfbox{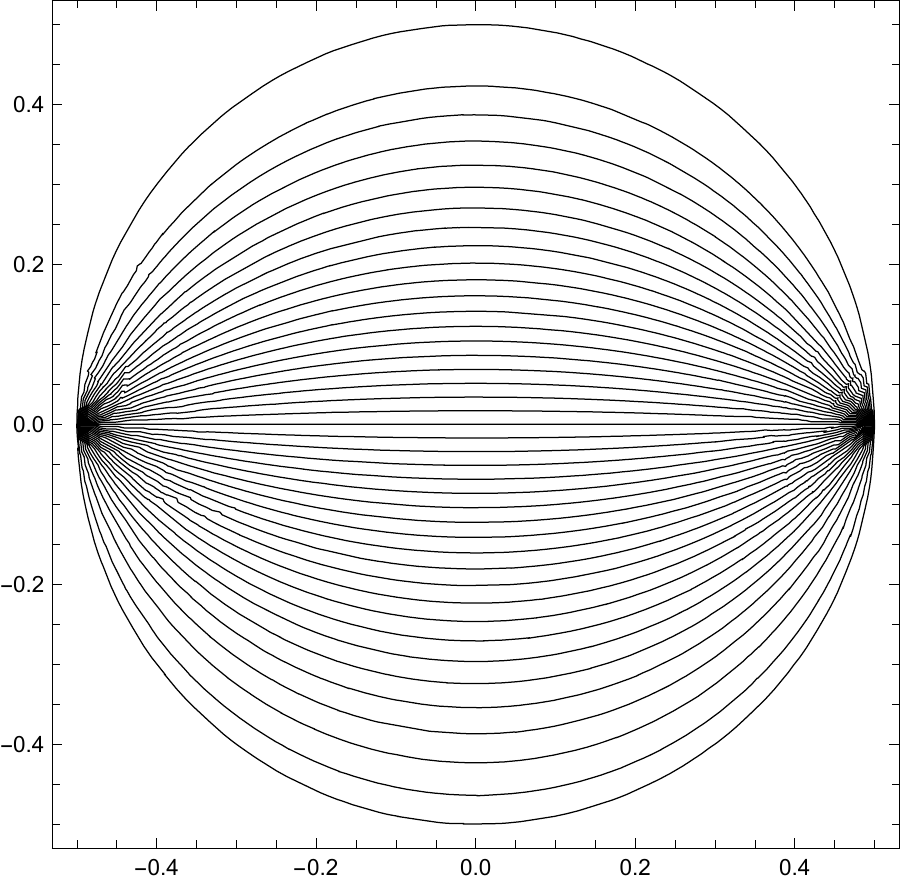}
\end{center}
\caption{\small Left:  Contour 
plot of $|h_{\phi_0}|$ (with $\phi_0=0$) defining the density
of systolic geodesics on the $|z| \leq {1\over 2}$ disk.
Right:  Systolic
geodesics at equal spacings of the
coordinate $\varphi_{\phi_0}$ that labels them.}
\label{ffjncl}
\end{figure}
In order to understand better the solution we consider a few figures.
In Figure~\ref{ffjnvmcl} we show the azimuthal departure angles
for $\phi_0$ geodesics leaving  the boundary point $x_{\phi_0}=0$, 
placed at the origin of the figure.
For $\tilde\varphi_{\phi_0}$, geodesics must
be imagined leaving at equal intervals of this variable.  This is shown
to the left: since $\tilde\varphi_{\phi_0}$ is the azimuthal angle 
relative to the departing point, the geodesics
are uniformly distributed in the interval $(-\tfrac{\pi}{2}, \tfrac{\pi}{2})$.  
To the right, we show the geodesics for the final dual coordinate 
$\varphi_{\phi_0}\sim \sin(\tilde \varphi_{\phi_0})$. Here geodesics equally spaced in this variable
are not equally spaced azimuthally.  They are squeezed near the vertical
axis and are set apart near the boundary.

In Figure~\ref{ffjncl} we consider again the solution
for a single band of systolic geodesics described by coordinates
$(x_{\phi_0},\varphi_{\phi_0})$ with $\phi_0=0$ so that the
antipodal points lie on the real axis of the $|z|\leq {1\over 2}$
disk.   To the left we show the function $|h_{\phi_0=0}|$, which using
(\ref{eq:|hphi|}) is proportional to 
\be
|h_{0}| \sim   {1- 16 (x^2+ y^2)^2\over 1 - 8 (x^2-y^2) + 16 (x^2 + y^2)^2} \,.
\ee 
The function $|h_{\phi_0}|$  has the interpretation of density of geodesics,
or number of geodesics per unit transverse length 
in the convention that 
geodesics are spaced by equal amounts in the $\varphi_{\phi_0}$ 
coordinate~\cite{Headrick:2018ncs}.  
This density function is  singular at the departing and ending points
$(\tfrac{1}{2}, 0)$ and $(-\tfrac{1}{2}, 0)$ where it becomes infinite as approached from the interior.  The function $|h|$ vanishes at the rest of the boundary. 
To the right we show the geodesics themselves, plotted on the $z$-disk 
for equal 
spacings of the $\varphi_{\phi_0}$ coordinate.  Recall that the $(x,y)$ 
coordinates on the disk are nontrivially related to the Cartesian 
coordinates $(\hat x, \hat y)$, that together with $\hat z$ describe
the hemisphere as $\hat x^2 + \hat y^2 + \hat z^2= 1/\pi^2$, with
$\hat z \geq 0$.

\section{Calabi's variational principle} \label{cal_var_pri}

In this section we review the variational principle put forth by Calabi~\cite{calabi} to deal with regions covered by three bands of systolic geodesics and 
give the associated Euler-Lagrange equations. 
This variational principle was further investigated in reference~\cite{bryant},
that proposed an extension to deal with regions covered by more than
three bands of systolic geodesics.  We find, however, some difficulties
with the proposal and its conclusions.  We show that the bands of
geodesics relevant to the metric on $\mathbb{RP}_2$ do not satisfy
the anticipated equation of motion.   We will give our 
modified  
variational principle in section~\ref{iso-var-prin-mult-fol}.

\subsection{The variational principle for regions with three systolic bands}
\label{the-var-pri-for-reg-thr-ban}

Calabi considered the extremal isosystolic metric in a region $U_3$ covered by exactly three bands of geodesics, each band 
calibrated by a closed one-form of unit norm.  
Locally, the calibrating one-form can be written as the differential of the length parameter that Calabi calls `geodesic potential function'.
The calibrating one-forms are thus given by $u^\alpha= d X^\alpha$, with
$\alpha = 1,2,3$.   
We take $(X, Y) \equiv (X^1,X^2)$ as coordinates on $U_3$ and consider the third geodesic potential $Z\equiv X^3$ as a function  $Z\bkt{X,Y}$
of the coordinates. 
In coordinates $\bkt{X, Y }$, the condition that 
$u^1= dX$ and $u^2=dY$  
are one-forms of unit norm fixes the cotangent space metric $g^{-1}$ to be:
\be\label{eq:metinvx1x2}
g^{-1} =\begin{pmatrix}
1 & f \\
f & 1
\end{pmatrix}.
\ee
Here $f$ is a function that appears in the inner product of $dX $ and $dY $: 
\be
\langle u^1, u^2 \rangle = \langle dX , dY  \rangle = f  =  \langle \hat u^1 , \hat u^2 \rangle \,.
\ee
Here $\hat u^1$ and $\hat u^2$ are the unit 
vectors associated to the one-forms
 $u^1$
and $u^2$, defined as usual by the relation 
$u^i (v) = \langle \hat u^i, v \rangle$, with $v$
an arbitrary vector.  Since the unit vectors $\hat u^1$ and $\hat u^2$ are, respectively,
tangent to the first and second band of geodesics, their inner product is
simply the cosine of the angle $\theta_{12}$ between these first two geodesic
bands:
\be
\label{angle-between-geod}
\cos \theta_{12} = f \,. 
\ee
This is the interpretation of $f$.  When $f$ vanishes the bands associated
with $X$ and $Y$ are orthogonal. 
We also require the last one-form,  $u^3= dZ$, to     
be of unit magnitude.  For this we note that
\be
dZ = Z_X   dX  + Z_Y  dY\, ,  \qquad Z_X\equiv{\p Z\over \p X} \,, 
\quad Z_Y\equiv{\p Z\over \p Y}\,.
\ee
Then the norm equal one condition gives
\be
| dZ |^2 = \langle dZ , dZ  \rangle 
=Z_X^2+Z_Y^2+2fZ_XZ_Y = 1 \,.
\ee
We can now solve for $f$ in terms of derivatives of the potential $Z$:
\be
\label{f-from12}
f=\frac{1-Z_X^2-Z_Y^2}{2\,  Z_X Z_Y }\, ,  \ee
showing that the existence  
of the third band determines the metric.
Associated with the cotangent metric $g^{-1}$ above  the metric
$g$ is
\be\label{eq:metx1x2}
g=\frac{1}{1-f^2}\begin{pmatrix}
1 & -f \\
-f & 1
\end{pmatrix} \quad \to \quad \sqrt{\det g} = {1\over \sqrt{1-f^2}}\,.
\ee
The variational principle posits that the area, defined as the integral $I$
of the area form $\sqrt{\det g} \ dX  \wedge dY$,  is stationary under
{\em local} variations $\delta Z$ of the potential $Z$ for the third band.
The local variation, by definition does not change the period of the potential,
which controls the total length of the geodesics.   We thus have
\be
\label{eq:LagCalabiBryant}
I = \int_{U_3}  L(Z_X, Z_Y)  \, dX  \wedge dY \,,  \qquad  L\  = \ {2 |Z_X Z_Y|\over \sqrt{ (2 Z_X Z_Y)^2 - (Z_X^2 + Z_Y^2 -1)^2} }\,. 
\ee
The Euler-Lagrange equation is readily obtained and takes the form
\be
{\partial \over \partial X}  \Bigl[ {\partial L \over \partial Z_X} \Bigr] 
+ {\partial \over \partial Y}  \Bigl[ {\partial L \over \partial Z_Y} \Bigr] = 0 \,.
\ee
To get an explicit version we first find that 
\be
\begin{split}
{\partial L \over \partial Z_X}  \ = \ &  \bigl[ Z_X^4 - (Z_Y^2-1)^2 \bigr]  {2 Z_Y\over \Delta^{3/2}} \,, \\
{\partial L \over \partial Z_Y} \ = \  & \bigl[ Z_Y^4 - (Z_X^2-1)^2 \bigr]  {2 Z_X\over \Delta^{3/2}} \,, 
\end{split}
\ee
where we defined
\be
\Delta \equiv (2 Z_X Z_Y)^2 - (Z_X^2 + Z_Y^2 -1)^2.
\ee 
The final form of the differential equation is long, but worth
recording\footnote{We have dropped an overall factor of $\frac{4}{\Delta^{5/2}}\,$ in arriving at this final form.}:
\be\label{eq:eomfinal}
  f(Z_X, Z_Y) \, Z_{XX}  +  g(Z_X, Z_Y) \, Z_{XY}  +  \tilde f(Z_X, Z_Y)  \, Z_{YY}  \ = \ 0 \,,
\ee
where $f, g$, and $\tilde f$ are functions of $Z_X$ and $Z_Y$ that ordered
by the number of derivatives take the form
\be
\begin{split}
f =\  & \ Z_X Z_Y (\,  3    - 5 Z_X^2 - 3 Z_Y^2  + \ \,  
Z_X^4 + 10 Z_X^2 Z_Y^2  - 3 Z_Y^4   +  
\   Z_X^6 + Z_X^4 Z_Y^2 -5 Z_X^2 Z_Y^4  + 3 Z_Y^6 ) \,, \\
g  = \ &\ 1 - 2 Z_X^2 - 2 Z_Y^2  + 16 Z_X^2 Z_Y^2 \ + 2 Z_X^6 - 10 Z_X^4 Z_Y^2  - 10 Z_X^2 Z_Y^4  + 2 Z_Y^6  \\
& \ \  -Z_X^8 -  4 Z_X^6 Z_Y^2 + 10 Z_X^4 Z_Y^4  - 4 Z_X^2 Z_Y^6 - Z_Y^8 \,, \\
\tilde f  =\  & \ Z_X Z_Y (\,  3   - 5 Z_Y^2 - 3 Z_X^2  
 \ \, + Z_Y^4  + 10 Z_X^2 Z_Y^2 
 - \  3 Z_X^4  \ \, +   Z_Y^6     
  +Z_X^2  Z_Y^4  - 5Z_X^4  Z_Y^2      + \,   3 Z_X^6   ) \,. 
\end{split}
\ee
We note that
\be
\tilde f (Z_X, Z_Y) = f (Z_Y, Z_X)\,, \quad \hbox{and} \quad  g(Z_X , Z_Y) = g(Z_Y ,Z_X) \,.
\ee

\subsection{Connecting the conformal and isosystolic formalisms} \label{conn_iso_form}
\begin{figure}[!ht]
\leavevmode
\begin{center}
\epsfysize=6cm
\epsfbox{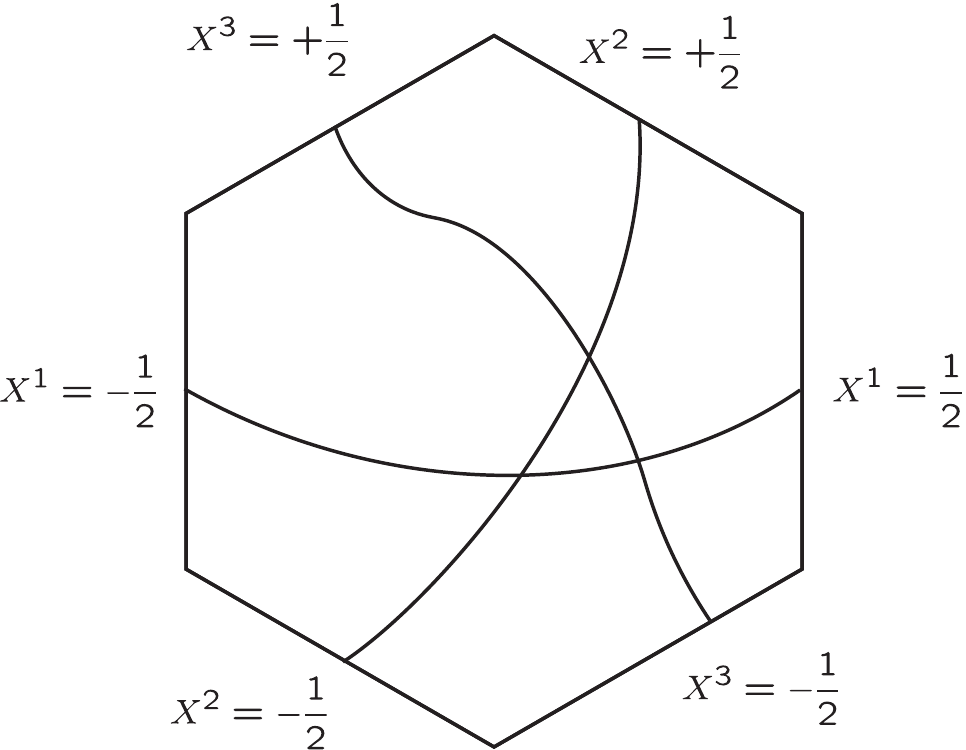}
\end{center}
\caption{\small  The length parameters $X^1 , X^2 ,$ and $X^3 $ along geodesics in classes $C_1,C_2,$ and $C_3$ respectively. The parameter $X^\alpha $ attains the value $-\tfrac12$ on the edge $e_\alpha$ and $+\tfrac12$ on the edge $\tilde{e}_\alpha$. }
\label{fig:6gonBryantPot}
\end{figure}

Our numerical results for the hexagon~(section~\ref{cal_bi_hex})
can be used to obtain the
potential functions $X^1, X^2, X^3$ 
 and to describe the 
area form in terms of these coordinates.
The coordinate $X^\alpha $ attains the value $-\tfrac12$ on the edge $e_\alpha$ and $+\tfrac12$ on the edge $\tilde{e}_\alpha$, as shown for $P_6$ in Figure~\ref{fig:6gonBryantPot}. 

Our program gives a numerical solution for $\phi^\alpha$ everywhere on 
the hexagon and hence the functions $X^1, X^2, X^3$ can also be computed everywhere on the hexagon, described with fiducial coordinates $(x,y)$. The interpretation of $X^\alpha$ 
as the length parameter, however, only holds in the region where the corresponding band of geodesics $C_\alpha$ exists. In relating our results to Calabi's, we will focus on regions of the 
hexagon where the geodesics in class $C_1$ and $C_2$ exist. We shall use the length parameters $\bkt{X^1 ,X^2 }$ along geodesics in these two bands as coordinates. 

Our convex program gives the extremal metric on $P_{6}$ in the original
fiducial coordinates $(x,y)$:
\be\label{eq:dAour}
dA= \rho^2\bkt{x,y} dx\wedge dy.
\ee
Given the $X^\alpha (x,y)$, 
as well as the alternative area form $dA = {1\over \sqrt{1-f^2}}  \, dX^1  \wedge dX^2 $, with $f$ given in (\ref{f-from12}), 
we must have:
\be
\label{eq:M=Fgen}
dA = {1\over \sqrt{1-f^2}}  \, \biggl| \frac{\p \bkt{X^1 ,X^2 }}{\p \bkt{x,y}}
\biggr| \,
dx  \wedge dy\ = \ \tilde \rho^2 (x,y)  \, dx \wedge dy\,.  
\ee
The $\tilde \rho$ so determined in the Riemannian language
must coincide with the conformally determined~$\rho$ in the
region $U_3$.  Of course, since we are using our conformal data throughout,
this is just a consistency check. 
Our data is not good enough to test that 
$Z(X,Y)$ satisfies the differential equation (\ref{eq:eomfinal}) following from the variational principle.

 Our numerical data can be used to compute the local angle $\theta_{12}$ between geodesics.  Using the conformal metric, we have
\be
\cos \theta_{12}=\langle dX^1 ,dX^2 \rangle=\rho^{-2} \bkt{\frac{\p X^1 }{\p x} \frac{\p X^2 }{\p x}+\frac{\p X^1 }{\p y}\frac{\p X^2 }{\p y}}.
\ee
The angle can also be computed from the Riemannian approach 
(\ref{angle-between-geod}) as 
$\sin \tilde\theta_{12} = \sqrt{1-f^2}$, valid  
in the region where the geodesic band associated with $Z$ exists. 
We expect agreement $\theta_{12}= \tilde \theta_{12}$ in $U_3$.

\bigskip

\begin{figure}[!ht]
\leavevmode
\begin{center}
\epsfysize=7cm
\epsfbox{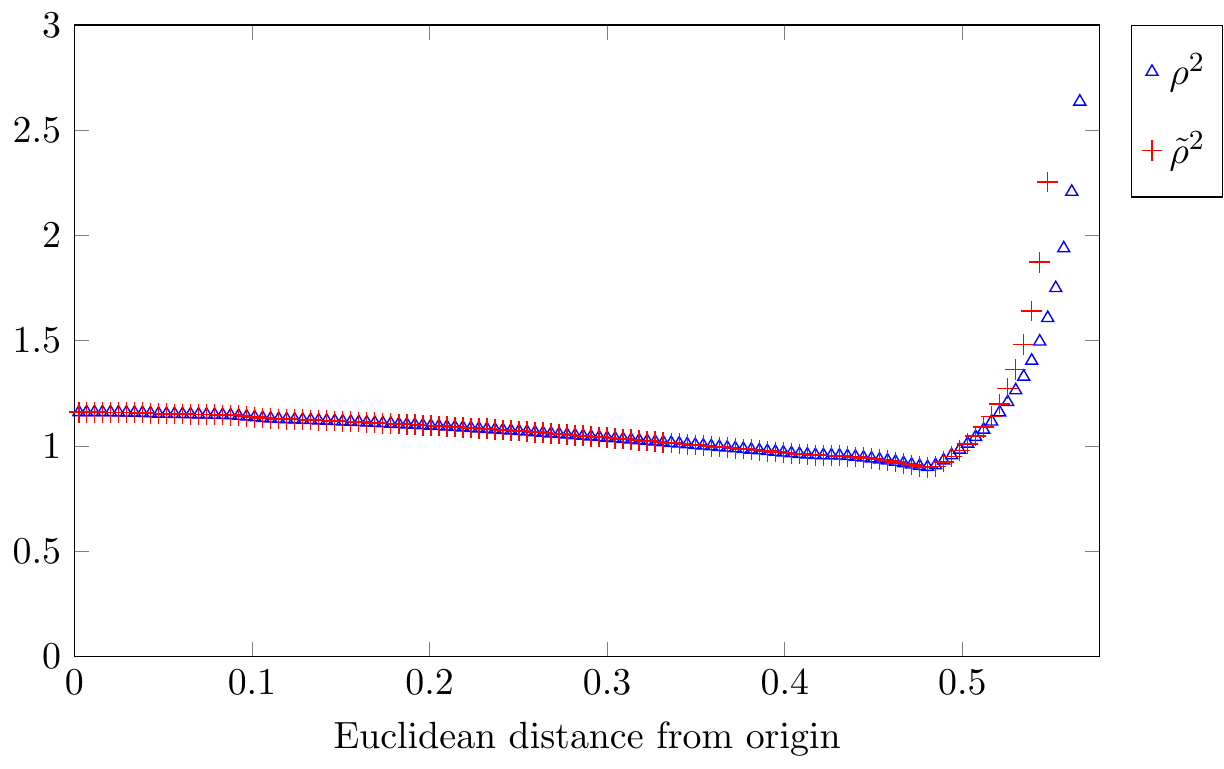}
\end{center}
\caption{\small  The conformal metric $\rho^2$ (blue) and
the alternative $\tilde \rho^2$ (red) from the Riemannian computation.
 Both plotted along the hypotenuse of $T_6$ and agreeing until
 the end of the $U_3$ region at $r\simeq 0.49$.}
\label{fig:MF3DiaHexa}
\end{figure}

Figure~\ref{fig:MF3DiaHexa} shows $\rho^2$ in blue and
$\tilde \rho^2$ in red along the hypotenuse of $T_6$ 
parameterized with the fiducial distance $r = \sqrt{x^2+ y^2}$. 
This line runs from the center of the hexagon to the corner whose neighborhood
is covered by the first and second geodesic bands, but not the third.
As would be expected, we
 find excellent agreement until we cross the transition point
where $U_3$ ends along this line. 
The transition point is the kink in the graph which occurs at 
$r\simeq0.49$. This kink signals line curvature singularity 
as expected on the boundary of $U_3$. As expected, there is 
no agreement beyond $U_3$. 
(The maximum euclidean distance along this hypothenuse is $0.577$.)

\begin{figure}[!ht]
\leavevmode
\begin{center}
\epsfysize=7cm
\epsfbox{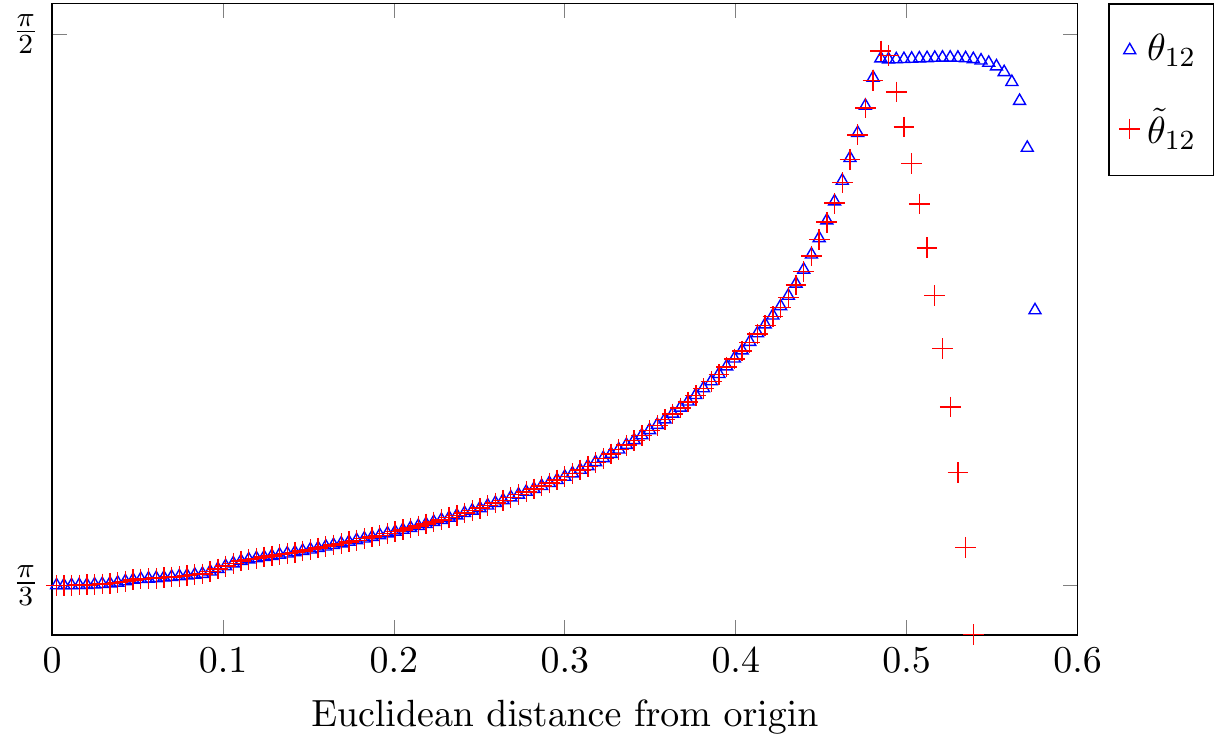}
\end{center}
\caption{\small  The local angle between the first and second bands of 
geodesics plotted along the hypothenuse of $T_6$: 
$\theta_{12}$ (conformal, in blue) and  $\tilde \theta_{12}$ (Riemannian, in red)
As expected, we see agreement until the angle reaches $\frac{\pi}{2}$ at $r\simeq 0.49$,
where $U_3$ ends, and no agreement beyond this point.}
\label{fig:anglesHexagon}
\end{figure}

Figure~\ref{fig:anglesHexagon} shows two evaluations of
the local angle between geodesics of the first
and second classes along the diagonal of $T_6$: $\theta_{12}$ (conformal, in blue) and
$\tilde \theta_{12}$ (Riemannian, in red). 
The boundary
of $U_3$ on the diagonal of $T_6$  occurs when the 
angle $\theta_{12}$ reaches $\pi/2$, consistent with a continuous
transition into $U_2$ where systolic bands must be orthogonal.
Beyond this point the Riemannian determination is not valid, while the
conformal one remains valid and relatively constant at $\pi/2$ until one
approaches the corner, where numerical errors become large.
Consistent with the previous figure, the transition point is reached 
at about $r\simeq 0.49$.

\subsection{Comments on a suggested extension for multiple bands}

Following some remarks in~\cite{calabi}, reference~\cite{bryant} investigated Calabi's variational principle 
and considered the equations that would apply for a
region $U_m$ covered by exactly $m>3$ bands of geodesics.
Choosing potential functions $(X, Y) \equiv (X^1,X^2)$ as coordinates on $U_m$,
one has functions $Z^i\bkt{X,Y}$ for $i=1,\cdots, m-2$ for the remaining
systolic bands.  As we have seen, the condition that $dZ^i$ is of unit norm
determines the cotangent metric $g^{-1}$ 
\be
g^{-1} = \begin{pmatrix} 1 & f \\ f & 1 \end{pmatrix}\,,
\qquad f={1-\bkt{Z_X^i}^2-\bkt{Z_Y^i}^2\over 2 Z_X^i Z_Y^i}\,.
\ee 
Since the metric is unique, all $Z^i$ with $i=1, \cdots, m-2$ must give
the same value of $f$.  In this proposal the Lagrangian $L$ is 
not changed from (\ref{eq:LagCalabiBryant}), one simply demands
that the equation of motion for $Z(X,Y)$ have $m-2$ inequivalent
solutions $Z^i(X,Y)$, each of which gives the same $f$ above.  
Since the Lagrangian
is itself just a function of $f$, all solutions evaluate to the same
Lagrangian.  It was claimed in~\cite{bryant}
that for any non-constant $L$, there are at most two inequivalent solutions 
to this prescribed Lagrangian problem and concluded that 
regions with five or more geodesics are necessarily flat.
We are unable to confirm such variational
principle for $U_m$ with $m>3$.  Our proposal, given
in section~\ref{iso-var-prin-mult-fol},  makes use of 
Lagrange multipliers to implement the constraint of a single $f$,
and ends up modifying the equation of motion of the $Z^i$.

In the rest of this section we point out a difficulty 
with the proposal of~\cite{bryant}.  We consider the minimal area metric on $\mathbb{RP}_2$, a metric with an infinite
number of bands of systolic geodesics.  Not only is the extremal metric
not flat,  we can find the general $Z^i(X,Y)$ that describes the geodesic
bands that leads to a single $f$. 
This $Z^i$, however,  does not satisfy the suggested 
equation of motion~\cite{bryant}. 

We choose 
the potentials $X = x_0$ and $Y= x_{\pi/2}$ as coordinates on the
hemisphere.  These are the potentials for the geodesic bands that originate
at $\phi=0$ and at $\phi= \pi/2$ and travel to the antipodal points.  
Using (\ref{Xphidef}), they are given
by
\be
\begin{split}
X (\theta, \phi) =\ &  \tfrac{1}{\pi}  \cos^{-1} \bigl( \sin \theta \cos \phi \bigr) \,,\\
Y (\theta, \phi) =\ &  \tfrac{1}{\pi}   \cos^{-1} \bigl( \sin \theta \sin\phi \bigr) \,.
\end{split}
\ee
It follows that 
\be
\begin{split}
\cos \pi X  =\ &  \sin \theta \cos \phi \,, \\
\cos \pi Y  =\ &   \sin \theta \sin\phi  \,.
\end{split}
\ee
The ranges of the coordinates $X, Y$ can be determined (although they
play no role in the calculation that follows). It turns out that the condition
$\sin^2 \theta =  \cos^2 \pi X + \cos^2 \pi Y$
following from the two equations above imply that the subdomain $D$ of $X, Y \in [0,1]$
that describes the hemisphere is the set of values for which
$\cos^2 \pi X + \cos^2 \pi Y  \leq 1$.
The domain $D$ is in fact the square formed by the midpoints of the edges in the unit square.

Starting from the round extremal metric on the hemisphere:
\be
ds^2  = \tfrac{1}{\pi^2} (d\theta^2 + \sin^2 \theta d\phi^2 )  \,,
\ee
a short calculation shows that in ($X,Y$) coordinates the components of the
cotangent metric are
\be
g^{XX} = g^{YY} = 1 \,,  \ \ \ g^{XY} = - \cot \pi X  \, \cot \pi Y \,. 
\ee
We can now consider an arbitrary  band departing from $\phi=\phi_0$
and with potential $Z^{\phi_0} = x_{\phi_0}$.  Using  (\ref{Xphidef}) we have:
\be
\label{zzee-phi0}
Z^{\phi_0}(\theta, \phi) =  \tfrac{1}{\pi}   \cos^{-1} \bigl( \sin \theta \cos (\phi-\phi_0) \bigr) \,.
\ee
Therefore
\be
\cos \pi Z^{\phi_0} =  \sin \theta \cos (\phi-\phi_0)  = \sin \theta \cos \phi \cos\phi_0 
+ \sin \theta \sin\phi \sin\phi_0 
\ee
and using the earlier results we get
\be
\label{sol-z}
\cos \pi Z^{\phi_0} \ = \  (\cos \pi X) \cos \phi_0  + (\cos \pi Y) \sin \phi_0
\ee
This defines the function $Z^{\phi_0}(X,Y)$.  We now verify
that for all values of $\phi_0$,  $Z^{\phi_0}$ leads to the same $f$ that
happens to be equal to $g^{XY}$.  
 Indeed, taking derivatives of  (\ref{sol-z})
one quickly finds 
\be\label{eq:ZXZY}
Z^{\phi_0}_X \ =\   {\sin \pi X\over \sin \pi Z^{\phi_0}} \cos \phi_0  \,, \quad
Z^{\phi_0}_Y \ = \  {\sin \pi Y\over \sin \pi Z^{\phi_0}} \sin \phi_0  \,.
\ee
A short calculation then gives
\be
\label{f-determined}
f=  {1- (Z^{\phi_0}_X)^2 - (Z^{\phi_0}_Y)^2\over 2 Z^{\phi_0}_X Z^{\phi_0}_Y}=  -\cot\pi X \cot \pi Y =  g^{XY}\,.
\ee
As expected, the $\phi_0$ dependence drops out.
 This confirms the consistency of the setup: 
after choosing two foliations to define coordinates on the surface, any additional
foliation can be used to define the metric completely, providing the off-diagonal
element. 

On the other hand, this does not guarantee that $Z^{\phi_0} (X,Y)$ 
solves the differential equation (\ref{eq:eomfinal}) 
associated to the Lagrangian $L$. 
The Lagrangian $L$ at the extremum must coincide with the volume
form, but this does not guarantee that the chosen Lagrangian is 
correct. To verify this we compute second partial derivatives.
\be\label{eq:ZXXZYYZXY}
\begin{split}
Z^{\phi_0}_{XX} =\ &  \Bigl( {\cos \pi X \over \sin \pi Z^{\phi_0}} -  \sin^2 \pi X {\cos \pi Z^{\phi_0}\over \sin^3 \pi Z^{\phi_0}} 
 \cos \phi_0\Bigr) \pi \cos \phi_0\,, \\
Z^{\phi_0}_{XY} \ = \ &  - \Bigl( \sin \pi X \, \sin \pi Y 
 {\cos \pi Z^{\phi_0}\over \sin^3 \pi Z^{\phi_0}}  \Bigr)  
 \pi \cos\phi_0 \sin \phi_0  \,,  \\
Z^{\phi_0}_{YY} \  = \  &  \Bigl( {\cos \pi Y \over \sin \pi Z^{\phi_0}} -  \sin^2 \pi Y {\cos \pi Z^{\phi_0}\over \sin^3 \pi Z^{\phi_0}} 
\sin\phi_0\Bigr) \pi \sin \phi_0 \,.
\end{split}
\ee
Evaluating the left-hand side of (\ref{eq:eomfinal})
at various points we see that it does not vanish, and 
therefore the differential equation 
  is not satisfied by $Z^{\phi_0}$.   
This can also be seen
with modest effort by a perturbative analysis near the north pole of the sphere
for a choice $\phi_0 = \pi/4$ as we demonstrate now.
Near the north pole we have $X=Y=\frac12$. 
We do a change of variable $X\to X+\frac12, Y\to Y+\frac12$. We want to evaluate the left-hand side of equation~(\ref{eq:eomfinal}) to first order in $X$ and $Y$. We start by perturbatively expanding $Z^{\pi/4}\bkt{X,Y}$ 
 near the north pole to cubic order in coordinates
\be
Z^{\pi/4}\bkt{X,Y}= \frac12+\frac{X+Y}{\sqrt{2}}+\frac{\pi^2 X Y \bkt{X+Y}}{4 \sqrt{2}} +  \ldots\,. 
\ee
From this one readily computes various partial derivatives of $Z\bkt{X,Y}$:
\be
\begin{split}
Z^{\pi/4}_X&= {1\over \sqrt{2}} + {\pi^2\over 4 \sqrt{2}} Y (2X+ Y) \ \, , \quad Z^{\pi/4}_Y=
{1\over \sqrt{2}} + {\pi^2\over 4 \sqrt{2}} X (X+ 2Y) \ , \\
Z^{\pi/4}_{XX}&=\frac{\pi^2 Y}{2\sqrt2}\, , 
\qquad Z^{\pi/4}_{XY}=\frac{\pi^2 \bkt{X+Y}}{2\sqrt{2}}\, , 
\qquad Z^{\pi/4}_{YY}=\frac{\pi^2 X}{2\sqrt2}\, .
\end{split}
\ee
From this we compute various functions appearing in the EoM. To leading order, we have
\be
f\bkt{Z^{\pi/4}_X,Z^{\pi/4}_Y}=\tilde{f}\bkt{Z^{\pi/4}_X,Z^{\pi/4}_Y}=
\tfrac{1}{2} +{\cal O}\bkt{X^2}\,\ \quad g\bkt{Z^{\pi/4}_X,Z^{\pi/4}_Y}=1+{\cal O}\bkt{X^4}.
\ee 
The left-hand side of the equation of motion (\ref{eq:eomfinal}) then evaluates to
\be
\frac{3\pi^2\bkt{ X+Y}}{4\sqrt{2}}\neq 0,
\ee
explicitly demonstrating that the equation
 is not satisfied.  We will see in the following 
 section that the modified variational principle gives an equation
 that is satisfied by the potentials in the $\mathbb{RP}_2$ 
 extremal metric.

\section{Isosystolic variational principle with multiple foliations}
\label{iso-var-prin-mult-fol}

We propose here a variational principle to determine the 
extremal minimal area Riemannian in a region $U_m$ covered
by $m$ bands of systolic geodesics.  For the case $m=3$ the 
answer was provided by Calabi~\cite{calabi}.  We provide here
a variational principle that works for $m\geq 3$, and differs from
the proposal in~\cite{bryant} for $m>3$.   

After introducing this modified variational approach,  
we apply it to the extremal
Riemannian metric in $\mathbb{RP}_2$, showing that the equations
of motion following from the variational principle are satisfied.
The computations are somewhat analogous to those that showed
that the extremal metric in $\mathbb{RP}_2$ satisfy the equations
of the dual convex program.

\subsection{Extension of the variational principle}
\label{ext-var-principle}

Assume that we have a region $U_m$ which is covered by exactly $m\geq 3$ bands of geodesics. We choose any two geodesics potential functions $X$ and $Y$ as coordinates on $U_m$ and express the remaining others as functions 
 $Z^i(X,Y)$, with $i = 1, \ldots , m-2$.  
With cotangent space metric $g^{XX} = g^{YY} = 1$ (to guarantee $|dX| = |dY| =1$) and $g^{XY} \equiv f$,  
to enforce the constraints $|dZ^i|=1$ for all $i$ we require Lagrange multipliers
$\lambda_i$.   We thus have the variational principle:
\be
I = \int  dX dY   {1\over \sqrt{1-f^2} }\Biggl(1 + \sum_i  \lambda_i \Bigl[  (Z^i_X)^2+(Z^i_Y)^2+2 f Z^i_X Z^i_Y-1 \Bigr] \Biggr) \,.
\ee
Here $f, \lambda_i,$ and $Z^i$ must be viewed as functions of $X$ and $Y$
and are all to be varied to find stationary points of the functional $I$. 
As usual $Z_X^i = \partial_X Z^i$ and $  Z_Y^i = \partial_Y Z^i$.

The equations of motion for the Lagrange multipliers $\lambda_i$ are:
\be
\label{eq:lambda-eom}
(Z^i_X)^2 + (Z^i_Y)^2+2 f Z^i_X Z^i_Y-1 \ =   \ 0 \ \,, \quad i= 1, \cdots\,,  n-2 \,. 
\ee
Indeed this sets $f$ equal to the desired value 
\be\label{eq:fdetermination}
f\ = \ {1 -(Z^i_X)^2 - (Z^i_Y)^2\over 2  Z^i_X Z^i_Y} \,,
\ee
and makes the determination of $f$ from each of the bands equal.
Moreover, when this holds for all $i$ the action $I$ reduces to the
integral of the area form (see~(\ref{eq:metx1x2})). 
 The variation of $f$ gives an equation of motion
that using (\ref{eq:lambda-eom}) becomes
\be
\label{fff-eom}
{f\over 1-f^2} +2  \sum_i \lambda_i  Z_X^i Z_Y^i \ =   \ 0 \,.
\ee
Finally, the variation of $Z^i$ gives the equation of motion
\be\label{eq:Zeom}
\, \partial_X \biggl[ \frac{\lambda_i}{\sqrt{1-f^2}}  \bkt{Z^i_X+f Z^i_Y}\biggr]  
+\partial_Y  \biggl[ \frac{\lambda_i}{\sqrt{1-f^2}}  \bkt{Z^i_Y+f Z^i_X}\biggr]  
\ = \ 0 \,. 
\ee
Recalling that 
\be
g^{-1} = \begin{pmatrix} 1 & f \\ f & 1 \end{pmatrix}  \,, \qquad   \sqrt{g} = {1\over \sqrt{1-f^2}} \,,
\ee
and letting $\mu, \nu$ indices run over $X, Y$, 
the last equation can be written covariantly as
\be\label{eq:eomZ}
\p_\mu \bkt{
\lambda_i \sqrt{g} \, g^{\mu\nu} \partial_\nu Z^i} = 0  \,. 
\ee
Using differential form notation the equation reads
\be\label{eq:newEoMvm}
d*(\lambda_i \, d Z^i)=0\, .
\ee
Evaluating the derivatives in  (\ref{eq:Zeom}) and collecting
terms we find a more explicit form of the $Z^i$ equation of motion:
\be
\label{explicit-kj-Z-phi-eom}
\begin{split}
0 \, = \ &    
\  Z^i_{XX} + Z^i_{YY}  + 2 f\, Z^i_{XY}  \\
& \hskip-5pt + Z^i_X \biggl(  {\lambda_{iX} + f\lambda_{iY}\over \lambda_i }
+  {f_Y + f f_X\over 1-f^2}  \biggr) 
\, + Z^i_Y \biggl(  {\lambda_{iY} + f\lambda_{iX}\over \lambda_i }
+  {f_X + f f_Y\over 1-f^2}  \biggr) \,. 
\end{split}
\ee
It is straightforward to verify that for $n=3$ the new equations of motion are equivalent to those of~\cite{calabi,bryant}, which we reviewed in~\ref{the-var-pri-for-reg-thr-ban}.

\subsection{Extremal Riemannian metric on $\mathbb{RP}_2$ from
new variational principle}
\label{ext-ire-met-new-var}

In this section we
 show that the extremal Riemannian metric on $\mathbb{RP}_2$ satisfies the Euler-Lagrange equations obtained from the new variational principle.  
 We begin with some preliminary discussion of equation (\ref{eq:newEoMvm}),
\be\label{eq:newEoM}
d*(\lambda_i \, d Z^i)=0.
\ee
This is very similar to the equation of motion derived using conformal methods (see discussion after eq. (7.19) of~\cite{Headrick:2018ncs}). In that notation
we had
\be
\label{vmcljn}
u^\alpha =-  {*d\varphi^\alpha \over | d\varphi^\alpha| }\,, 
\ee
with the local metric taking the form 
\be
ds^2 = (dx^\alpha)^2  + {1\over h_\alpha^2} (d\varphi^\alpha)^2 \,. 
\ee
It follows that $|d\varphi^\alpha | = |h_\alpha|$, where $h_\alpha$ has
the interpretation of the density of $\alpha$ geodesics.  
In fact, $\varphi^\alpha$ are the functions that appear in the dual program
and are constants along $\alpha$-geodesics.  We also note
that $u^\alpha = dx^\alpha$ where $x^\alpha$ is the length parameter along 
$\alpha$-geodesics.  Therefore (\ref{vmcljn}) and $* (* u)= -u$, gives 
\be
\, d\varphi^\alpha =  \,   * (|h_\alpha|  dx^\alpha)\,.  
\ee
Taking another exterior derivative we get
\be
d*(|h_\alpha| dx^\alpha) =0.
\ee
We verified this equation explicitly in section~\ref{con_met_on_rp2} where
we showed that the metric on $\mathbb{RP}_2$ was a solution of the
dual (conformal) program.   This equation is exactly analogous to
equation (\ref{eq:newEoM}) with the Lagrange multipliers identified
with the functions that determine geodesic density.  
 The $Z^{\phi_0}$ in (\ref{zzee-phi0}) 
are in fact length coordinates
of the $x^\alpha$ type.  As before
 $\phi_0$ is the label for a systolic band.  
 In this language we would have a metric defined via the band:
\be 
\label{met-z-tilde-varphi}
ds^2 =  (dZ^{\phi_0})^2 +    {1\over \tilde h_{\phi_0}^2} (d\tilde\varphi_{\phi_0})^2 \,. 
\ee
It follows that if we choose $\lambda_{\phi_0}\sim  |\tilde h_{\phi_0}|$ with 
$|\tilde h_{\phi_0}|$ given in equation~(\ref{tilde-h-first}) then the 
equation of motion (\ref{eq:newEoM}), reading now
$d*(\lambda_{\phi_0} \, d Z^{\phi_0})=0$ is guaranteed to be satisfied. 
We will anyway check the equation explicitly. 

We note, however, the analog of a reparameterization ambiguity 
discussed around (\ref{h-repara}) for the dual-program analysis of the 
$\mathbb{RP}_2$ metric. 
If a pair $(\lambda_i,Z^i)$ satisfies equation (\ref{eq:eomZ}) then 
a pair 
\be
( \chi(\varphi^i) \lambda_i\, ,\,  Z^i)  \,,
\ee
 with a rescaled $\lambda_i$
will also satisfy the equation if
 \be
 g^{\mu\nu} \p_\mu Z^i \p_\nu \varphi^i=0 \,. 
 \ee
But this holds when $(Z^i, \varphi^i)$ are a pair of orthonormal 
coordinates for a systolic band.  In our case, the 
coordinate pair is $(Z^{\phi_0}, \tilde\varphi_{\phi_0})$ and the metric (\ref{met-z-tilde-varphi}) indeed implies that
\be
0= \langle dZ^{\phi_0},  d\tilde\varphi_{\phi_0} \rangle  
=  g^{\mu\nu} \p_\mu Z^{\phi_0}  \p_\nu \tilde\varphi_{\phi_0} \,. 
\ee
We thus have the freedom to replace solutions as
\be
\label{rep-sol-final} 
(\lambda_{\phi_0} \,, \, Z^{\phi_0} ) \ \to  \  
( \,  \chi(\tilde \varphi_{\phi_0}) \,  \lambda_{\phi_0} \,, \, Z^{\phi_0} ) \,. 
\ee

To confirm that we have a solution for the extremal metric on $\mathbb{RP}_2$
we consider  an arbitrary  systolic band departing from $\phi=\phi_0$ for which the length parameter $Z^{\phi_0}$ is given in equation~(\ref{sol-z}):
 \be
\label{sol-z-again}
\cos \pi Z^{\phi_0} \ = \  (\cos \pi X) \cos \phi_0  + (\cos \pi Y) \sin \phi_0\,.
\ee
We also recall that the extremal metric  has (\ref{f-determined}): 
\be
f=-\cot\pi X\cot\pi Y \,.
\ee
These two, $Z^{\phi_0}$ and $f$, with values so specified must solve
the equations of motion.  The only unknown is the value of the Lagrange
multipliers $\lambda_{\phi_0}$.  From the above remarks, however,  
we can choose $\lambda_{\phi_0}$ equal to the first choice $\tilde h$ in
(\ref{tilde-h-first}): 
\be
\label{lambda-phi-0-f-c}
\lambda_{\phi_0}=\frac{1}{\sin\pi Z^{\phi_0}}, 
\ee
keeping in mind that we are free to multiply it with an arbitrary function of $\tilde\varphi_{\phi_0} $.  Our task now is to verify explicitly that we have a solution, and
in doing so we will find the final form of $\lambda_{\phi_0}$.

The equation of motion for $Z^{\phi_0}$ is (\ref{explicit-kj-Z-phi-eom}), with 
$i\to \phi_0$.  With a little reordering it reads: 
\be
\begin{split}
 0\  =  \ &  \ \  Z^{\phi_0}_{XX} + Z^{\phi_0}_{YY}  + 2 f\, Z^{\phi_0}_{XY} 
  +{1\over \lambda_{\phi_0} }
  \Bigl(  
  Z^{\phi_0}_X 
  (  \lambda_{{\phi_0}X} + f\lambda_{{\phi_0}Y} )   + Z^{\phi_0}_Y  (\lambda_{{\phi_0}Y} + f\lambda_{{\phi_0}X}) \Bigr) 
\\
& \hskip-5pt +{ 1\over 1-f^2} \, \Bigl[ 
f_X\bkt{Z^{\phi_0}_Y+f Z^{\phi_0}_X}
+ f_Y\bkt{Z^{\phi_0}_X+f Z^{\phi_0}_Y}\Bigr]  \,.
\end{split}
\ee
We now claim that for $Z^{\phi_0}$ the first three terms in the above
equation add up to zero:
\be
 Z^{\phi_0}_{XX} + Z^{\phi_0}_{YY}  + 2 f\, Z^{\phi_0}_{XY}  = 0\,.
\ee
This is readily verified with a bit of work 
using the second derivatives computed in
(\ref{eq:ZXXZYYZXY}). 
Since $\lambda_{\phi_0}$ is just a function of $Z^{\phi_0}$, using 
chain rule we quickly see that
\be
\begin{split}
\hskip-10pt{1\over \lambda_{\phi_0} }\biggl( Z^{\phi_0}_X 
  (  \lambda_{{\phi_0}X} + f\lambda_{{\phi_0}Y} )   + Z^{\phi_0}_Y  (\lambda_{{\phi_0}Y} + f\lambda_{{\phi_0}X})\biggr)
   = & \ \  {1\over\lambda_{\phi_0} } {d\lambda\  \over dZ^{\phi_0}}  [ (Z^{\phi_0}_X)^2 + (Z^{\phi_0}_Y)^2  
   + 2 f Z^{\phi_0}_X  Z^{\phi_0}_Y  ] \\
   = & \ \  {1\over\lambda_{\phi_0} }{d\lambda\  \over dZ^{\phi_0} } = - \pi \cot \pi 
   Z^{\phi_0}  \,. 
\end{split}
\ee
In passing to the second line we recalled that the expression in 
brackets is $|dZ^{\phi_0}|^2 = 1$. 
The equation of motion has thus simplified to
\be
\label{just-remains}
0\ = \ - \pi \cot \pi 
   Z^{\phi_0}+{ 1\over 1-f^2} \, \Bigl[ 
f_X\bkt{Z^{\phi_0}_Y+f Z^{\phi_0}_X}
+ f_Y\bkt{Z^{\phi_0}_X+f Z^{\phi_0}_Y}\Bigr]  \,.
\ee
We use the first derivatives~(\ref{eq:ZXZY}) to compute:
\be
\begin{split}
f_X\bkt{Z^{\phi_0}_Y+f Z^{\phi_0}_X}
& =
{\pi \cos\pi Y\over \sin\pi Z^{\phi_0}\sin \pi X\sin \pi Y}\bkt{{\sin \pi Y \over \sin\pi X}\sin\phi_0 -\cot\pi X\cot\pi Y\cos\phi_0} \,,\\
f_Y\bkt{Z^{\phi_0}_X+f Z^{\phi_0}_Y}
& =
{\pi \cos\pi X\over \sin\pi Z^{\phi_0}\sin \pi X\sin \pi Y}\bkt{{\sin \pi X\over \sin\pi Y}\cos\phi_0 -\cot\pi X\cot\pi Y\sin\phi_0}\, . 
\end{split}
\ee
Adding above two terms and simplifying we find
\be
\begin{split}
f_X\bkt{Z^{\phi_0}_Y+f Z^{\phi_0}_X}
+ f_Y\bkt{Z^{\phi_0}_X+f Z^{\phi_0}_Y}= & \ \pi \cot\pi Z^{\phi_0} \bkt{1-\cot^2\pi X\cot^2\pi Y}\\ 
= & \ \pi \cot\pi Z^{\phi_0}\bkt{1-f^2}\,. 
\end{split}
\ee
Back in (\ref{just-remains}) we have a complete cancellation. This
completes
the verification that the equation of motion for $Z^{\phi_0}$ is satisfied.

We now verify the equation of motion for $f$, given in (\ref{fff-eom}) and 
rearranged as follows: 
\be
\label{integral-for-fvm}
2\,  \frac{1-f^2}{f}\sum_{\phi_0} \lambda_{\phi_0} Z_X^{\phi_0} Z_Y^{\phi_0}=\, -1 \,.
\ee
Multiplying by $\Delta \phi_0 = {\pi/n}$,
for the case of $n$ systolic bands, with $n$ large, 
we turn the left-hand side into an integral:
\be
\label{integral-for-f}
2 \, \frac{1-f^2}{f}\int_0^\pi d\phi_0\,  \lambda_{\phi_0} 
\, Z_X^{\phi_0} Z_Y^{\phi_0}\ =\ -{\pi\over n}\,.  
\ee
It is easiest to verify the equation in fiducial disk coordinates $(x,y)$
using the parameters $a$ and $b$ defined in equation~(\ref{eq:abconstraints}). These different representations are obtained by straightforward  computations, using the expressions for $X, Y$, and $Z$ 
in terms of $(x,y)$ given in~(\ref{pi-x-phi-def}).  We first consider separately
the various ingredients of the above left-hand side.   For the prefactor, 
\be
\begin{split}
{ 1-f^2 \over f}&\ =\  \cot\bkt{\pi X} \cot\bkt{\pi Y}-\tan\bkt{\pi X}\tan\bkt{\pi Y} \\
 &\ =\ 
-\,  {\bkt{1-16|z|^4}^2\over 16 x y\sqrt{\bkt{1-16|z|^4}^2+256 x^2 y^2}}
\,  =\, 
-{ 1-a^2-b^2\over b\sqrt{1-a^2}}\,. 
\end{split}
\ee
We now consider $\lambda_{\phi_0}$, already given in (\ref{lambda-phi-0-f-c}).
Recall that we are free to scale $\lambda_{\phi_0}$ by an arbitrary function which is constant along $\phi_0$-geodesics without spoiling the equation of motion for $Z^{\phi_0}$. We choose to make $\lambda_{\phi_0}$ 
proportional  to $|h_{\phi_0}|$ in equation~(\ref{h-reparam}).  For this
we multiply our original choice by $\cos{\tilde\varphi_{\phi_0}}$ and 
a constant $\beta$ to be determined: 
 \be
\lambda_{\phi_0}= \beta |h_{\phi_0}|=\beta{\cos\tilde\varphi_{\phi_0}\over \sin\bkt{\pi Z^{\varphi_0}}}\,. 
\ee
Recall that 
$\tilde{\varphi}_{\phi_0}$ in disk coordinates is given in~(\ref{eq:varphiDisk})
and is, by construction,  constant along $\phi_0$-geodesics.   
In terms of the disk coordinates and parameters $a$ and $b$ we have
\be
\begin{split}
\lambda_{\phi_0}&= \   {  \beta(1-16|z|^4)\over 1+16|z|^4-8\bkt{x^2-y^2}\cos 2\phi_0-16 xy\sin 2\phi_0}
=\, 
{\beta \sqrt{1-a^2-b^2}\over 1-a \cos 2 \phi _0-b \sin 2 \phi _0}\,.
\end{split}
\ee
We now evaluate the product $Z_X^{\phi_0} Z_Y^{\phi_0}$, beginning with equation~(\ref{eq:ZXZY}):
\be
\begin{split}
 Z_X^{\phi_0} Z_Y^{\phi_0} &\ = \ 
  {\sqrt{\bkt{1+16|z|^4}^2+256 x^2 y^2} \ \sin 2\phi_0\over 
  2 \bkt{1+16|z|^4-8\bkt{x^2-y^2} \cos 2\phi_0 -16 x y \sin 2 \phi_0}}\,, 
 \\
 &=
 \frac{\sqrt{1-a^2} \sin \left(2 \phi _0\right)}{2 \left(1-a\cos 2 \phi _0
 -b \sin 2 \phi _0\right)}\,. 
\end{split}
\ee
Collecting the above ingredients we now get:
\be
2 \, \frac{1-f^2}{f} \lambda_{\phi_0} Z_X^{\phi_0} Z_Y^{\phi_0}
=
-{\beta\over b} \left(1-a^2-b^2\right)^{3/2}
\frac{\  \sin 2 \phi _0}{ \left(1-a \cos 2 \phi _0-b \sin 2 \phi _0\right){}^2}.
\ee
We perform the integral required in (\ref{integral-for-f}): 
\be
\label{integral-for-f-jn}
-{\beta\over b} \left(1-a^2-b^2\right)^{3/2} \int_0^\pi \frac{\  
\sin 2 \, \phi _0 d\phi_0}{ \left(1-a \cos 2 \phi _0
-b \sin 2 \phi _0\right){}^2}=-{\pi\over n}\,.   
\ee
The integral can be evaluated by relating 
 it to an already computed integral in~(\ref{eq:confint}):
\be
\begin{split}
\int_0^{\pi } \frac{ \sin 2 \phi_0\, d\phi _0 }{\left(1-a\cos 2\phi_0-b\sin 2\phi_0
\right)^2}\, 
&=
\frac{d}{d b} \int_0^\pi {d\phi_0\over \left(1-a\cos 2\phi_0
-b\sin 2\phi_0\right)} \\
&= \frac{d}{db }
 {\pi \over \sqrt{1 - (a^2 + b^2 )} }={\pi b
 \over \bkt{1-a^2-b^2}^{3/2}
 }\,.
\end{split}
\ee
Back in (\ref{integral-for-f-jn}) all the position dependence cancels, showing
that the equation is satisfied  when the constant $\beta$ is chosen to be
\be
- \beta\pi=-{\pi\over n}\quad \to \quad  \beta=\frac{1}{n}\, .
\ee
This completes the verification that the extremal metric in
$\mathbb{RP}_2$  satisfies the Euler-Lagrange equations coming from the new variational principle.

\section*{Acknowledgements}
We are indebted to Matthew~Headrick for his help in constructing
the argument in section~\ref{riem_conf} and for conversations that
led to the variational principle of section~\ref{iso-var-prin-mult-fol}.
We thank Michael~Wolf for instructive discussions and Liam
Cohen for discussions on triangulations
at the very early stages of this work. 
This material is based upon work supported by the U.S. Department of Energy, Office of Science, Office of High Energy Physics of U.S. Department of Energy under grant Contract Number  DE-SC0012567. U.N also acknowledges the fellowship by the Knut and Alice Wallenberg Foundation,
Stockholm, Sweden.

\begin{appendix}
\section{Discretization of polygons for numerical analysis}\label{discretization_appendix}

\subsection{Hexagon}

\begin{figure}[ht]
\leavevmode
\begin{center}
\epsfysize9cm
\epsfbox{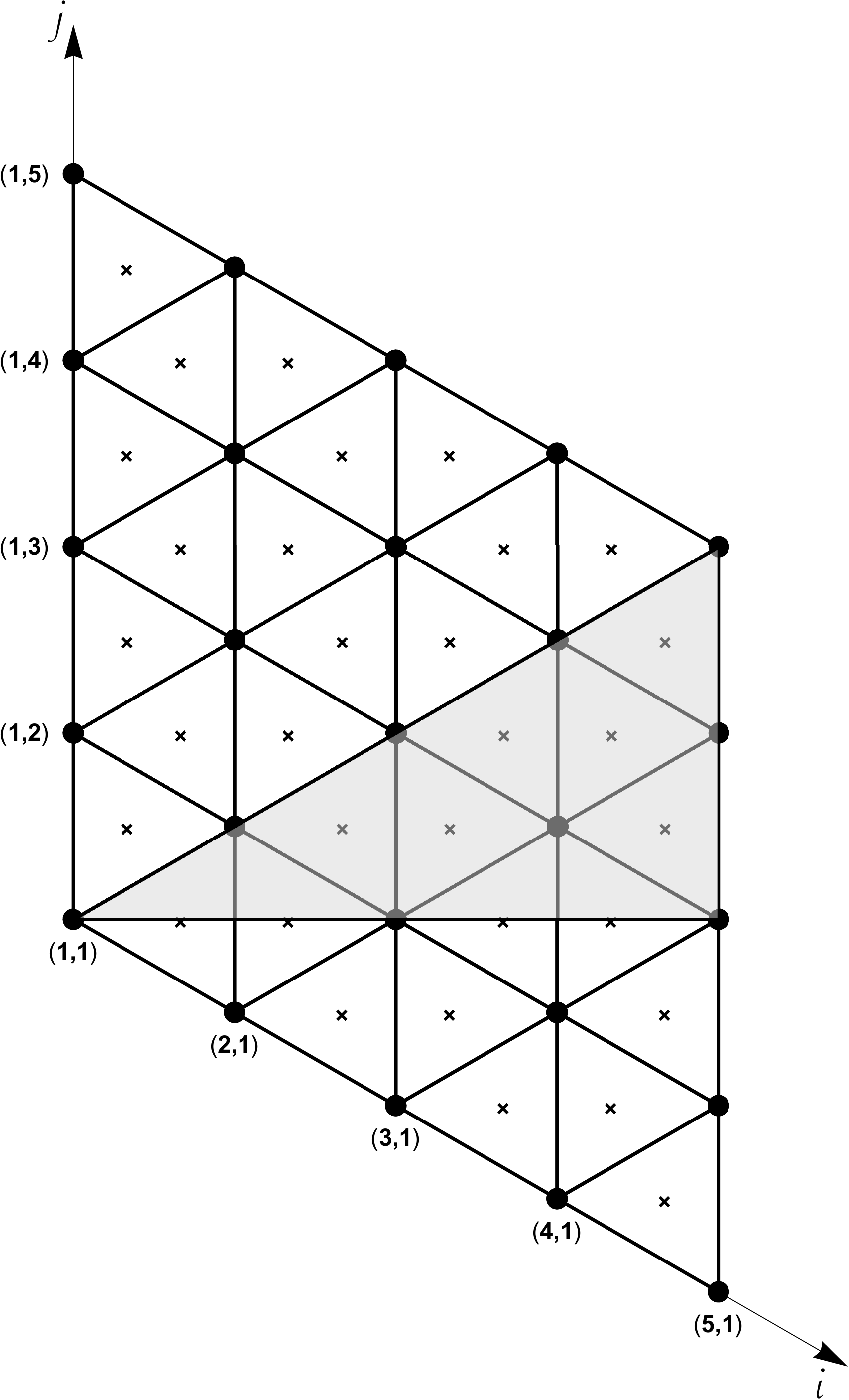}
\end{center}
\caption{\small Discretization of the parallelogram for $N_c=4$. The shaded region is the fundamental domain $T_6$ discussed in section~\ref{the_for_gen_pol}. 
The fuunctions $\phi^1$ and $\varphi^1$ is defined at the vertices of plaquettes, denoted by `dots'.  Their derivatives are evaluated at the centroid of each plaquette, denoted by a `cross'. }
\label{fig:Hex_Disc}
\end{figure}

For the discretization of the hexagon  the plaquettes are equilateral triangles 
that fit seamlessly into the parallelogram with vertices at $(0,0),(1,-\frac{1}{\sqrt{3}}),(1,\frac{1}{\sqrt{3}})$ and $(0,\frac{2}{\sqrt{3}})$, as shown in Figure~\ref{fig:Hex_Disc}.
We denote by $N_c$ the number of plaquette edges on each side of the parallelogram. 
For any fixed $N_c$, the region is triangulated by $2 N_c^2$ plaquettes.
The $(i,j)$ labels at the corners of the plaquettes are based on axes running along
the edges of the parallelogram.  We have $i,j=1,\cdots , N_{c}+1$.
Figure~\ref{fig:Hex_Disc} shows a discretization for $N_c=4$.

\begin{figure}[ht]
\leavevmode
\begin{center}
\epsfysize7cm
\epsfbox{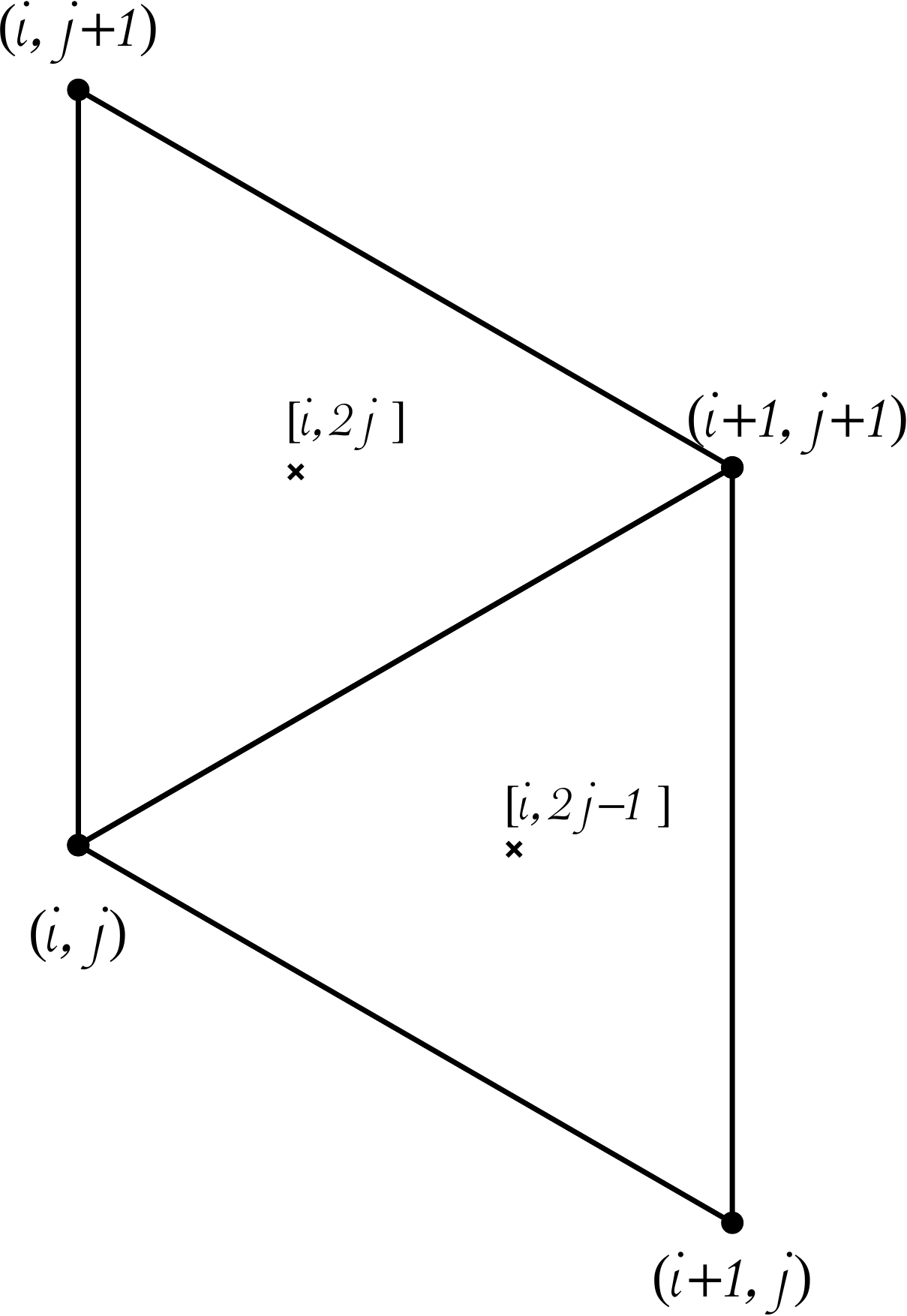}
\end{center}
\caption{\small  Labeling the centroids of  plaquettes
 with two integers $\sbkt{\tilde{i},\tilde{j}}$. 
 We use the left-lower vertex $(i,j)$ to generate labels for two centroids. Both are assigned the value $\tilde{i}=i$. For the centroid to the right of the vertex $(i,j)$ we assign $\tilde{j}=2j-1$ and for  the other $\tilde{j}=2j$.  }
\label{fig:Hex_lab}
\end{figure}

As discussed in section~\ref{the_for_gen_pol}, one only needs to use the 
region $Q_6$ of the hexagon in the first quadrant as the fundamental region for the primal variable $\phi^1$ and for the dual variable $\varphi^1$.  Moreover, one 
effectively works with the region $T_6$ (shaded on the figure) which is 
the fundamental domain
of the metric. 
The functions $\phi^\alpha$ and $\varphi^\alpha$ are defined on the vertices
of the plaquettes. Derivatives are defined at the centroid of each plaquette
 in terms of the value of the function at the vertices of the plaquette. In figure~\ref{fig:Hex_Disc}, centroids are marked by a $\times$.  We label the centroids by $\sbkt{i,j}$ where $i=1,2,\cdots, N_c$ and $j=1,2,\cdots 2 N_c$ as shown in Figure~\ref{fig:Hex_lab}. 

\subsection{Octagon and higher polygons}

\begin{figure}[!ht]
\leavevmode
\begin{center}
\epsfysize5cm
\epsfbox{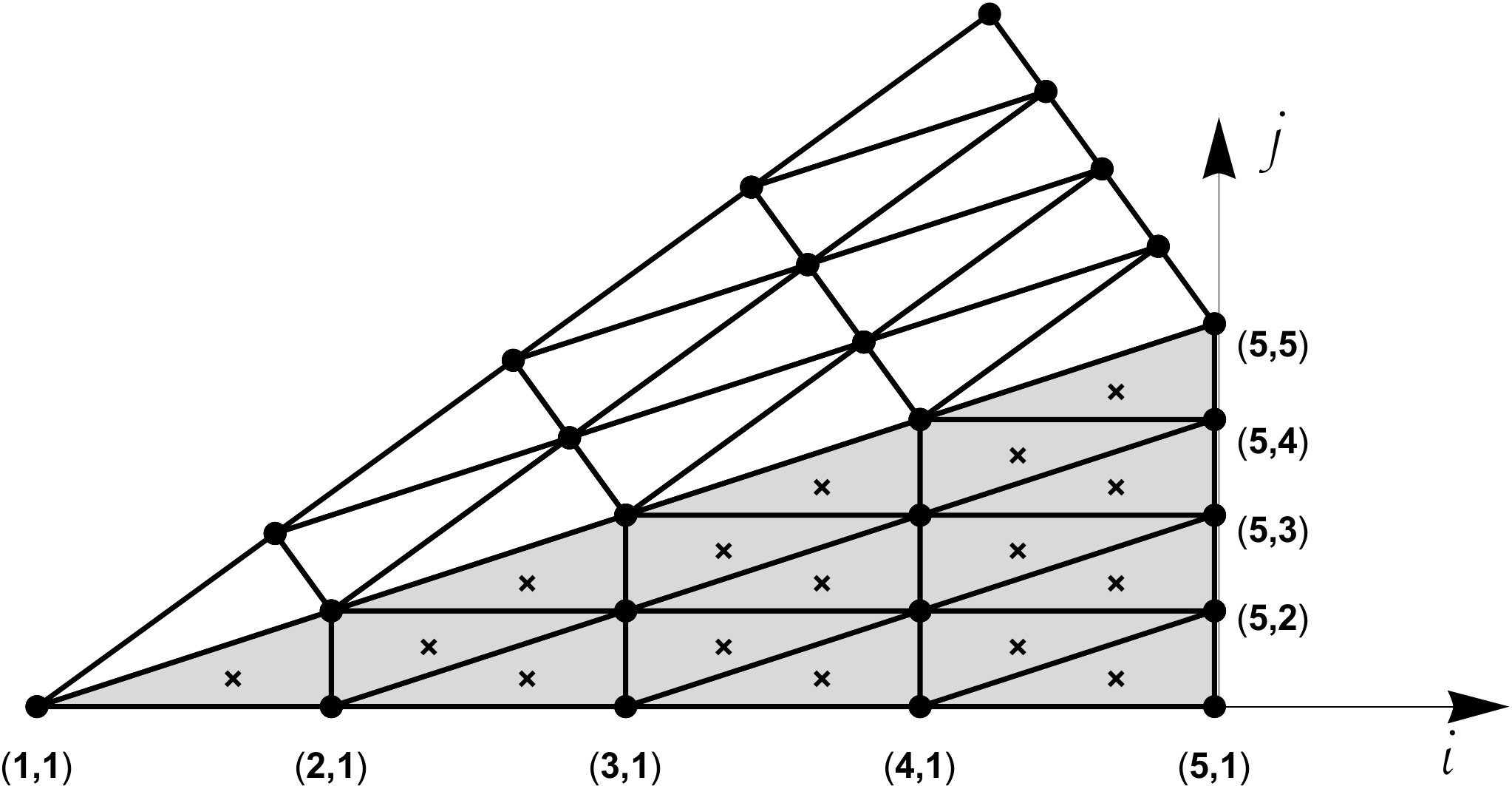}
\end{center}
\caption{\small  Discretization scheme for $P_{2n}$ with $N_c=4$ subdivisions of the apothem. The shaded region denotes the fundamental triangle $T_{2n}$. Each triangular plaquette is \emph{similar} to $T_{2n}$. Functions $\phi^\alpha$ and $\varphi^\alpha$ are defined on vertices of plaquettes and their derivatives are defined on centroids. The triangulation is extended across the hypotenuse of $T_{2n}$ by reflection. }
\label{fig:2n-gon_Disc}
\end{figure}

For any polygon $P_{2n}$ with $n \geq 4$, we discretize
the fundamental region $T_{2n}$ using the same strategy: 
we discretize  into smaller triangles \emph{similar} to $T_{2n}$. 
 $N_c$ is the number of subdivisions of the 
apothem as well as the number of subdivisions of half of the polygon edge $\tilde e_1$.
This gives $N_c^2$ plaquettes on $T_{2n}$ and $nN_c^2$ plaquettes on $Q_{2n}$. 
 Figure~\ref{fig:2n-gon_Disc} shows the case $N_c=4$, with 
$T_{2n}$ shaded, and a copy of it immediately above, whose discretization
is obtained by reflection across the hypotenuse of $T_{2n}$. 
These two triangles
 define the building block that is then successively 
 rotated rigidly to provide a triangulation of $Q_{2n}$.

\begin{figure}[!ht]
\leavevmode
\begin{center}
\epsfysize3.0cm
\epsfbox{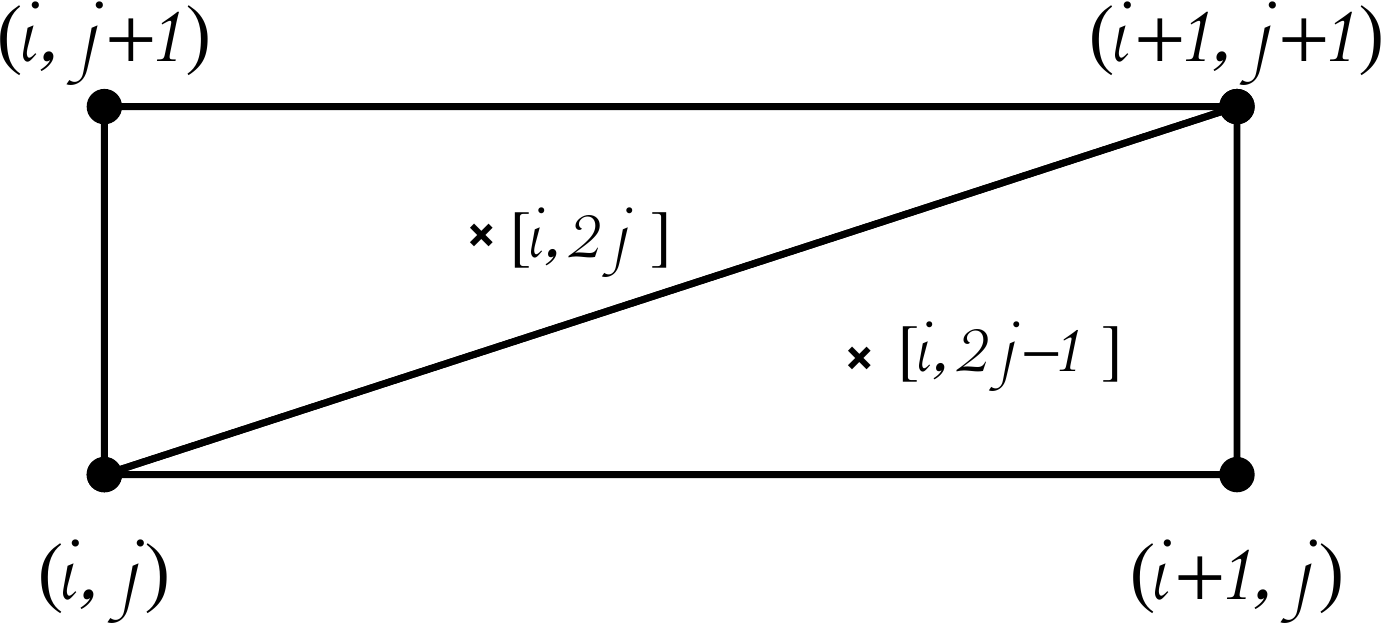}
\end{center}
\caption{\small  Labeling the vertices and centroids of triangular plaquettes. As in the case of the hexagon, the left-lower vertex $(i,j)$ is used to label centroids of two triangles. The centroid of the lower triangle  is labeled
 $\sbkt{i,2j-1}$  and the centroid of the upper triangle is labeled $[i, 2j]$. }
\label{fig:2n_Lab}
\end{figure}

The functions $\phi^\alpha$ 
 and $\varphi^\alpha$, relevant to the primal and dual programs, respectively, 
 are defined on the vertices of the triangular plaquettes. 
 These vertices, over $T_{2n}$, are labelled by integers $(i,j)$, 
 with $i=1,\cdots,N_c+1$ 
 and $j=1,\cdots , i$.   To extend the labeling over $Q_{2n}$ one can simply let
 $j$ run over a larger set of integers. 
Derivatives  are calculated at the {\em centroids} of the triangular plaquettes,
marked by  $\times$. These points are labelled by two integers $\sbkt{i,j}$, as explained in Figure~\ref{fig:2n_Lab}. These integers take values $i=1,2,\cdots N_c$ and $j=1,2,\cdots 2i-1$.
\end{appendix}


\clearpage

\begin{thebibliography}{99}

 \bibitem{strebel}
  K.~Strebel, {\em Quadratic differentials,}  Springer-Verlag Berlin Heidelberg 1984.


\bibitem{ahlfors} 
L.~V.~Ahlfors,  {\em  Conformal Invariants:  Topics in Geometric Function Theory},  AMS Chelsea Publishing (1973). 


\bibitem{Zwiebach:1990nh} 
  B.~Zwiebach,
  ``How Covariant Closed String Theory Solves A minimal-area problem,''
  Commun.\ Math.\ Phys.\  {\bf 136}, 83 (1991).
  doi:10.1007/BF02096792



\bibitem{Zwiebach:1992ie} 
  B.~Zwiebach,
  ``Closed string field theory: Quantum action and the B-V master equation,''
  Nucl.\ Phys.\ B {\bf 390}, 33 (1993)
  doi:10.1016/0550-3213(93)90388-6
  [hep-th/9206084].


\bibitem{Wolf:1992bk} 
  M.~Wolf and B.~Zwiebach,
  ``The Plumbing of minimal area surfaces,''
 Journal of Geometry and Physics {\bf 15} (1994) 23-56.
  [hep-th/9202062].
  

\bibitem{gromov} 
 M. Gromov, ``Filling Riemannian Manifolds,"  
 J.\ Differential Geom. {\bf 18} (1983) 1-147.
 
  \bibitem{m_katz}     
  M. G. Katz,  {\em Systolic Geometry and Topology,}   Mathematical Surveys and Monographs, Volume 137.  American Mathematical Society 2007.

\bibitem{mberger}
M. Berger, ``A Panoramic View of Riemannian Geometry," 
Springer-Verlag Berlin Heidelberg 2003.

\bibitem{guth}
L.~Guth, ``Metaphors in systolic geometry",  arXiv:1003.4247.
``Systolic inequalities and minimal hyper surfaces", arXiv:0903.5299.

 
 \bibitem{bavard}
 C. Bavard,  ``In\'egalit\'es isosystoliques conformes,"
Comment. Math. Helv. {\bf 67} (1992) 146-166. 

\bibitem{katz}
Mikhail Katz and St\'ephane Sabourau, ``An optimal systolic inequality for CAT(0) metrics in genus two." 
      Pacific Journal of Mathematics, vol. {\bf 227} (2006), no. 1, 95-107.  \hfill

\bibitem{katz2} 
Mikhail Katz and St\'ephane Sabourau,  ``Systolically extremal nonpositively curved 
surfaces are flat with finitely many singularities,"  to appear in 
Journal of Topology and  Analysis.  DOI: 10.1142/S1793525320500144 [arXiv:1904.00730] 


\bibitem{calabi}
E. Calabi,  ``Extremal isosystolic metrics for compact surfaces."
pp. 167-204 in {\em Actes
de la Table Ronde de G\'eom\'etrie Diff\'erentielle} (Luminy, 1992), Soc. Math. France, Paris, 1996.  




\bibitem{Headrick:2018ncs} 
  M.~Headrick and B.~Zwiebach,
  ``Convex programs for minimal-area problems,''
  arXiv:1806.00449 [hep-th].
  
  
\bibitem{Headrick:2018dlw} 
  M.~Headrick and B.~Zwiebach,
  ``Minimal-area metrics on the Swiss cross and punctured torus,''
  arXiv:1806.00450 [hep-th].
  
    \bibitem{boyd}
  S. Boyd and L. Vandenberghe, {\em Convex Optimization}, 
  Cambridge University Press (2004).   
  Available online at: https://web.stanford.edu/~boyd/cvxbook/

  
 \bibitem{bryant}
R.~L.~Bryant, ``On extremals with prescribed Lagrangian densities,"
in ``Manifolds and Geometry" (Pisa, 1993), Cambridge University Press, 1996.
arXiv:dg-ga/9406001. 


 
  \bibitem{pu}
P.~M.~Pu, ``Some inequalities in certain nonorientable Riemannian 
manifolds", Pacific. J. ~Math.~{\bf 2} (1952) 55-71.
 



    
      
  \end{thebibliography}
\end{document}